\documentclass[10pt]{amsart}
\usepackage{verbatim}
\usepackage{xspace}

\makeindex

\usepackage{syntonly}
\usepackage{amsmath,amsthm,amssymb,amscd}
\usepackage[all]{xy}
\usepackage{hyperref}

\newsavebox{\cm}
\sbox{\cm}{ \begin{picture}(12,8)
              \put(6,4){\oval(8,8)[b]}
              \put(6,4){\oval(8,8)[r]}
              \put(6,8){\vector(-1,0){2}}
              \end{picture}  }

\newcommand{\inc}{\ensuremath{\hookrightarrow}}
\newcommand{\f}{\to}
\newcommand{\ff}{\ensuremath{\longmapsto}}
\newcommand{\ov}{\overline}

\newcommand{\Z}{\ensuremath{\mathbb{Z}}}

\newcommand{\R}{\ensuremath{\mathbb{R}}}
\newcommand{\C}{\ensuremath{\mathbb{C}}}

\newcommand{\pr}[2]{\ensuremath{\langle {#1},{#2}\rangle}}
\newcommand{\norm}[1]{\ensuremath{\|#1\|}}

\newcommand{\sd}[1]{\rtimes_{#1}}

\newcommand{\ma}{\mathcal{A}}
\newcommand{\mb}{\mathcal{B}}
\newcommand{\mc}{\mathcal{C}}

\newcommand{\me}{\mathcal{E}}
\newcommand{\mf}{\mathcal{F}}

\newcommand{\mi}{\mathcal{I}}

\newcommand{\ml}{\mathcal{L}}
\newcommand{\mo}{\mathcal{O}}
\newcommand{\mr}{\mathcal{R}}
\newcommand{\mv}{\mathcal{V}}

\newcommand{\rief}{\mathsf{R}}

\renewcommand{\r}[1]{\ensuremath{\big|_{#1}}}

\newcommand{\link}[1]{\ensuremath{\mathbb{L}({#1})}}

\newcommand{\adj}[1]{\mathcal{L}({#1})}
\renewcommand{\k}[1]{\mathcal{K}({#1})}
\newcommand{\nuc}[1]{\Bbbk ({#1})}
\newcommand{\nucc}[1]{\Bbbk_c({#1})}
\newcommand{\nucr}[1]{\Bbbk_r({#1})}

\newcommand{\sM}{\ensuremath{\stackrel{{}_M}{\sim}}}

\newcommand{\bts}[3]{\ensuremath{{#1}\bigotimes_{#2}{#3}}}
\newcommand{\bats}[2]{\ensuremath{{#1}\bigodot{#2}}}

\newcommand{\bc}{\begin{center}}
\newcommand{\ec}{\end{center}}
\newcommand{\be}{\begin{enumerate}}
\newcommand{\ee}{\end{enumerate}}
\newcommand{\bi}{\begin{itemize}}
\newcommand{\ei}{\end{itemize}}
\newcommand{\bd}{\begin{description}}
\newcommand{\ed}{\end{description}}
\newcommand{\beq}{\begin{equation}}
\newcommand{\eeq}{\end{equation}}
\newcommand{\beqa}{\begin{eqnarray}}
\newcommand{\eeqa}{\end{eqnarray}}
\newcommand{\bfr}{\begin{flushright}}
\newcommand{\efr}{\end{flushright}}
\newcommand{\bfl}{\begin{flushleft}}
\newcommand{\efl}{\end{flushleft}}
\newcommand{\bp}{\begin{picture}}
\newcommand{\ep}{\end{picture}}

\theoremstyle{plain}
\newtheorem{thm}{Theorem}[section]
\newtheorem{prop}[thm]{Proposition}
\newtheorem{lem}[thm]{Lemma}
\newtheorem{cor}[thm]{Corollary}          
\newtheorem{df}[thm]{Definition}
\theoremstyle{definition}

\newtheorem{ex}[thm]{Example}

\theoremstyle{remark}
\newtheorem{rk}[thm]{Remark}

\newcommand{\al}{\ensuremath{\alpha}} 
\newcommand{\del}{\ensuremath{\delta}} 
\newcommand{\ga}{\ensuremath{\gamma}}
\newcommand{\la}{\ensuremath{\lambda}}

\newcommand{\om}{\ensuremath{\omega}}
\newcommand{\Om}{\ensuremath{\Omega}}
\newcommand{\Ga}{\ensuremath{\Gamma}}
\newcommand{\La}{\ensuremath{\Lambda}}

\newcommand{\Del}{\ensuremath{\Delta}} 
\newcommand{\ka}{\ensuremath{\kappa}}

\newcommand{\env}[1]{ \ensuremath{ {#1}^{\mathsf{e}} } }
\newcommand{\hal}{\ensuremath{\env{\al}}}
\newcommand{\nal}{\ensuremath{\check{\al}}}
\newcommand{\nbeta}{\ensuremath{\check{\beta}}}

\newcommand{\tal}{\ensuremath{\tilde{\al}}}
\newcommand{\hex}{\ensuremath{\env{X}}}

\newcommand{\hh}{\ensuremath{\env{H}}}
\newcommand{\hu}{\ensuremath{\env{u}}}
\newcommand{\hphi}{\ensuremath{\env{\phi}}}
\newcommand{\hpi}{\ensuremath{\env{\pi}}}
\newcommand{\id}{\vartriangleleft}

\newcommand{\ve}{\vspace*{.01cm}}

\newcommand{\cs}{$C^*$-algebra\xspace}
\newcommand{\css}{$C^*$-algebras\xspace}
\newcommand{\ct}{$C^*$-tring\xspace}
\newcommand{\cts}{$C^*$-trings\xspace}
\newcommand{\ea}{enveloping action\xspace}
\newcommand{\eas}{enveloping actions\xspace}
\newcommand{\fb}{Fell bundle\xspace}
\newcommand{\fbs}{Fell bundles\xspace}
\newcommand{\hm}{homomorphism\xspace}
\newcommand{\hms}{homomorphisms\xspace}
\newcommand{\ilt}{inductive limit topology\xspace}
\newcommand{\lc}{locally compact\xspace}
\newcommand{\mea}{Morita enveloping action\xspace}
\newcommand{\meas}{Morita enveloping actions\xspace}
\newcommand{\pa}{partial action\xspace}
\newcommand{\pas}{partial actions\xspace}

\newcommand{\rep}{representation\xspace}
\newcommand{\reps}{representations\xspace}
\newcommand{\rr}{regular \rep\xspace}

\DeclareMathOperator{\gen}{span}

\DeclareMathOperator{\p}{Prim}
\DeclareMathOperator{\dom}{dom}
\DeclareMathOperator{\gr}{Gr}

\AtEndDocument{\multipasswarning}
\newcommand{\multipasswarning}{%
 \clearpage
 \typeout{%
 ***************************************************************************}
 \typeout{%
 Note: This document needs to run through LaTeX three times, instead of }
 \typeout{%
 the usual two, to resolve indirect cross-references}
 \typeout{%
 ***************************************************************************}
 }

\begin{document}

\title[Enveloping Actions and Takai Duality]
      {Enveloping Actions and Takai Duality for Partial Actions}
\author[Fernando Abadie]{Fernando Abadie}

\address{Centro de Matem\'atica\\
         Facultad de Ciencias\\
         Universidad de la Rep\'ublica\\
         Igu\'a 4225\\
         CP 11400\\
         Montevideo, Uruguay}
\email{fabadie@cmat.edu.uy}
\date{July 14, 2000}
\thanks{This work was partially financied by Fapesp, Brazil, 
        Processo No. 95/04097-9}
\begin{abstract}
We show that any continuous partial action $\al$ on a topological space $X$ 
is the restriction of a suitable continuous global action $\env{\al}$, that is 
essentially unique. We call this action $\env{\al}$ the enveloping action of 
$\al$, and the space $\env{X}$ where $\env{\al}$ acts is called the 
enveloping space of $X$. $\env{X}$ is Hausdorff if and only if $X$ is 
Hausdorff and the graph of $\al$ is closed. 
\par In the case of \css, we prove that any \pa has a unique 
enveloping action up to Morita equivalence, and that the corresponding 
reduced crossed products are Morita equivalent. The study of the enveloping 
action up to Morita equivalence reveals the form that Takai duality takes for 
partial actions.
\par By applying our constructions, we prove that any \pa of a 
connected group on a unital \cs must be a global action. We also prove that 
the reduced crossed product of the reduced cross sectional algebra of a Fell 
bundle by the dual coaction is liminal, postliminal, or nuclear, if and only 
if the unit fiber of the bundle is liminal, postliminal, or nuclear, 
respectively.     
\end{abstract}
\maketitle
\tableofcontents
\section{Introduction}\label{sec:intro}
\par Partial actions on \css were gradually introduced in \cite{excirc}, 
\cite{mc} and \cite{extwist}. Since then, several classes of \css have been 
described as crossed products by \pas. This is the case of approximately 
finite, Bunce--Deddens and Cuntz--Krieger algebras, among others. 
(\cite{exaf}, \cite{exbd}, \cite{el} and \cite{elq}). In addition, the 
description of a \cs as a crossed product by a \pa has proved to be 
useful to describe its structure and sometimes to compute its 
K-theory. 
\par In the present paper we deal with \eas of \pas. That is, we discuss   
the problem of deciding whether or not a given \pa is the restriction of some
global action, and the uniqueness of this global action.
\par The exact statement of the problem depends on the category under 
consideration, as much as the definition of a \pa does. For instance, 
in the category of topological spaces and 
continuous maps, we say that an action $\beta$, acting on $Y$, is an \ea of 
the \pa $\al$, acting on $X$, if $X$ is an open subset of $Y$, 
$\al=\beta\r{X}$, and $Y$ is the $\beta$--orbit of $X$. In the category of 
\css and their \hms, we say that $(\beta,B)$ is an \ea of the \pa $\al$, 
acting on the \cs $A$, if 
$A\id B$, $\al=\beta\r{A}$, and $B$ is the closed linear $\beta$--orbit of 
$A$. In this paper we discuss \eas in both categories. 
In the first one, the \ea always exists and is unique. In 
the second one it is unique when it exists, but in general this is not 
the case. For this reason we consider a weaker notion of \ea in the context 
of \css, called \mea. We show that any \pa on a \cs has a \mea, which is 
unique up to Morita equivalence. Moreover, the corresponding reduced crossed 
products are strongly Morita equivalent. It turns out that \meas of \pas 
are intimately related to Takai duality: if $\al$ is a \pa of the group 
$G$ on a \cs $A$, $\del$ is the dual coaction of $G$ on $A\sd{\al,r}G$ and 
$\hat{\del}$ is the dual action of $G$ on $A\sd{\al,r}G\sd{\del, r}\hat{G}$, 
then $\hat{\del}$ is the \mea of $\al$.    
\par The structure of the paper is as follows. 
\par In Section \ref{sec:top} we study \pas in a topological context. In 
this case we show that for every \pa there exists a unique \ea, which is 
characterized by a universal property (\ref{thm:main1}). We exhibit an example 
where the \pa acts on a Hausdorff space while its \ea acts on a 
non--Hausdorff space (\ref{ex:z2}). This implies that in the category 
of \css, the problem of existence of \eas does not have a solution in general. 
Those \pas whose \eas act on a Hausdorff space are precisely those with 
closed graph (\ref{prop:t2}).
\par Section \ref{sec:c*} is devoted to consider the problem of \eas in the 
category of \css. In view of the results of Section \ref{sec:top}, there is 
in general no \ea for a given \pa on a \cs. So the main matter of this 
section is the uniqueness of the 
\ea. It is shown in Theorem~\ref{thm:unienv} that the \ea is unique when it 
exists. In the remainder of the section we study some relations between 
the \css where a \pa and its \ea, respectively, act.
\par In Section \ref{sec:cross} we discuss the 
relation between the reduced crossed products $A\sd{\al,r}G$ 
and $B\sd{\beta,r}G$, where $(\beta,B)$ is the enveloping action of $(\al,A)$. 
In Theorem \ref{thm:equivmor} we prove that they are Morita equivalent 
(In this paper Morita equivalence means strong Morita equivalence). 
As an application of this result, we show in \ref{prop:dil} that 
any partial \rep of a discrete amenable group may be dilated to a unitary 
\rep of $G$, so in particular it is a positive definite map.
\par In order to overcome the negative result obtained in 
Section~\ref{sec:cross} about the 
existence of \eas, we introduce in Section~\ref{sec:mea} the weaker notion 
of \mea. This concept involves Morita equivalence of \pas, which is defined 
and studied in this section. In particular, we show that the reduced crossed 
products of Morita equivalent \pas are Morita equivalent 
(\ref{prop:equivprods}). This allows us to  
deduce that the reduced crossed product by a \pa and the reduced crossed 
product by a corresponding \mea are Morita equivalent (\ref{prop:equivmor2}).
\par In the sixth section we pave the way for proving the main result of the 
paper, namely the existence and 
uniqueness of the \mea, which is achieved in Section \ref{sec:exun} 
(\ref{thm:super}, \ref{prop:super}). 
With this goal in mind, we consider a \cs, $\nuc{\mb}$, that is the completion of a certain *-algebra of integral operators naturally associated to a \fb. The algebra $\nuc{\mb}$ carries a canonical action of the group~$G$, which turns out to be a Morita enveloping action of the partial action $\alpha$ when $\mb$ is the Fell bundle associated to $\alpha$.  
In the last section we will see that, if the \fb is associated to a \pa $\al$ 
on a \cs $A$, then the algebra $\nuc{\mb}$ also agrees with the double crossed 
product $A\!\sd{\al,r}\!G\!\sd{\del,r}\!\hat{G}$, where $\del$ is the dual 
coaction on $A\!\sd{\al,r}\!G$ (\ref{prop:nadual}).
\par This paper corresponds to the second part of my doctoral thesis 
(\cite{fav0}). It is a pleasure to express my  
gratitude to my advisor Ruy Exel for his guidance and       
several conversations about \eas, which greatly enriched this work. 
\section{Enveloping actions: the topological case}\label{sec:top}
\par In this section we consider the problem of \eas in the category of 
topological spaces and continuous maps. In the first part we give the 
necessary definitions and some examples. In the second one we show that any 
\pa has a unique \ea, which is characterized by a universal property. This 
result implies, in particular, that if $\mathbf{v}$ is a vector field on a 
smooth manifold $X$, then there exist a smooth manifold $Y$, a vector 
field $\mathbf{w}$ on $Y$, and an inclusion $\iota:X\to Y$, such that 
$\iota(X)$ is open in $Y$ and $\mathbf{v}=\mathbf{w}\iota$. We show that 
if $(\al,X)$ is a \pa, where $X$ is a Hausdorff space, and if $(\beta,Y)$ 
is its \ea, then $Y$ is a Hausdorff space if and only if the graph of $\al$ 
is closed. 
\ve
\subsection{Partial actions: basic facts and examples}
\begin{df}\label{df:pa}
A partial action of a topological group $G$ on a topological 
space $X$ is a pair $\al =\big(\{X_s\}_{s\in G},\{\al_s\}_{s\in G}\big)$ such 
that:
\be
 \item $X_t$ is open in $X$, and $\al_t:X_{t^{-1}}\f X_t$ is a homeomorphism, 
       $\forall t\in G$.
 \item The set $\Ga_{\al}=\big\{ (t,x)\in G\times X:\, t\in G, 
       x\in X_{t^{-1}}\big\}$ is open in $G\times X$, and the function 
       (also called $\al$) $\al:\Ga_{\al}\f X$ given by $(t,x)\ff\al_t(x)$ 
       is continuous.
 \item $\al$ is a partial action, that is, $X_e=X$, and $\al_{st}$ is an 
       extension of $\al_s\al_t$, $\forall s,t\in G$.
\ee  
If $\al =\big(\{X_t\}_{t\in G},\{\al_t\}_{t\in G}\big)$ and $\beta
=\big(\{Y_t\}_{t\in G},\{\beta_t\}_{t\in G}\big)$ are 
\pas of $G$ on $X$ and $Y$, we say that a continuous function 
$\phi: X\f Y$ is a morphism $\phi :\al\f\beta$ if $\phi(X_t)\subseteq Y_t$, 
 and the following diagram commutes, $\forall t\in G$:     
\[ \xymatrix
{X_{t^{-1}}\ar@{->}[r]^-{\phi}\ar[d]_-{\al_t}&Y_{t^{-1}}\ar[d]^-{\beta_t}\\
X_t\ar@{->}[r]_-{\phi}&Y_t} \] 
If we forget the topological structures of $G$ and $X$, we say that $\al$ is a 
set theoretic partial action; note that in this case condition 2. is  
superfluous, and condition 1. amounts to saying that each $\al_t$ is a 
bijection.
\end{df}

Condition 3. above is equivalent to the following set of conditions 
(see Lemma 1.2 of \cite{qr}):
\be
 \item $\al_e=id_X$ and $\al_{t^{-1}}=\al_t^{-1}$, $\forall t\in G$.
 \item $\al_t(X_{t^{-1}}\cap X_s)=X_t\cap X_{ts}$, $\forall s,t\in G$.
 \item $\al_s\al_t :X_{t^{-1}}\cap X_{t^{-1}s^{-1}}\f X_s\cap X_{st}$ 
       is a bijection, and $\al_s\al_t(x)=\al_{st}(x)$, $\forall x\in 
       X_{t^{-1}}\cap X_{t^{-1}s^{-1}}$ and $\forall s,t\in G$.
\ee
\begin{ex}\label{ex:rest}
Let $\beta:G\times Y\f Y$ be a continuous global action
and let $X$ be an open subset of $Y$. Consider $\al=\beta\r{X}$, the 
``restriction'' of $\beta$ to $X$, that is: $X_t=X\cap\beta_t(X)$, and 
$\al_t:X_{t^{-1}}\f X_t$ such that $\al_t(x)=\beta_t(x)$, $\forall t\in G$, 
$x\in X_{t^{-1}}$. It is easy to verify that $\al$ is a \pa on $X$. In fact, 
the main result of this section shows that any \pa arises in this way. Note
that, in particular, $\beta$ may be identified with the \pa $\beta\r{Y}$.
\end{ex}
\begin{ex}\label{ex:flow}
The flow of a differentiable vector field is a partial action. More precisely, 
consider a smooth vector field $\mathbf{v}:X\f TX$ on a manifold $X$, and for 
$x\in X$ let $\ga_x$ be the corresponding integral curve through $x$ (i.e.: 
$\ga_x(0)=x$), defined on its maximal interval $(a_x,b_x)$. Let us define, 
for $t\in\R$: $X_{-t}=\big\{ x\in X:\, t\in (a_x,b_x)\big\}$, 
$\al_t:X_{-t}\f X_t$ such that $\al_t(x)=\ga_x(t)$, and 
$\al =\big(\{X_t\}_{t\in\R},\{\al_t\}_{t\in\R}\big)$. Then $\al$ is a partial 
action of $\R$ on $X$.
\end{ex}
\par It is well known that the integral curves of a vector field on a 
compact manifold $X$ are defined on all of $\R$. This is a particular case 
of the next result, which in turn may be generalized to a theorem about 
partial actions on \css to be proved later in Section \ref{sec:apps} 
(Corollary \ref{cor:conti2}) 

\begin{prop}\label{prop:conti}
Let $\al$ be a \pa\ of $G$ on a compact space $X$. Then there exists
an open subgroup $H$ of $G$ such that $\al$ restricted to $H$ is a  
global action. In particular, if $G$ is connected, $\al$ is a global action.
\end{prop}
\begin{proof}
Let $A_x=\{ t\in G:\, x\in X_{t^{-1}}\}$, and $A=\cap_{x\in X}A_x$. It is 
clear that $e\in A$ and $st\in A$ whenever $s$, $t\in A$; that is, $A$ is a
submonoid of $G$. For every $x\in X$ there exist open neighborhoods 
$U_x\subseteq X$ of $x$ and $V_x\subseteq G$ of $e$ such that 
$V_x\times U_x\subseteq \Ga_{\al}$, and $V_x=V_x^{-1}$. Since $X$ is compact, 
there exist $x_1,\ldots ,x_n\in X$ such that $X=\cup_{j=1}^n U_{x_j}$. 
Consider now the neighborhood $V=\cap_{j=1}^nV_{x_j}$. Since $V$ is  
symmetric and $V\subseteq A$, we have that $H=\cup_{n=1}^\infty V^n$ is an 
open subgroup of $G$ contained in $A$. 
\par As for the last assertion, just recall that the unique open subgroup of 
a connected group is the group itself.   
\end{proof}
\ve
\subsection{Existence and uniqueness of enveloping actions}
\begin{thm}\label{thm:main1}
Let $\al$ be a \pa of $G$ on $X$. Then there exists a pair 
$(\iota,\env{\al})$ such that $\env{\al}$ is a continuous action of $G$ on 
a topological space $\env{X}$, and $\iota :\al\f\env{\al}$ is a morphism, 
such that for any morphism $\psi:\al\f\beta$, where $\beta$ is a continuous 
action of $G$, there exists a unique morphism $\env{\psi}:\env{\al}\f\beta$ 
making the following diagram commutative:  
 \[ \xymatrix{
       \al\ar[rr]^{\iota}
\ar[dr]_{\psi}&\ar @{}[d]
|{\circlearrowleft}&\hal\ar@{-->}[dl]^{\env{\psi}}\\ &\beta&}\]
\par Moreover, the pair $(\iota,\env{\al})$ is unique up to canonical 
isomorphisms, and:
\be
 \item $\iota(X)$ is open in $\env{X}$.
 \item $\iota:X\f\iota(X)$ is a homeomorphism.
 \item $\env{X}$ is the $\env{\al}$--orbit of $\iota(X)$.
\ee
\end{thm}
\begin{proof}
Let us consider the action $\ga:G\times (G\times X)\f G\times X$ such that 
$\ga_s(t,x)=(st,x)$, $\forall s,t\in G$, $x\in X$. We endow $G\times X$ with 
the product topology, so $\ga$ is a continuous action. Moreover, $\ga$ is 
compatible with respect to the equivalence relation $\sim$ on $G\times X$ 
given by: $(r,x)\sim (s,y)\iff x\in X_{r^{-1}s}$ and $\al_{s^{-1}r}(x)=y$. 
Thus $\ga$ induces an action $\env{\al}$ of $G$ by homeomorphisms of the quotient  
topological space $\env{X}=(G\times X)/\sim$. Let $q:G\times X\f\env{X}$ be 
the quotient map, and define $\iota:X\f\env{X}$ such that $\iota(x)=q(e,x)$. 
Since the inclusion $X\inc G\times X$ given by $x\ff (e,x)$ is continuous, 
we have that $\iota$ also is. Note that $\iota$ is 
clearly injective, so to prove (1) and (2) it suffices to show that it is an open map. Let $U\subseteq X$ be an open subset. We have 
to show that $q^{-1}\big(\iota(U)\big)$ is open in $G\times X$. 
But $q^{-1}\big(\iota(U)\big)
=\{ (t,x):\, (t,x)\sim (e,y)\text{ for some }y\in U \} 
=\{ (t,x):\, \al_t(x)\in U\} =\al^{-1}(U)$, which is open in $\Ga_{\al}$ 
because $\al$ is continuous, and hence open in $G\times X$ because $\Ga_{\al}$ 
is open in $G\times X$. Statement (3) is clear, because $q(t,x)=\env{\al}_t\big(\iota(x)\big)$.
\par Let us see that the action $\alpha^e$ is continuous and $\iota:\alpha\to\alpha^e$ is a morphism. Note first that $q$ is an open map, for if $U\subseteq G$ and $V\subseteq X$ are open subsets, we have 
$q(U\times V)
  =\cup_{t\in U}q(\gamma_t(\{e\}\times V))
  =\cup_{t\in U}\alpha_t^e(\iota (V))$,
which is open because $\iota$ is open and every $\alpha_t^e$ is a homeomorphism.  Now let $W\subseteq X^e$ be an open subset. Since $\alpha^e\circ(id\times q)=q\circ\gamma:G\times(G\times X)\to X^e$ is continuous, we have that $(id\times q)^{-1}((\alpha^e)^{-1}(W))=\gamma^{-1}(q^{-1}(W))$ is an open subset of $G\times(G\times X)$. Therefore $(\alpha^e)^{-1}(W)$ is open in $G\times X^e$, because $id\times q$ is an open and surjective map. We see that $\iota$ is a morphism: if $x\in X_{t^{-1}}$: 
\[\iota\big(\al_t(x)\big)=q\big(e,\al_t(x)\big)=q(t,x)
                         =q\big(\ga_t(e,x)\big)
                         =\env{\al}_t\big(q(e,x)\big)
                         =\env{\al}_t\big(\iota(x)\big),\]
\par We must prove now that the pair $(\iota,\env{\al})$ has the claimed 
universal property. Note first that if $\beta:G\times Y\f Y$ is a continuous 
action and $\psi:X\f Y$ is any continuous function, then the map 
$\psi':G\times X\f Y$ such that $\psi'(t,x)=\beta_t\big(\psi(x)\big)$ is a 
morphism $\ga\f\beta$. Moreover, if $\psi:\al\f\beta$ is also a morphism, 
then $\psi'$ is compatible with $\sim$: if $(r,x)\sim (s,y)$ in $G\times X$, 
since $\al_{s^{-1}r}(x)=y$, we have:
\[ \beta_{s^{-1}}\big(\psi'(r,x)\big)
                    =\beta_{s^{-1}}\big(\beta_r(\psi (x))\big)
                    =\beta_{s^{-1}r}\big(\psi(x)\big)
                    =\psi\big(\al_{s^{-1}r}(x)\big)
                    =\psi(y), \]
and therefore: $\psi'(r,x)=\beta_s\big(\psi (y)\big) =\psi'(s,y)$. Thus 
$\psi'$ induces a continuous map $\env{\psi}:\env{X}\f Y$, such that 
$\env{\psi}\big(q(t,x)\big)=\beta_t(\psi(x))$, $\forall t\in G$, $x\in X$.
We have that $\env{\psi}\iota(x)=\env{\psi}\big(q(e,x)\big)=\psi (x)$, and 
it is also clear that $\env{\psi}:\env{\al}\f\beta$ is a morphism, uniquely 
determined by the relation $\env{\psi}\iota=\psi$. 
\par Since the pair $(\iota,\env{\al})$ is characterized by a universal 
property, it is unique up to isomorphisms (In categorical terms, 
$\env{\al}$ is a universal from $\al$ to $\mathfrak{F}$, where 
$\mathfrak{F}:\ma\to\mathcal{P}\ma$ is the forgetful functor from the category 
of actions to the category of \pas; see \cite{rowen} for details).  
\end{proof}
\begin{df}\label{df:envaction}
Let $\al$ be a \pa of $G$ on $X$. We say that the action $\env{\al}$ provided 
by Theorem \ref{thm:main1} is an enveloping action of $\al$. We will also say 
that $\env{X}$ is the enveloping space of $X$, $\env{\psi}$ is the enveloping 
morphism of $\psi$, etc.  
\end{df}
\begin{rk}\label{rk:susp}
Assume that $h:X\f X$ is a homeomorphism, so we have an action of $\Z$ on 
$X$. We may think of this action as a partial action of $\R_d$ on $X$, where 
$\R_d$ denotes the real numbers with the discrete topology. Indeed, define 
$X_s=X$ if $s\in \Z$, $X_s=\emptyset$ if $s\notin\Z$, and 
$\al_s:X_{-s}\f X_s$ as $\al_s=h^s$ if $s\in\Z$, $\al_s=\emptyset$ otherwise. 
Note that $\al$ is not a partial action of $\R$ on $X$, because $\Z\times X$ 
is not open in $\R\times X$. However, we can imitate the construction of the 
enveloping action made in the proof of \ref{thm:main1} above, using $\R$ 
instead of $\R_d$, to obtain 
 a global continuous action 
$\beta:\R\times(\R\times X)/\sim\f (\R\times X)/\sim$, such that 
$\beta_n(x)=\al_n(x)$, $\forall n\in\Z$, $x\in X$. This action $\beta$ 
is called the {\em suspension of $h$}, and its construction is well known in 
dynamical systems theory (see \cite{tomi}, page 45).
\end{rk}
\par From now on we will suppose, as we can, that $X\subseteq\env{X}$. Since 
$\env{X}$ is the $\hal$--orbit of $X$, we see that $X$ and $\hex$ share the 
same local properties. However, their global properties may be very different, 
as shown in the next two examples. 
\begin{ex}\label{ex:pi1}
Consider the action $\beta:\Z\times S^1\f S^1$ given by the 
rotation by an irrational angle $\theta$: 
$\beta_k(z)=e^{2\pi ik\theta}z$, $\forall k\in\Z$, $z\in S^1$.
Let $U$ be a nonempty open arc of $S^1$, $U\neq S^1$, and consider the \pa 
given by the restriction $\al$ of $\beta$ to $U$ (see Example \ref{ex:rest}). 
Since the action $\beta$ is minimal, it follows that $\beta$ is the 
enveloping action of $\al$. This example shows that, even when $X$ and 
$\hex$ are similar locally, their global properties may be deeply different.
In this case, for instance, the first homotopy groups of $U$ and $S^1$ are 
different.
\end{ex}
\begin{ex}\label{ex:z2}
Consider the \pa $\al$ of $\Z_2$ on the unit interval $X=[0,1]$, given by  
$\al_1=id_X$, $\al_{-1}=id_V$, where $V=(a,1]$, $a>0$. Let 
$\hal :G\times\hex\f\hex$ be the enveloping action of $\al$. 
Consider $J=J^-\cup J^+\subseteq\R^2$ with the relative topology, where 
$J^\pm=\{\pm 1\}\times [0,1]$. It is not difficult to see that 
$\hex$ is the topological quotient space obtained from $J$ by identifying 
the points $(1,t)$ and 
$(-1,t)$, $\forall t\in (a,1]$. Therefore, $\hex$ is not a Hausdorff space: 
$(1,a)$ and $(-1,a)$ do not have disjoint neighborhoods. Note also that 
$\hal_{-1}$ permutes $(1,t)$ and $(-1,t)$ for $t\in [0,a]$, and is the 
identity in the rest of $\hex$.
\end{ex}
\par The obstruction for the enveloping space to be Hausdorff
is made clear in the next proposition.
\begin{prop}\label{prop:t2}
Let $\al$ be a \pa\ of $G$ on the Hausdorff space $X$. Let $\gr(\al )$
be the graph of $\al$, that is 
$\gr(\al )=\{(t,x,y)\in G\times X\times X:\ x\in X_{t^{-1}},\, \al_t(x)=y \}.$ 
Then $\hex$ is a Hausdorff space if and only if $\gr(\al )$ is a closed 
subset of $G\times X\times X$.
\end{prop}
\begin{proof}
Let us suppose that $\hex$ is a Hausdorff space, and let  
$\big(t_i,x_i,\al (t_i,x_i)\big)\f (t,x,y)\in G\times X\times X$. 
In particular, $\al (t_i,x_i)\f y\in X$. Since $\hal$ is continuous, 
$\hal (t_i,x_i)\f\hal (t,x)$, and it must be $y=\al (t,x)$ because of the  
uniqueness of limits in Hausdorff spaces. 
\par Conversely, assume that $\text{Gr}(\al )$ is closed in  
$G\times X\times X$, and let $\env{x},\env{y}\in\hex$. We want to show 
that if there does not exist disjoint open sets in $\hex$, each of them 
containing $\env{x}$ or $\env{y}$, then $\env{x}=\env{y}$. 
Since $\hal_t$ is a  homeomorphism of $\hex$, by 3. of \ref{thm:main1} 
we may suppose that $\env{x}=x\in X$. Let $y\in X$, $t\in G$ be such that 
$\hal_t(y)=\env{y}$. If every neighborhood of $x$ intersects every 
 neighborhood of $\env{y}$, then for any pair $(U,V)$ of neighborhoods in $X$ 
of $x$ and $y$ respectively, there exists 
$x_{_{U,V}}\in U\cap\hal_t(V)$, say $x_{_{U,V}}=\hal_t(y_{_{U,V}})$, 
with $y_{_{U,V}}\in V$. Consider the net
$\{ (t,y_{_{U,V}},x_{_{U,V}})\}_{U,V}\subseteq\text{Gr}(\al)$: it converges 
to $(t,y,x)$, so $\al_t(y)=x$, because $\gr(\al)$ is closed. Hence $x=\env{y}$,
and $\hex$ is Hausdorff.  
\end{proof}
\begin{rk}\label{rk:t2} 
if $G$ is a discrete group, then $\gr(\al )$ is closed in 
$G\times X\times X$ if and only if  
$\gr(\al_t)$ is closed in $X\times X$, $\forall t\in G$.
\end{rk}
\begin{rk}\label{rk:flow}
As already seen in \ref{ex:flow}, the flow of a smooth vector field on a 
manifold is a \pa, indeed a smooth \pa. The enveloping space inherits a 
natural manifold structure, although not always Hausdorff, by translating 
the structure of the original manifold through the enveloping action. It 
would be interesting to characterize those vector fields whose flows have 
closed graphs. For such a vector field, one obtains a Hausdorff manifold that 
contains the original one as an open submanifold, and a vector field whose 
restriction to this submanifold is the original vector field. Note, however,  
that the inclusion of the original manifold in its enveloping one could be 
a bit complicated. 
\par It is possible to exhibit examples of flows with closed graphs and flows 
with non-closed graphs.
\end{rk}
\ve
\subsection{On the dynamical properties of the enveloping action}
Before closing this section we would like to make some brief 
comments about the dynamical behavior of the enveloping action.
\par Many of the algebraic and even dynamical notions related to global 
actions may be easily extended to the context of \pas. For instance, it is 
possible to make sense of expressions such as transitive \pas or minimal \pas.
To give an example, we say that a \pa $\al$ on a topological space $X$ is 
minimal when each $\al$--orbit is dense in $X$, that is, when 
$X=\ov{\{\al_t(x):\, t\in X_{t^{-1}}\}}$, $\forall x\in X$. It is not
difficult to show that the dynamical properties of $\al$ and $\env{\al}$ are 
in general the same, although we will not do it in this work. For instance, 
it is not hard to see that $\al$ is minimal if and only if $\env{\al}$ is 
minimal.  

\section{Enveloping actions: the $C^*$--case}\label{sec:c*}
\par In this section we consider partial actions on \css. We begin by 
recalling the definition of a \pa in this context, and then we introduce a  
notion of enveloping action that corresponds, in the case of commutative 
\css, to the concept of enveloping action treated in Section \ref{sec:top}. 
Next, we discuss the existence and uniqueness of enveloping 
actions, and we close the section by studying some properties of the 
enveloping \css. Throughout the rest of this paper, $G$ will denote a \lc 
Hausdorff group. 
\par The most general definition of a \pa is the one given in \cite{extwist}, 
where the reader is referred to for more information. We recall it in 
\ref{df:c*pa} below. 
\begin{df}\label{df:cf}
Let $E$ be a Banach space, $X$ a topological space, 
and for each $x\in X$, let $E_x$ be a Banach subspace of $E$. We say that  
$\{ E_x\}_{x\in X}$ is a continuous family if for any open subset $U$ of $E$, 
the set $\{ x\in X:\ U\cap E_x\neq\emptyset\}$ is open in $X$.  
\end{df}
\par If $\me =\{ (x,v)\in X\times E:\ v\in E_x\}$, with the product topology, 
and $\pi :\me\f X$ is the natural projection, then $\{E_x\}_{x\in X}$ is a 
continuous family if and only if $\pi$ is open; in this case, $(\me,\pi)$ is a 
Banach bundle (see \cite{extwist}, where the notion of continuous family was 
introduced). For details about Banach bundles and \fbs, we refer the reader to 
\cite{fd}. Note that our notation differs sometimes from that of \cite{fd}.
Also observe that \fbs are called $C^*$-algebraic bundles in that book. 
\begin{df}\label{df:c*pa}
Let $G$ be a locally compact group and 
$\al =(\{D_t\}_{t\in G},\{\al_t\}_{t\in G})$ a set theoretic \pa of $G$ on 
the \cs $A$, where each $D_t$ is an ideal of $A$ and each $\al_t$ is 
an isomorphism of \css. Consider 
$ \mb^{-1}=\{ (t,b)\in G\times B:\, b\in D_{t^{-1}}\}\subseteq G\times A$ with 
the product topology. We say that $\al$ is a \pa of $G$ on $A$ if 
$\{D_t\}_{t\in G}$ is a continuous family and the map (also called) 
$\al :\mb^{-1}\f A$ such that $(t,b)\ff \al_t(b)$ is continuous  
(note that, being $\{D_t\}_{t\in G}$ a continuous family, 
$\mb^{-1}$ is a Banach bundle). 
If $\al'=(\{D_t'\}_{t\in G},\{\al_t'\}_{t\in G})$ is a \pa of $G$ on $A'$, 
a morphism $\phi:\al\to\al'$ is a \hm $\phi:A\to A'$ such that 
$\phi(A_t)\subseteq A_t'$, $\forall t\in G$.
\end{df} 
\par In \cite{extwist}, Exel has shown that if 
$\al =(\{D_t\}_{t\in G},\{\al_t\}_{t\in G})$ is a \pa of $G$ on the 
\cs $A$, the Banach bundle $\mb_{\al}=\{ (t,x):\ 
x\in D_t\}\subseteq G\times A$, with the relative topology, is a 
Fell bundle with product and involution given by 
(here $x_t\delta_t:=(t,x_t)\in\mb_{\al}$):
$(x_t\delta_t)*(x_s\delta_s)=
           \al_t\left(\al_t^{-1}(x_t)x_s\right)\delta_{ts}$ 
and $(x_t\delta_t)^*=\al_t^{-1}(x_t^*)\delta_{t^{-1}}$
respectively. The bundle $\mb_\al$ is called the {\em semidirect product} of $A$ and $G$. 
We will also say that $\mb_\al$ is the \fb \textit{associated} with $\al$. 
The cross sectional \cs $C^*(\mb_{\al})$ of $\mb_{\al}$ is called 
\textit{crossed product} of $A$ by $\al$, and is denoted by $A\sd{\al}G$.
\begin{ex}\label{ex:rest*}
If $\beta:G\times B\f B$ is a continuous action and $A\id B$, then  
the restriction $\beta\r{A}$ of $\beta$ to $A$ (see \ref{ex:rest}) is a \pa 
of $G$ on $A$. In particular, $\beta$ may be identified with the \pa 
$\beta\r{B}$. 
\end{ex}
\par Let us concentrate for a moment on the case where $A=C_0(X)$, for some 
locally compact Hausdorff space $X$. Suppose that $\{X_t\}_{t\in G}$ is a 
family of open subsets of $X$, so $\{D_t\}_{t\in G}$ is a family of ideals 
in $A$, where $D_t=C_0(X_t)$, $\forall t\in G$. 
If $G$ is a locally compact Hausdorff space, one 
can show that $\{D_t\}_{t\in G}$ is a continuous family if and only if the 
set $\Ga=\{ (t,x):\, x\in X_t\}\subseteq G\times X$ is open with the product 
topology. 
Suppose in addition that $G$ is a group. To give an isomorphism 
$\al_t:D_{t^{-1}}\f D_t$ is equivalent to give a homeomorphism 
$\hat{\al}_t:X_{t^{-1}}\f X_t$. Now, it is possible to show that a given  
family of isomorphisms $\{\al_t:D_{t^{-1}}\f D_t\}_{t\in G}$ is a \pa on $A$ 
if and only if the corresponding family of homeomorphisms 
$\{\hat{\al}_t:X_{t^{-1}}\f X_t\}_{t\in G}$ is a \pa on $X$ (\cite{fav0}, 
\cite{fav4}).   
\par In the situation above, if the \pa 
$\hat{\al}=(\{X_t\}_{t\in G},\{\hat{\al}_t\}_{t\in G})$ has an enveloping 
action $\hat{\beta}=\env{\hat{\al}}$ acting on the enveloping space $Y$, 
then $A$ is an ideal of $B=C_0(Y)$, and the action $\beta$ induced by 
$\hat{\beta}$ on $B$ satisfies: $\beta\r{A}=\al$. Moreover, the 
$\beta$--linear orbit $[\beta(A)]:=\gen\{\beta_t(a):\, a\in A,t\in G\}$ of 
$A$ is dense in $B$, by the Stone--Weierstrass theorem. These facts justify 
the following definition.
\begin{df}\label{df:env*}
Let $\al =(\{D_t\}_{t\in G},\{\al_t\}_{t\in G})$ be a \pa of $G$ on the \cs 
$A$, and let $\beta$ be a continuous action of $G$ on a \cs $B$ that 
contains $A$. We say that $(\beta ,B)$ is an enveloping action of $(\al ,A)$ 
(in the category of \css ), if the following three properties are fulfilled:
\be
\item $A$ is an ideal of $B$ (two--sided and closed, of course).
\item $\al =\beta\r{A}$, that is $D_t=A\cap\beta_t(A)$, and 
      $\al_t(x)=\beta_t(x)$, $\forall t\in G$ and $x\in D_{t^{-1}}$.
\item $B=\ov{[\beta(A)]}$, where 
      $[\beta(A)]:=\gen\big\{\beta_t(x):\, t\in G, x\in A\big\}$.
\ee  
We then say that $B$ is an enveloping \cs of $A$. 
\end{df} 
\begin{prop}\label{prop:envab}
Let $\al =(\{D_t\}_{t\in G},\{\al_t\}_{t\in G})$ be a \pa of $G$ on a 
commutative \cs $A$, and let $\hat{\al}$ the corresponding \pa of $G$ on 
$\hat{A}$. Then the following assertions are equivalent:
\be
\item $\gr(\hat{\al} )$ is closed in  $G\times\hat{A}\times\hat{A}$.
\item $\al$ has an enveloping action in the category of commutative \css.
\item $\al$ has an enveloping action in the category of \css.
\ee 
\end{prop}
\begin{proof}
It is clear that 2. implies 3., and 1. and 2. are equivalent by Proposition  
\ref{prop:t2}, because of the categorical equivalence between 
locally compact Hausdorff spaces and commutative \css. To see that
3. implies 2, let us suppose that $(\beta ,B)$ is an enveloping action 
of the \pa $\al$ on a commutative \cs $A$. Since $A$ is a commutative 
ideal of $B$, $A$ is contained in the center $Z(B)$ of $B$: if $a\in A$, 
$b\in B$, and $(e_i)_{i\in I}\subseteq A$ is an approximate unit of $A$:
\[ ab=\lim_i(ab)e_i=\lim_ia(be_i)=\lim_i(be_i)a=\lim_ib(e_ia)=ba.\]
Since $Z(B)$ is invariant, it follows that $\beta_t(A)\subseteq Z(B)$, 
$\forall t\in G$, thus $\gen\{\beta_t(a):\, t\in G, a\in A\}
\subseteq Z(B)$. But then 
$B=\ov{\gen}\{\beta_t(a):\, t\in G, a\in A\}\subseteq Z(B)\subseteq B$, and 
therefore $B=Z(B)$, so $B$ is commutative.     
\end{proof}
\par In \cite{elq} several algebras generated by isometries satisfying certain 
relations are studied, and they are shown to be crossed products by \pas. 
At the topological level, all these \pas have \eas acting on Hausdorff spaces, 
and therefore they have also \eas at the \cs level by \ref{prop:envab}.
\ve
\subsection{On the uniqueness of enveloping actions}
\par Proposition \ref{prop:envab} and Example \ref{ex:z2} show that not 
every \pa on a \cs has an enveloping action. We will prove that at least 
the enveloping action is unique when it does exist. Note that we already 
know this in the commutative case.
 
\begin{lem}\label{lem:unienv}
Let $\{J_{\la}\}_{\la\in\La}$ be a family of ideals of a \cs $A$, and consider 
$\norm{\cdot}_{\La}:A\f\R$ such that 
$\norm{a}_{\La}=\sup_{\la\in\La}\{\norm{ax}:\, x\in J_\la,\norm{x}\leq 1\}$. 
Then $\norm{\cdot}_{\La}$ is a $C^*$-seminorm on $A$, such that 
$\norm{\cdot}_{\La}\leq\norm{\cdot}$, and $\norm{\cdot}_{\La}$ is a norm iff
$\ov{\gen}\{ x\in J_{\la}:\la\in\La\}$ is an essential ideal of $A$. In this 
case, $\norm{\cdot}_{\La}=\norm{\cdot}$.
\end{lem}
\begin{proof}
Let $B=\prod_{\la\in\La}A_{\la}$, where $A_{\la}= A$, 
$\forall\la\in\La$, and 
consider $J\!\!=\{x= (x_{\la})\in B:\, x_{\la}\in J_{\la},
\forall\la\in\La\}$. 
Then $B$ is a \cs with $\norm{b}=\sup_{\la\in\La}\norm{b_{\la}}$, and $J$ is 
an ideal of $B$. In particular, $J$ may be considered as a right Hilbert 
$B$-module with the inner product: $\pr{x}{y}=x^*y$, 
so there is a \hm $\eta:B\f\adj{J}$ given by 
$\eta(b)x=bx$. On the other hand, we have an inclusion $\iota:A\inc B$ given 
by $\iota(a)_{\la}=a$, $\forall \la\in\La$. Thus, we get a \hm 
$\tilde{\eta}=\eta\iota$, and therefore $\norm{\tilde{\eta}(a)}\leq\norm{a}$, 
$\forall a\in A$. But $\norm{\tilde{\eta}(a)}=\sup\{\norm{\eta(a)x}:\, x\in J, 
\norm{x}\leq 1\}=\norm{a}_{\La}$. Finally, it is clear that $\tilde{\eta}$ is 
injective if and only if $\ov{\gen}\{ x\in J_{\la}:\la\in\La\}$ is an 
essential ideal of $A$. 
\end{proof}

\begin{lem}\label{lem:betaga}
Let $\al =(\{D_t\}_{t\in G},\{\al_t\}_{t\in G})$ be a \pa of $G$, and 
assume that $(\beta,B)$ and $(\ga,C)$ 
are enveloping actions of $\al$. Then, $\forall a,b\in A$, $t\in G$, 
we have: $\beta_t(a)b=\ga_t(a)b$ (note that both of these products belong to 
$D_t$).  
\end{lem}
\begin{proof}
Let $(u_i)$ be an approximate unit of $D_{t^{-1}}$. Then  
$u_ia\in D_{t^{-1}}$, $\forall i$, and since $\al$ is both the restriction 
of $\beta$ and $\ga$ to $A$, we have: 
$\beta_t(a)b
=\lim\beta_t(u_ia)b
=\lim\al_t(u_ia)b
=\lim\ga_t(u_ia)b
=\ga_t(a)b.$
\end{proof}

\begin{thm}\label{thm:unienv}
Let $(\al,A)$ be a \pa of $G$, and assume that $(\beta,B)$ and $(\ga,C)$ 
are enveloping actions of $\al$. Then there exists a unique isomorphism  
$\phi:B\f C$ such that $\phi\beta_t=\ga_t\phi$, $\forall t\in G$, and  
$\phi\r{A}=id_A$.
\end{thm}
\begin{proof}
For $s\in G$, let $\norm{\cdot}_s:B\f\R$ and $\norm{\cdot}^s:C\f\R$ 
be given by $\norm{b}_s:=\sup\{\norm{bx}:\, x\in\beta_s(A),\, \norm{x}\leq 1\}$
and $\norm{c}^s:=\sup\{\norm{cy}:\, y\in\ga_s(A),\, \norm{y}\leq 1\}$. By  
Lemma \ref{lem:unienv}, $\norm{\cdot}_s$ and $\norm{\cdot}^s$ are  
$C^*$-seminorms, and $\norm{\cdot}_B=\sup_s\norm{\cdot}_s$, $\norm{\cdot}_C=
\sup_s\norm{\cdot}^s$.
\par Let $t_1,\dots ,t_n\in G$ and $a_1,\ldots ,a_n\in A$. We want to show 
that $\norm{\sum_i\beta_{t_i}(a_i)}_B=\norm{\sum_i\ga_{t_i}(a_i)}_C$. 
For this, it is enough to prove that $\norm{\sum_i\beta_{t_i}(a_i)}_s
=\norm{\sum_i\ga_{t_i}(a_i)}^s, \forall s\in G$. Let $s\in G$ and $a\in A$. 
By Lemma \ref{lem:betaga} we have:
\[\norm{\sum_i\beta_{t_i}(a_i)\beta_s(a)}
=\norm{\beta_s\big(\sum_i\beta_{s^{-1}t_i}(a_i)a\big)}
=\norm{\ga_s\big(\sum_i\ga_{s^{-1}t_i}(a_i)a\big)}
=\norm{\sum_i\ga_{t_i}(a_i)\ga_s(a)}. \]
It follows that $\norm{\sum_i\beta_{t_i}(a_i)}_s
=\norm{\sum_i\ga_{t_i}(a_i)}^s$, $\forall s\in G$, and hence that  
$\norm{\sum_i\beta_{t_i}(a_i)}_B=\norm{\sum_i\ga_{t_i}(a_i)}_C$.   
Thus, $\phi:[\beta(A)]\f [\ga(A)]$ such that  
$\phi\big(\sum_i\beta_{t_i}(a_i)\big)=\sum_i\ga_{t_i}(a_i)$ is an 
isometry of a *-dense ideal of $B$ onto a *-dense ideal of $C$, and 
therefore it extends uniquely to an isomorphism 
$\phi:B\f C$, which clearly satisfies $\ga_t\phi=\phi\beta_t$, 
$\forall t\in G$, and $\phi\r{A}=id_A$. Moreover, it is clear that these 
conditions determine $\phi$.
\end{proof}
\ve
\subsection{Some properties of the enveloping algebra}
\par To close this section we study some properties that are 
shared by a \cs and its enveloping algebra. In what follows it will be 
assumed that $(\env{\al},\env{A})$ is an enveloping action of $(\al,A)$.

\begin{prop}\label{prop:class}
Let $\mc$ be a class of \css that is closed by ideals, isomorphisms, 
and such that 
any \cs $B$ has a largest ideal $\mc(B)$ that belongs to $\mc$
($\mc$ may be, for instance, one of the following classes of \css: 
nuclear, type $I_0$, liminal, postliminal, antiliminal). Then $A\in\mc\iff\env{A}\in\mc$, and $\mc(A)=0\iff\mc(\env{A})=0$. 
\end{prop}
\begin{proof}
Note that $\mc(\env{A})$ is $\env{\al}$--invariant, and since 
$\mc(\env{A})\cap A=\mc(A)$, we have that  
$\mc(\env{A})=\ov{\gen}\{\env{\al}_t\big(\mc(A)\!:\, t\in G\big)\}$. 
From this, the result follows immediately.
\end{proof}
\begin{prop}\label{prop:envaf}
Any separable \cs has a largest ideal that is approximately finite.
If $G$ is a separable group, then $A$ is approximately finite iff $\env{A}$ is 
approximately finite.
\end{prop}
\begin{proof}
Let $AF(B)$ be the set of AF--ideals of a \cs $B$. Since $0\in AF(B)$, 
we have that 
$AF(B)\neq\emptyset$. Let $I,J\in AF(B)$, and consider the following 
exact sequence of \css:
 \[ \xymatrix@1
   {0\ar[r]&I\ar[r]^-{\iota}
           &I+J\ar[r]^-{\pi}
           &J/(I\cap J)\ar[r]&0},\]
where $\iota$ is the inclusion and $\pi$ is the quotient map. Since the class 
of AF--\css is closed by ideals, quotients and extensions, it follows that 
$I+J\in AF(B)$. Suppose in addition that $B$ is a separable \cs, and let 
$D=\{d_n\}_{n\geq 1}$ be 
a countable and dense subset of $\gen\{x\in J:\, J\in AF(B)\}$. Since  
$J_1+\cdots +J_k\in AF(B)$, whenever $J_1,\ldots ,J_k\in AF(B)$, there exists 
an increasing sequence $\{J_n\}_{n\geq 1}\subseteq AF(B)$ such that $d_n\in J_n$, 
$\forall n\geq 1$. It follows that 
$J=\ov{\cup_{n\geq 1}J_n}$ is an AF--ideal that contains any ideal of $AF(B)$, 
and this proves our first assertion. As for the second one, note that, since 
$G$ is separable, then $A$ is separable if and only if so is $\env{A}$. Now, by 
the first part, the result is proved as in Proposition~\ref{prop:class}. 
\end{proof}

\section{Enveloping actions and crossed products}\label{sec:cross}
\par The main result of this section, Theorem \ref{thm:equivmor}, is the 
following: if a \pa $\al$ of $G$ on a \cs $A$ has an enveloping action 
$(\env{\al},\env{A})$, then $A\sd{\al,r}G$ and 
$\env{A}\sd{\env{\al},r}G$ are Morita equivalent. This theorem will be 
obtained as a consequence of a more general result on Fell bundles.  
\par We begin by recalling the definition of the \rr and the reduced cross 
sectional algebra of a Fell bundle. Then we prove \ref{prop:increds}, a result 
that will be of great importance later to prove \ref{thm:equivmor}. 
Finally we apply Theorem \ref{thm:equivmor} to show that any partial \rep 
of an amenable discrete group $G$ is the compression of some unitary \rep of 
$G$. 
\ve
\subsection{Preliminaries on Fell bundles and crossed products by 
partial actions}
\par In \cite{exeng}, the authors defined the reduced cross sectional algebra 
of a Fell bundle, generalizing the definition given by Exel in \cite{examen} 
for bundles over discrete groups. The definition is the following: 
$C^*_r(\mb )$ is the closure of $\La \left(L^1(\mb )\right)$ in 
$\mathcal{L}\big(L^2(\mb)\big)$ , where $\La_f\xi =f*\xi$, $\forall f\in 
C_c(\mb )\subseteq L^1(\mb )$, $\xi\in C_c(\mb )\subseteq L^2(\mb )$. 
When $G$ is a non--discrete group, it is not immediate that $\La$, 
defined by means of convolution of continuous functions, extends to a \rep 
of $L^1(\mb )$. The proof of this fact given below seems to us more direct 
than the one in \cite{exeng}. Next, we show that if $\ma$ is 
a sub--Fell bundle (\ref{df:subfdf}) of $\mb$, 
then $C^*_r(\ma )$ may be considered to be a sub-\cs of $C^*_r(\mb )$.
\par For the general theory of \fbs, also called $C^*$--algebraic bundles, and 
their \reps, the reader is referred to \cite{fd}. Let us suppose that  
$\mb =(B_t)_{t\in G}$ is a Fell bundle over $G$. If $K\subseteq G$ is a compact
subset, then $C_K(\mb)$ denotes the Banach space of continuous sections of 
$\mb$, with the supremum norm. $C_c(\mb)$ denotes the *--algebra of 
continuous sections of $\mb$ that have compact support, with the locally 
convex \ilt defined by the natural inclusions 
$C_K(\mb)\stackrel{\iota_K}{\inc}C_c(\mb)$. 
If $B_e$ is the fiber of $\mb$ over the unit $e$ of $G$, then $L^2(\mb )$ is 
the right Hilbert $B_e$--module obtained by completing $C_c(\mb )$ with 
respect to the $B_e$--inner product: 
$\pr{\xi}{\eta}=\int_G\xi (s)^*\eta(s)ds$.  
If $b_t\in B_t$, let $(\La_{b_t}\xi )\r{s}=b_t\xi (t^{-1}s),\ 
\forall \xi\in C_c(\mb ).$ 
We have that $\La_{b_t}\xi \in C_c(\mb )$, and 
$\text{supp}(\La_{b_t}\xi )\subseteq t\,\text{supp}(\xi )$. In addition:
$ \pr{\La_{b_t}\xi }{\La_{b_t}\xi }
                =\int_G\xi (t^{-1}s)^*b_t^*b_t\xi (t^{-1}s)ds
                   \leq\int_G\xi (t^{-1}s)^*\norm{b_t^*b_t}\xi (t^{-1}s)ds
                =\norm{b_t}^2\pr{\xi}{\xi}$, 
and hence $\La_{b_t}$ may be extended to $L^2(\mb)$. In fact,  
it is easy to see that $\La_{b_t}$ is adjointable, and  
$\La_{b_t}^*=\La_{b_t^*}$.
\par Let us define $\La :\mb\f\mathcal{L}\big(L^2(\mb)\big)$ by  
$b\ff \La_b$. It is immediate that $\La\r{B_t}$ is a bounded linear map, 
$\forall t\in G$, and that $\La_b\La_c=\La_{bc}$, $\forall b,c\in\mb$. 
We will see that $\La$ is also continuous (that is: $\La$ satisfies 
\ref{prop:contrr}).  

\begin{lem}\label{lem:contrr}
If $\xi\in C_c(\mb )$ and $b_t\in\mb$, then for any $\epsilon > 0$ there 
exists an open $U\subseteq\mb$, with $b_t\in U$, such that if $b\in U$, then  
$\norm{\La_{b}\xi-\La_{b_t}\xi}_{\infty}<\epsilon$.
\end{lem}
\begin{proof}
Suppose that there exists $\epsilon >0$ such that 
for any open neighborhood $U$ of $b_t$ there exist $b_{r_U}\in U$ 
and $s_U\in G$ such that  
$\norm{(\La_{b_{r_U}}\xi)(s_U)-(\La_{b_t}\xi)(s_U)}\geq\epsilon$, 
that is, $\norm{b_{r_U}\xi(r_U^{-1}s_U)-b_t\xi(t^{-1}s_U)}\geq\epsilon$. 
Note that $\text{supp}(\La_{b_{r_U}}\xi )\subseteq r_U\,
\text{supp}(\xi )$, $\text{supp}(\La_{b_t}\xi )\subseteq 
t\,\text{supp}(\xi )$. Since $b_{r_U}\f\, b_t$, then $r_U\f\, t$.
Thus there exist a compact set $K\subseteq G$ and a neighborhood  
$U_0$ of $b_t$ such that $\text{supp}(\La_{b_{r_U}}\xi )\subseteq K$, 
$\forall U\subseteq U_0$. Then the net $(s_U)_{U\subseteq U_0}\subseteq K$,
and hence it must have a subnet that is convergent to some $s_0\in K$. We may 
assume without loss of generality that the net itself converges to $s_0$. 
But this is a contradiction, because:
\[
0=\norm{b_{t}\xi(t^{-1}s_0)-b_t\xi(t^{-1}s_0)}=\lim_U 
\norm{b_{r_U}\xi(r_U^{-1}s_U)-b_t\xi(t^{-1}s_U)}\geq\epsilon .
\]    
The contradiction implies that the Lemma is true. 
\end{proof}
\par Recall that $\mathfrak{L}^2(\mb )$ is the completion of $C_c(\mb )$ with  
respect to the norm  
$\norm{\xi}_2=\left(\int_G\norm{\xi (s)}^2ds\right)^{1/2}.$ 
So if $\xi_n\f\xi$ in $\mathfrak{L}^2(\mb )$, with  
$\xi_n,\xi\in C_c(\mb )$, we have that $\xi_n\f\xi$ in $L^2(\mb )$, 
because $\norm{\xi}\leq\norm{\xi}_2$. 
In fact, $\norm{\xi}=\norm{\int_G\xi (s)^*\xi (s)ds}^{1/2}$,
and since $\int_G\xi (s)^*\xi (s)ds$ is a positive element of $B_e$, 
there is a state $\varphi$ of $B_e$ such that $\norm{\int_G\xi(s)^*\xi(s)ds}=
\varphi (\int_G\xi(s)^*\xi(s)ds)$. Then:
\[\left\|\int_G\xi (s)^*\xi (s)ds\right\|
             =\varphi\left(\int_G\xi (s)^*\xi (s)ds\right)
             =\int_G\varphi \left(\xi (s)^*\xi (s)\right)ds
             \leq\int_G\norm{\xi (s)^*\xi (s)}ds
             =\norm{\xi}_2^2\]
On the other hand, it is clear that if $b_r\f b_t$, then   
$\La_{b_r}\xi\f\La_{b_t}\xi$ in $\norm{\cdot}$, because  
$\norm{\xi}_2\leq m\big(\text{supp}(\xi )\big)\norm{\xi}_{\infty}$ 
(here $m$ is the left Haar measure on $G$), and hence, by Lemma 
\ref{lem:contrr}, $\La_{b_r}\xi\f\La_{b_t}\xi$ in $\norm{\cdot}_{\infty}$,
so $\La_{b_r}\xi\f\La_{b_t}\xi$ in $\norm{\cdot}_2$ and $\norm{\cdot}$.
\begin{prop}\label{prop:contrr}
Let $\xi\in L^2(\mb )$. Then the map $\mb\f L^2(\mb )$, given by 
$b\ff\La_b\xi$, is continuous.
\end{prop}
\begin{proof}
Let us fix $b\in\mb$, and let $b_j\f b$ in $\mb$. Given $\epsilon >0$, let 
$\xi_{\epsilon}\in C_c(\mb )$ be such that $\norm{\xi-\xi_{\epsilon}}
<\epsilon$, and let $j_0$ such that $\norm{b_j},\, \norm{b}\leq c$, 
$\forall j\geq j_0$ and some constant $c$. Then, if $j\geq j_0$:
\[\norm{\La_{b_j}\xi -\La_{b}\xi}
    \leq \norm{b_j}\, \norm{\xi -\xi_{\epsilon}} +
     \norm{ \La_{b_j}\xi_{\epsilon}-\La_b\xi_{\epsilon} }+
     \norm{b}\, \norm{\xi_{\epsilon}-\xi}
    < 2c\varepsilon + 
       \norm{ \La_{b_j}\xi_{\epsilon}-\La_b\xi_{\epsilon} }.\]
It follows that $\limsup_i\norm{\La_{b_j}\xi -\La_b\xi}
\leq 2c\epsilon$, $\forall \epsilon >0$, and therefore  
$\La_{b_j}\xi\f\La_b\xi$ in $\norm{\cdot}$.  
\end{proof}
\begin{df}\label{df:rrfib}
(cf. \cite{exeng})
The \rep $\La:\mb\f\mathcal{L}\big(L^2(\mb)\big)$ defined above is  
called the \rr of the Fell bundle $\mb$. That is, $\La_{b_s}$ is the unique 
continuous extension to all of $L^2(\mb)$ of the map $C_c(\mb)\f C_c(\mb)$ 
such that, if $\xi\in C_c(\mb)$, $t\in G$, then 
$\La_{b_s}(\xi)\r{t}=b_s\xi(s^{-1}t)$.
\end{df}   
\begin{thm}\label{thm:regrep}
(cf. \cite{exeng})
There exists a unique non--degenerate \rep 
$\La :L^1(\mb )\f\mathcal{L}\big(L^2(\mb)\big)$, 
given by $f\ff\La_f$, where $\La_f(\xi) =f*\xi$, $\forall f\in C_c(\mb )
\subseteq L^1(\mb )$, $\xi\in C_c(\mb )\subseteq L^2(\mb )$. 
\end{thm}
\begin{proof}
Proposition \ref{prop:contrr} tells us that $\La :\mb\f 
\mathcal{L}\big(L^2(\mb)\big)$ is a Fr\'echet \rep, in the sense of VIII-8.2 
of \cite{fd}. So we may apply VIII-11.3 of \cite{fd}, and conclude 
that $\La$ is integrable. That is, there exists a \rep 
$\La :C_c(\mb )\f B\big(L^2(\mb )\big)$ such that
$\varphi (\La_f )=\int_G\varphi (\La_{f(s)})ds,\ \  \forall f
\in C_c(\mb ),\ \varphi\in B\big(L^2(\mb )\big)'.$
Moreover, $\La$ is unique. We set $\La_f =\int_G\La_{f(s)}ds$. 
\par We want to see that $\forall f\in C_c(\mb )$, it is $\La_f
\in\mathcal{L}\big(L^2(\mb)\big)$. Now, for $\xi ,\eta\in L^2(\mb )$, 
we have that $\pr{\xi}{\La_f(\eta)}=\int_G\pr{\xi}{\La_{f(s)}(\eta)}ds$.
In particular, since $f^*(s)=\Delta(s^{-1})f(s^{-1})^*$, 
\[
\pr{\xi}{\La_{f^*}(\eta )}
        =\int_G\pr{\xi}{\La_{\Delta (s^{-1})f(s^{-1})^*}(\eta )}ds
        =\int_G\Delta (s^{-1})\left(\Delta (s)
                       \pr{\La_{f(s)}(\xi )}{\eta }\right) ds
        =\pr{\La_f(\xi)}{\eta}. \]
Thus $\La_f^*=\La_{f^*}$, and therefore 
$\La_f\in\mathcal{L}\big(L^2(\mb)\big)$.
Moreover, the \rep $\La :C_c(\mb )\f\mathcal{L}\big(L^2(\mb)\big)$ 
is continuous in the norm $\norm{\cdot}_1$: 
$\norm{\pr{\xi}{\La_f(\eta)}}
            \leq\int_G\norm{\xi}\,\norm{\La_{f(t)}\eta}dt
            \leq\int_G\norm{f(t)}\,\norm{\xi}\,\norm{\eta}dt
            =\norm{f}_1\,\norm{\xi}\,\norm{\eta}.$
It follows that $\norm{\La_f}\leq\norm{f}_1$, and hence we may extend $\La$  
by continuity to a \rep of $L^1(\mb )$. This \rep is non--degenerate, because 
$C_c(\mb )*C_c(\mb )$ is dense in $C_c(\mb )$ in the 
inductive limit topology, and therefore also in $L^2(\mb)$. 
\end{proof}
\begin{df}\label{df:reduced}
(cf. \cite{exeng})
The \rep $\La :L^1(\mb )\f\mathcal{L}\big(L^2(\mb)\big)$ defined in Theorem  
\ref{thm:regrep} is called the \rr of $L^1(\mb )$, and 
$C^*_r(\mb )
:=\ov{\La\big(L^1(\mb)\big)}\subseteq\mathcal{L}\big(L^2(\mb)\big)$ 
is called the reduced \cs of $\mb$. If $\mb_{\al}$ is the Fell bundle 
associated with a \pa $\al$ of $G$ on a \cs $A$, then its  
reduced \cs is called the reduced crossed product of $A$ by $\al$, and it is
denoted by $A\sd{\al ,r}G$.   
\end{df}
\begin{rk}\label{rk:note} 
It is shown in \cite{exeng} that if $\al$ is a global action of $G$ 
on $A$, then the usual reduced crossed product agrees with the one defined
in Definition~\ref{df:reduced} (see also \cite{examen}, 3.8). Moreover, the authors
show that if $\pi:\mb\f B(H)$ is a non--degenerate \rep of $\mb$ 
and if $\pi_{\la}$ is the \rep 
$\pi_{\la}:\mb\f B\left(L^2(G,H)\right)$, given by 
$\pi_{\la}(b_t)=\la_t\otimes\pi(b_t)$, where $\la$ is the left \rr
of $G$ on $L^2(G)$, then the integrated \rep of $\pi_{\la}$ 
defines a \rep of $C^*(\mb)$, that factors through $C^*_r(\mb )$. In 
addition, if $\pi\r{B_e}$ is faithful, we have that 
$C^*_r(\mb)\cong\pi_{\la}\big(C^*(\mb)\big)$
(2.15 of \cite{exeng}). Note that if $\pi$ is a degenerate \rep, the result is 
also true, because in this case $\pi=\rho\oplus 0$, the direct sum of the non 
degenerate part of $\pi$ with a zero \rep, and therefore 
$\pi_{\la}=\rho_{\la}\oplus 0$. 
\end{rk}
\par Let us recall the definition of amenable Fell bundle (\cite{examen},
\cite{exeng}):
\begin{df}\label{df:amen}
The \rr $\La$ induces a \rep  of $C^*(\mb )$, also called 
regular and denoted by $\La$. When $\La :C^*(\mb )\f C^*_r(\mb )$ is an 
isomorphism, we say that $\mb$ is amenable. If $\mb_{\al}$ is the Fell bundle 
associated with the \pa $\al$, we say that $\al$ is amenable when 
$\mb_{\al}$ is amenable.
\end{df}
\begin{df}\label{df:subfdf}
If $\mb=(B_t)_{t\in G}$ is a Banach bundle over the Hausdorff space 
$G$ (or a Fell bundle over $G$), we say that 
$\ma\subseteq\mb$ is a sub--Banach bundle of $\mb$ (respectively: a sub--Fell 
bundle of $\mb$) if it is a Banach (respectively: Fell) bundle over $G$ with 
the structure inherited from $\mb$. 
\end{df}
\begin{prop}\label{prop:increds}
If $\ma$ is a sub--Fell bundle of $\mb$, then 
$C^*_r(\ma )\subseteq C^*_r(\mb )$. More precisely: the closure of 
$C_c(\ma)$ in $C_r^*(\mb)$ is naturally isomorphic to $C_r^*(\ma)$.
\end{prop}
\begin{proof}
Let $\pi:\mb\f B(H)$ be a \rep of $\mb$ on the Hilbert space $H$, such that 
$\pi\r{B_e}$ is faithful. Then $\rho:=\pi\r{\ma}:\ma\f B(H)$ is a 
\rep of $\ma$, such that $\rho\r{A_e}$ is faithful. Let  
$\pi_{\la}:\mb\f B\big(\bts{L^2(G)}{}{H}\big)$
such that $\pi_{\la}(b_t)=\la_t\otimes\pi (b_t)$, where  
$\la :G\f B\big(L^2(G)\big)$ is the left \rr of $G$;
that is: $\forall \xi\in L^2(G)$, $\la_t(\xi)\r{s}=\xi (t^{-1}s)$. 
Similarly, define $\rho_{\la}:\ma\f B\big(\bts{L^2(G)}{}{H}\big)$. 
It is clear that $\rho_{\la}=\pi_{\la}\r{\ma}$. Integrating 
$\pi_{\la}$ and $\rho_{\la}$, we obtain \reps of $C^*(\mb)$ and $C^*(\ma)$, 
which we also call $\pi_{\la}$ and $\rho_{\la}$ respectively, and it is clear
again that $\rho_{\la}\r{L^1(\ma )}$ agrees with $\pi_{\la}\r{L^1(\ma)}$. 
Now, by 2.15 of \cite{exeng} (see also the end of Remark \ref{rk:note}), 
we have isomorphisms $\tilde{\pi}_{\la}:  
C^*_r(\mb)\f 
\ov{\pi_{\la}\big(C^*(\mb)\big)}$, and  
$\tilde{\rho}_{\la}:C^*_r(\ma)\f\ov{\rho_{\la}\big(L^1(\ma)\big)}
\cong C^*_r(\ma)$. 
Therefore, $\tilde{\pi}_{\la}^{-1}\tilde{\rho}_{\la}:C_r^*(\ma)\f 
\ov{C_c(\ma)}\subseteq C^*_r(\mb)$ is an isomorphism. Thus, 
$C^*_r(\ma)$ is naturally identified with $\ov{C_c(\ma)}$ in  
$C^*_r(\mb)$.    
\end{proof}
\begin{rk}\label{rk:increds}
When $G$ is a discrete group, there is a shorter proof of Proposition~\ref{prop:increds},  
because in this case the reduced cross sectional algebra is characterized in 
terms of conditional expectations (\cite{examen}).
\end{rk}
\begin{df}\label{df:farrows}
Let $\ma=(A_t)_{t\in G}$ and $\mb =(B_t)_{t\in G}$ be Banach bundles. A \hm  
$\phi :\ma\f\mb$ is a continuous map such that $\phi(A_t)\subseteq B_t$,   
and $\phi\r{A_t}:A_t\f B_t$ is linear and bounded, $\forall t\in G$. If 
$\ma$ and $\mb$ are Fell bundles, we also require that: 
$\phi(xy)=\phi(x)\phi(y)$, $\phi(x^*)=\phi(x)^*$, $\forall x,y\in\ma$.
\end{df}
\begin{rk}\label{rk:farrows}
Let $\phi :\ma\f\mb$ be a \hm. If $f\in L^1(\ma)$, we have that 
$\phi^1(f):G\f\mb$ such that $\phi^1(f)(t)=\phi\big(f(t)\big)$ is in 
$L^1(\mb)$, and $\norm{\phi^1(f)}_1\leq\norm{f}_1$. On the other hand, 
it is clear that $\phi^1$ is a \hm of Banach *--algebras, and therefore it 
extends to a \hm $C^*(\phi ):C^*(\ma)\f C^*(\mb)$. This way we obtain 
a functor from the category of Fell bundles over $G$ to the category 
of \css. A morphism $\phi:\al\to\beta$ between \pas induces a \hm 
$\phi:\mb_{\al}\to\mb_{\beta}$ between 
the corresponding Fell bundles, given by 
$\phi\big((t,a_t)\big)=\big(t,\phi(a_t)\big)$. Thus we have a functor from 
the category of \pas to the category of \css.    
\end{rk}
\ve
\subsection{Morita equivalence between the reduced crossed products}
\begin{df}\label{df:idfdf}
If $\mb=(B_t)_{t\in G}$ is a Fell bundle,  
we say that a sub--Banach bundle $\ma$ of $\mb$ (\ref{df:subfdf}) is a right 
ideal of $\mb$ if 
$\ma\mb\subseteq\ma$, and that it is a left ideal if $\mb\ma\subseteq\ma$; 
$\ma$ is said to be an ideal of $\mb$ if it is both a right and a left ideal 
of $\mb$.    
\end{df}
\par Consider a Fell bundle $\mb=(B_t)_{t\in G}$. If $R$ is a right ideal of 
$B_e$ and we define $R_t:=\ov{\gen}RB_t$, $\forall t\in G$, we obtain a right 
ideal $\mr=(R_t)_{t\in G}$ of $\mb$. In a similar way we may use left ideals 
of $B_e$ to define left ideals of $\mb$. If $I\id B_e$, then  
$\mi=(I_t)_{t\in G}$, where $I_t=IB_t$, is a right ideal of $\mb$, but in 
general not an ideal. For this it is necessary and sufficient that $I$ is a 
$\mb$--\textit{invariant} ideal of $B_e$, that is, $IB_t=B_tI$, 
$\forall t\in G$. Conversely, if $\mi=(I_t)_{t\in G}$ is a given ideal of 
$\mb$, and $I=I_e$, then $I$ is a $\mb$--invariant ideal of $B_e$. Moreover, 
these correspondences establish an inclusion--preserving bijection between 
$\mb$--invariant ideals of $B_e$ and ideals of $\mb$.   
\par Let $\mathcal{V}$ be a local base of precompact and symmetric 
neighborhoods of the unit $e\in G$, directed by the relation: 
$V\geq V'\iff V\subseteq V'$. Then there exists an approximate unit 
$(f_V)_{V\in\mv}$ of $L^1(G)$ contained in $C_c(G)$, and such that 
$\text{supp}(f_V)\subseteq V$, $f_V\geq 0$, and $\int_Gf_V(s)ds=1$, 
$\forall V\in\mv$. 

\begin{lem}\label{lem:lema}
Let $\mb=(B_x)_{x\in X}$ be a Banach bundle, $G$ a group, with both $X$ and $G$ locally compact and Hausdorff. Suppose $b:G\times X\f\mb$ is a compactly  
supported continuous function such that 
$b(r,x)\in B_x$, $\forall (r,x)\in G\times X$.  
Let $(f_V)_{V\in\mv}$ be an approximate unit of $L^1(G)$ as described just 
above, and define, for each pair $(V,r)\in\mv\times G$, the function  
$b_{V,r}:X\f\mb$, such that $b_{V,r}(x)=\int_Gf_V(r^{-1}s)b(s,x)ds$.
Then:
\be 
 \item $b_{V,r}\in C_c(\mb)$, $\forall V\in \mv$, $r\in G$.
 \item $\lim_Vb_{V,r}=b_r$ in the inductive limit topology,  
       where $b_r:X\f\mb$ is given by $b_r(x)=b(r,x)$.
\ee
 \end{lem}
\begin{proof}
The function $\mu:G\times G\times X\f\mb$ such that  
$(r,s,x)\ff f_V(r^{-1}s)b(s,x)$ is continuous and $\mu(r,s,x)\in B_x$, 
$\forall r,s\in G$, $x\in X$. So by \cite{fd}, II-15.19, the function 
$G\times X\f\mb$ given by $(r,x)\ff\int_Gf_V(r^{-1}s)b(s,x)ds$ is continuous. 
In particular, $b_{V,r}$ is continuous, and thus (1) is proved. 
\par As for (2), let $K_1\subseteq G$, $K_2\subseteq X$ be compact sets such that 
$\text{supp}(b)\subseteq K_1\times K_2$. Since the function $b$ is continuous, 
each $b_r:X\f\mb$ is a continuous section with 
$\text{supp}(b_r)\subseteq K_2$, 
so we have a function $G\f C_c(\mb)$ defined by $r\f b_r$, that is supported  
in $K_1$. This map is continuous.
In fact, let $r_0\in G$. Since the function $G\times X\f\R$ that maps  
$(r,x)$ into $\norm{b(r_0,x)-b(r,x)}$ is continuous and equal to zero in 
every $(r_0,x)$, for $x\in K_2$, there exist open neighborhoods   
$U_x$ of $r_0$, $V_x$ of $x$, such that   
$\norm{b(r_0,y)-b(r,y)}<\epsilon$, $\forall r\in U_x$, $y\in V_x$.  
Since the $V_x$ cover the compact set $K_2$, there exists a finite 
subcovering $V_{x_1},\ldots ,V_{x_n}$ of $K_2$. 
Let $U=\bigcap_{i=1}^nU_{x_i}$ and pick $r\in U$, $x\in X$. 
Then either $x\notin K_2$, and then $b_r(x)=0$   
$\forall r\in G$, or $x\in K_2$, and therefore $x\in V_{x_i}$ for some $i$. 
In this case, $(r,x)\in U_{x_i}\times V_{x_i}$, and hence 
$\norm{b_{r_0}(x)-b_r(x)}<\epsilon$, so 
$\norm{b_{r_0}-b_r}_{\infty}<\epsilon$. 
It follows that $r\ff b_r$ is continuous, and hence uniformly continuous, 
because it has compact support. Thus, for any $\epsilon >0$, there exists
$V_0\in\mv$ such that if $r^{-1}s\in V_0$, then  
$\norm{b_r-b_s}_{\infty}<\epsilon$. So we have, for any $V\in\mv$ such that 
$V\geq V_0$, and for all $x$:
\[\norm{(b_r -b_{V,r})(x)}
                          =\norm{\int_Gf_V(r^{-1}s)\big(b(r,x)-b(s,x)\big)ds}
                          \leq\int_Vf_V(r^{-1}s)\norm{b(r,x)-b(s,x)}ds
                          <\epsilon \]
Therefore, since $\text{supp}(b_r)$, $\text{supp}(b_{V,r})\subseteq K_2$ and  
$\norm{b_r-b_{V,r}}_{\infty}<\epsilon$, $\forall V\geq V_0$, it follows that 
$\lim_Vb_{V,r}=b_r$ in the inductive limit topology, $\forall r\in G$. 
\end{proof}
\begin{thm}\label{thm:equimorl1}
Let $\mb=(B_t)_{t\in G}$ be a Fell bundle, $\me=(E_t)_{t\in G}$ a right ideal 
of $\mb$ (\ref{df:idfdf}), and $\ma=(A_t)_{t\in G}$ a sub--Fell bundle 
(\ref{df:subfdf})of $\mb$ contained in $\me$. 
If $\ma\me\subseteq\me$ and $\me\me^*\subseteq \ma$, we have:
\be
 \item $L^1(\ma)*L^1(\me)\subseteq L^1(\me)$
 \item $L^1(\me)*L^1(\mb)\subseteq L^1(\me)$.
 \item $L^1(\me)*L^1(\me)^*=L^1(\ma)$
 \item If $\gen(B_t\bigcap\me^*\me)$ is dense in $B_t$, $\forall t\in G$, 
       then $\ov{\gen}L^1(\me)^**L^1(\me)=L^1(\mb)$.
\ee 
\end{thm}
\begin{proof} 
It is straightforward to check that $C_c(\ma)*C_c(\me )\subseteq C_c(\me )$ 
and $C_c(\me )*C_c(\mb)\subseteq C_c(\me )$, and from these facts the two 
first inclusions follow easily. Let us prove the third assertion. Since 
$\ma\subseteq\me\subseteq\mb$, we have isometric inclusions 
$L^1(\ma)\subseteq L^1(\me )\subseteq L^1(\mb)$. $L^1(\ma)$ is a sub-*-Banach 
algebra with approximate unit of $L^1(\mb)$, and it is contained in the right 
ideal $L^1(\me )$ of $L^1(\mb)$. Thus, by 1. and the Cohen-Hewitt theorem,
      $ L^1(\ma)=L^1(\ma)*L^1(\ma)^*\subseteq 
         L^1(\me )*L^1(\me )^*.$
On the other hand, if $\xi$, $\eta\in C_c(\me )$, $t\in G$:
\[\xi *\eta^*(t) = \int_G\xi(s)\eta^*(s^{-1}t) ds
                 = \int_G\xi(s)\Delta(t^{-1}s)\eta(t^{-1}s)^* ds
                 = \int_G\Delta(t^{-1}s)\xi(s)\eta(t^{-1}s)^*ds\in A_t,\]
because $\xi (s)\eta(t^{-1}s)^*\in E_sE_{t^{-1}s}^*\in\ma\bigcap B_t=A_t$, 
$\forall s,t\in G$. It follows that $\xi *\eta^*\in C_c(\ma)$, and hence that 
$L^1(\me )*L^1(\me )^*\subseteq L^1(\ma)$.
\par Consider now $\xi$, $\eta\in C_c(\me )$, $t\in G$. We have:
\[\xi^**\eta (t)=\int_G\xi^*(s)\eta(s^{-1}t)ds
                =\int_G\Delta(s^{-1})\xi (s^{-1})^*\eta (s^{-1}t)ds
                =\int_G\Delta (s^{-1})\xi (s^{-1})^*\eta (s^{-1}t)ds.\]
\par Let $\{ (f_V,V)\}_{V\in\mv}$ be an approximate unit of $L^1(G)$ as in 
Lemma \ref{lem:lema}. For $\xi\in C_c(\me)$, $V\in\mv$, $r\in G$, define  
$\xi_{V,r}:G\f \mb$ by $\xi_{V,r}(s)=\Delta (s^{-1})f_V(r^{-1}s^{-1})\xi(s).$ 
Then $\xi_{V,r}\in C_c(\me )$, and we have:
\[\xi_{V,r}^**\eta (t)=\int_G\Delta (s^{-1})\Delta (s)
                       f_V(r^{-1}s)\xi(s^{-1})^*\eta (s^{-1}t)ds
                      =\int_Gf_V(r^{-1}s)\zeta_{\xi,\eta}(s,t)ds,\]
where $\zeta_{\xi,\eta}:G\times G\f\mb$ is such that  
$\zeta_{\xi,\eta}(s,t)=\xi(s^{-1})^*\eta(s^{-1}t)$. Note that 
$\zeta_{\xi,\eta}$ is continuous of compact support: 
$\text{supp}(\zeta_{\xi,\eta})\subseteq 
(\text{supp}(\xi))^{-1}\times (\text{supp}(\xi))^{-1}\text{supp}(\eta)$. By Lemma \ref{lem:lema}, 
we see that $\lim_V\xi_{V,r}^**\eta =\zeta_{\xi,\eta,r}$ in the inductive 
limit topology $C_c(\mb)$, and hence also in $L^1(\mb)$. So, we have that  
$\big(\ov{\gen}L^1(\me)^**L^1(\me)\big)\bigcap C_c(\mb)\supseteq Z$, 
where $Z=\gen\{\zeta_{\xi,\eta,r}:\ \xi,\eta\in C_c(\me),\, r\in G\}$.   
\par To see that $\ov{\gen}L^1(\me)^*L^1(\me)=L^1(\mb)$, it is sufficient  
to see that $Z$ is dense in $C_c(\mb)$ in the inductive limit topology. 
By II-14.6 of \cite{fd}, for this is enough to verify that:
(a) $Z(t)$ is dense in $B_t$, $\forall t\in G$, where 
$Z(t)=\{\zeta(t):\ \zeta\in Z\}$ and (b) if $g:G\f \C$ is continuous, 
then $g\zeta\in Z$, $\forall\zeta\in Z$.\\
(a) We have: $Z(t)\supseteq\{\zeta_{\xi,\eta,r}:\ \xi,\eta\in C_c(\me),
            r\in G\}=\{\xi(r^{-1})^*\eta(r^{-1}t):\ \xi,\eta\in C_c(\me),
            r\in G\}$. Therefore, $\ov{Z(t)}\supseteq\ov{\gen}
            \{E_{r^{-1}}^*E_{r^{-1}t}:\ r\in G\}=\ov{\gen}\big(B_t\bigcap 
            \me^*\me\big)=B_t$ by the hypothesis in 4.\\ 
(b) Let $g:G\f\C$ be a continuous function, $\xi$, $\eta\in C_c(\me)$, 
            $r,t\in G$. Defining $g^r:G\f\C$ as $g^r(s)=g(rs)$ we have:  
            $ (g\zeta_{\xi,\eta,r})(t)
                   =g(t)\xi(r^{-1})^*\eta(r^{-1}t)
                   =\xi(r^{-1})^*g^r(r^{-1}t)\eta(r^{-1}t)
                   =\zeta_{\xi,g^r\eta,r}(t).$
            Since $\zeta_{\xi,g^r\eta,r}\in C_c(\me)$, we conclude that 
            $\zeta_{\xi,g^r\eta,r}\in Z$, closing the proof.
\end{proof}
\begin{cor}\label{cor:equimorl1}
Let $\mb=(B_t)_{t\in G}$ be a Fell bundle, 
$\me=(E_t)_{t\in G}$ a right ideal of $\mb$, and $\ma=(A_t)_{t\in G}$ 
a sub--Fell bundle of $\mb$ contained in $\me$. If $\ma\me\subseteq\me$ and 
$\me\me^*\subseteq \ma$, we have that $C^*_r(\ma)$ is a hereditary sub-\cs 
of $C^*_r(\mb)$, and if $\gen\big(B_t\bigcap\me^*\me\big)$ is dense in $B_t$, 
for each $t\in G$, then $C^*_r(\ma)$ and $C^*_r(\mb)$ are Morita equivalent 
via the right ideal $C^*_r(\me):=\ov{C_c(\me)}\subseteq C^*_r(\mb)$ of 
$C^*_r(\mb)$. 
\end{cor}
\begin{proof}
By \ref{prop:increds}, $C^*_r(\ma)$ is naturally isomorphic to 
the closure of $C_c(\ma)$ in $C^*_r(\mb)$. By Theorem~\ref{thm:equimorl1} (2), 
$L^1(\me)$ is a right ideal of $L^1(\mb)$, and hence its closure 
$C^*_r(\me)$ in $C^*_r(\mb)$ is a right ideal of $C^*_r(\mb)$. Now, it 
follows from 3. of \ref{thm:equimorl1} that 
$C^*_r(\ma)=C^*_r(\me)C^*_r(\me)^*$, 
and therefore $C^*_r(\ma)$ is a hereditary sub-\cs of $C^*_r(\mb)$. Finally, 
the last assertion follows immediately from (4) of Theorem~\ref{thm:equimorl1}. 
\end{proof}
\begin{cor}\label{cor:amenab}
Let $\mb=(B_t)_{t\in G}$ be a Fell bundle, $\me=(E_t)_{t\in G}$ a right ideal 
of $\mb$, and $\ma=(A_t)_{t\in G}$ a sub--Fell bundle of $\mb$ contained 
in $\me$. If $\ma\me\subseteq\me$ and $\me\me^*\subseteq \ma$, and if 
$\gen\big(B_t\bigcap\me^*\me\big)$ is dense in $B_t$, for all  
$t\in G$, then $\mb$ is amenable whenever $\ma$ is amenable.
\end{cor}
\begin{proof}
Suppose that $\ma$ is amenable, and let $\norm{\cdot}_{\text{m\'ax}}$ 
be the norm on $C^*(\mb)$ and  $\norm{\cdot}_r$ the norm on  
$C^*_r(\mb)$. The closure of $C_c(\ma)$ in $C^*_r(\mb)$ is $C^*_r(\ma)$, 
by Proposition \ref{prop:increds}. We also have that the closure of 
$C_c(\ma)$ in $C^*(\mb)$ is $C^*(\ma)=C^*_r(\ma)$, because any \rep of  
$L^1(\mb)$ induces a \rep of $L^1(\ma)$ by restriction, and therefore the 
norm of $C^*(\ma)$ is greater or equal to $\norm{\cdot}_{\text{max}}$. 
The amenability of $\ma$ implies that these two norms are equal. 
Let $\xi\in C_c(\me )$. 
Since $C^*(\ma)=C^*_r(\ma)$ by assumption, and $\xi *\xi^*\in C_c(\ma)$ by 
\ref{thm:equimorl1}, we have: 
\[\norm{\xi}_r^2=\norm{\xi^**\xi}_r
                =\norm{\xi*\xi^*}_r
                =\norm{\xi*\xi^*}_{\text{m\'ax}}
                 =\norm{\xi^**\xi}_{\text{m\'ax}}
                =\norm{\xi}_{\text{m\'ax}}^2, \]
and therefore $\norm{\xi}_r=\norm{\xi}_{\text{m\'ax}}$. Then, the completions 
of $C_c(\me )$ with respect to $\norm{\cdot}_{\text{m\'ax}}$ and  
$\norm{\cdot}_r$ agree. Let us denote this completion by $E$. We have 
that $E$ is a full Hilbert module over both $C^*(\mb)$ and  
$C^*_r(\mb)$, and hence $C^*(\mb)= C^*_r(\mb)$. 
This shows that $\mb$ is amenable.
\end{proof}
\begin{thm}\label{thm:equivmor}
Let $\beta$ be a continuous action of $G$ on a \cs $B$, $I\id B$, 
$\al =\beta\r{I}$. Then $I\sd{\al ,r}G$ is a hereditary sub--\cs of 
$B\sd{\beta ,r}G$. If, in addition, $[\beta (I)]$ is dense in $B$, i.e. 
$\beta$ is the enveloping action of $\al$, then $I\sd{\al ,r}G$ and 
$B\sd{\beta ,r}G$ are Morita equivalent.
\end{thm}
\begin{proof}
Let $\mb_{\beta}=(B_t)_{t\in G}$ be the Fell bundle associated with $\beta$, 
$\mb_{\al}=(A_t)_{t\in G}$ the Fell bundle associated with $\al$, and 
$\me=(E_t)_{t\in G}$, where $E_t=\{(t,x)\in\mb_{\beta}:\ x\in I\}$. 
It is clear that $\mb_{\al}\subseteq\me\subseteq\mb_{\beta}$ as Banach
bundles, and that $\mb_{\al}$ is a sub--Fell bundle of $\mb_{\beta}$. 
Moreover, if $(r,a_r)\in A_r$, $(s,x_s)\in E_s$, $(t,y_t)\in E_t$, 
$(u,b_u)\in B_u$, we have: 
\bi
 \item $(r,a_r)(s,x_s)=(rs,\al_r\big(\al_{r^{-1}}(a_r)x_s\big))\in E_{rs}$, 
       because $\al_r\big(\al_{r^{-1}}(a_r)x_s\big)\in I$. Therefore: 
       $\mb_{\al}\me\subseteq\me$ 
 \item $(s,x_s)(t,y_t)^*=(s,x_s)(t^{-1},\beta_{t^{-1}}(y_t^*))
       =(st^{-1},x_s\beta_{st^{-1}}(y_t^*))\in A_{st^{-1}}$, because 
       $x_s\beta_{st^{-1}}(y_t^*))$ belongs to $I\bigcap\beta_{st^{-1}}(I)=
       D_{st^{-1}}$. Consequently, $\me\me^*\subseteq\mb_{\al}$. 
 \item $(s,x_s)(u,b_u)=(su,x_s\beta_s(b_u))\in E_u$ because $I$ is an ideal  
       and $x_s\in I$. Thus, $\me\mb_{\beta}\subseteq\me$. 
\ei
Thus, we may apply Corollary \ref{cor:equimorl1}, concluding that 
$I\sd{\al,r}G=C^*_r(\mb_{\al})$ is a hereditary sub-\cs of 
$C^*_r(\mb_{\beta})=B\sd{\beta,r}G$.  
\par Suppose now that $\beta$ is the enveloping action of $\al$, that is,  
$[\beta(I)]$ is dense in $B$. To see that $I\sd{\al,r}G\sM B\sd{\beta,r}G$, 
it is sufficient to show that $\gen\big(\me^*\me\bigcap B_t\big)$ is dense
in $B_t$, for all $t\in G$. Now: 
$(s,x_s)^*(t,y_t)
=(s^{-1},\beta_{s^{-1}}(x_s^*))(t,y_t)
=(s^{-1}t,\beta_{s^{-1}}(x_s^*y_t))$. Therefore 
$\gen(\me^*\me\bigcap B_t)=\{ (t,\beta_v(x'y')):\ v\in G,\, x',y'\in I\}$.
By the Cohen -Hewitt theorem, every $z\in I$ may be written in the form  
$z=x'y'$, for some $x',y'\in I$. So the set above is exactly 
$\{ t\}\times [\beta(I)]$, which is dense in $B_t$. 
\end{proof}
\begin{cor}\label{cor:equivmor}
If $\beta$ is the enveloping action of an amenable \pa $\al$ 
(\ref{df:amen}), then $\beta$ is amenable as well. 
\end{cor}
\begin{proof}
It follows immediately from \ref{cor:amenab}. 
\end{proof}
\ve
\subsection{An application: dilations of partial representations}
\par Closing this section, we apply our results to show that any partial 
\rep of an amenable discrete group may be dilated to a unitary \rep 
of the group. Recall from \cite{exinv} that a partial \rep of a 
discrete group $G$ on the Hilbert space $H$ is a map $u:G\f B(H)$ such that, 
for $t,s\in G$:\  
(i) $u_e=id_H$;\ \ 
(ii) $u_{t^{-1}}=u_t^*$;\ \  
(iii) $u_su_tu_{t^{-1}}=u_{st}u_{t^{-1}}$.
The conditions above imply that $u_t$ is a partial isometry, and also that 
$u_{s^{-1}}u_su_t=u_{s^{-1}}u_{st}$, $\forall s,t\in G$. The partial \reps 
of $G$ are in one to one correspondence with the non--degenerate \reps 
of its partial \cs $C^*_p(G)$, which is constructed as follows. Let 
$X_t=\{\om\in 2^G:\, e,t\in\om\}$ with the product topology, 
and $\al_t:X_{t^{-1}}\f X_t$ such that $\al_t(\om )=t\om$, 
$\forall \om\in X_{t^{-1}}$. Then $\al$ is a \pa on $X:=X_e$, and $C^*_p(G)$ 
is defined to be the corresponding crossed product $C(X)\sd{\al}G$, where 
we are also denoting by $\al$ the \pa induced on $C(X)$. 
\begin{prop}\label{prop:dil}
Let $G$ be a discrete amenable group, and $u:G\f B(H)$ a partial \rep. 
Then there exist a Hilbert space $\hh$, which contains $H$ as a Hilbert
subspace, and a unitary \rep $\hu :G\f B(\hh)$, such that  
$ u_t=P\hu_ti,\ \ \forall t\in G$, where $P:\hh\f H$ is the orthogonal
projection of $\hh$ on $H$, and $i:H\f \hh$ is the natural inclusion. 
In particular, the partial \reps of a discrete amenable group are positive 
definite maps.  
\end{prop}
\begin{proof}
Let $\al$ be the \pa described before. Since $G$ is amenable, we have that 
$C^*_p(G)=C(X)\sd{\al ,r}G$. First of all, note that $\al$ has an \ea 
$\env{\al}$ acting on $\env{X}=2^G\setminus\{\emptyset\}$, and given by the 
same formula as $\al$. Let $\mb_{\al}$ be the Fell bundle over $G$ 
associated with $\al$, $1_t$ the characteristic function of $X_t$ and, if 
$a_t\in C(X_t)$, let $a_t\delta_t\in C_c(\mb_{\al})$ be  
defined as $a_t\delta_t(s)=\delta_{t,s}a_t$, where $\delta_{t,s}$ is the 
Kronecker symbol. By 6.5 of \cite{exinv}, $u$ defines a unique non--degenerate 
\rep $\pi_u:C_p^*(G)\f B(H)$, such that $\pi_u(1_t\delta_t)=u_t$, 
$\forall t\in G$. In particular, $\pi_u(1_e\del_e)=id_H$. By Theorem 
\ref{thm:equivmor}, we have that $C^*_p(G)$ is a hereditary sub--\cs of 
$C(\hex)\sd{\hal ,r}G$, which is equal to $C(\hex)\sd{\hal}G$ because of the 
amenability of $G$. Moreover, $C^*_p(G)$ is Morita equivalent to 
$C(\hex)\sd{\hal}G$. Therefore, there exist a Hilbert space $\hh$ and a 
non--degenerate \rep $\hpi_u:C(\hex_G)\sd{\hal}G\f B(\hh )$, such that $H$ 
is a Hilbert subspace of $\hh$ and $\pi_u$ is the compression of $\hpi_u$ to 
$H$, i.e.: $\pi_u(x)=P\hpi_u(x)i$, $\forall x\in C^*_p(G)$, where 
$P:\hh\f H$ is the orthogonal projection and $i=P^*:H\f\hh$ is the 
natural inclusion (\cite{fd}, XI-7.6). 
\par Now, $\hpi_u=\hphi\times\hu$, for some covariant \rep  
$(\hphi ,\hu )$ of the dynamical system $(C(\hex),\hal ,G)$; in particular, 
$\hu :G\f B(\hh )$ is a unitary \rep of $G$. 
\par Note that $X$ is a clopen subset of $\hex$, so $1_e\in C(\hex)$, 
and we may compute, in $C(\hex)\sd{\hal ,r}G$: 
$(1_e\del_e)(1\del_t)(1_e\del_e)
=(1_e\del_t)(1_e\del_e)
=1_e\hal_t(1_e)\del_t
=1_t\del_t.$
Therefore:
\[u_t
=\pi_u(1_t\del_t)
=P\hpi_u(1_t\del_t)\r{H}
=P\hpi(1_e\del_e)\hpi(1\del_t)\hpi(1_e\del_e)\r{H}
=P\hpi(1_e\del_e)\hu_t\hpi(1_e\del_e)\r{H}.
\]
Observe now that $\hpi(1_e\del_e)$ is an orthogonal projection
such that $P\hpi(1_e\del_e)\r{H}=\pi(1_e\del_e)=id_H$, and hence 
$\hpi(1_e\del_e)$ is greater 
or equal to the orthogonal projection $Q\in B(\hh)$ with image $H$.  
Thus, we have that: $Q\hpi(1_e\del_e)=Q=\hpi(1_e\del_e)Q$.  
On the other hand, it is clear that $PQ=P$, $Q=Qi$, and 
consequently: 
\[ P\hpi(1_e\del_e)
=(PQ)\hpi(1_e\del_e)
=P\big(Q\hpi(1_e\del_e)\big)
=PQ
=P,\]
\[ \hpi(1_e\del_e)i
=\hpi(1_e\del_e)(Qi)
=\big(\hpi(1_e\del_e)Q\big)i
=Qi
=i.
\]
It follows that 
$ u_t
=P\hpi(1_e\del_e)\hu_t\hpi(1_e\del_e)i
=P\hu_ti$.
\end{proof}
\ve
\section{Morita equivalence of partial actions and Morita enveloping 
actions}\label{sec:mea}
\par We have seen previously that there exist \pas on \css that have no \eas. 
The aim of this section is to introduce a weaker notion of \ea, so that 
any \pa has a unique ``weak'' \ea, and Theorem \ref{thm:equivmor} is still
valid. 
\par To this end we define and study Morita equivalence of \pas, and show 
that the reduced crossed products of Morita equivalent \pas are Morita 
equivalent.  Then we define the so called Morita enveloping actions. If 
$\al$ is a \pa, we say that $\beta$ is an enveloping action up to Morita 
equivalence of $\al$, or just a \mea of $\al$, if $\beta$ is the \ea of 
a \pa that is Morita equivalent to $\al$. For this notion we have a result 
analogous to \ref{thm:equivmor}: Proposition \ref{prop:equivmor2}. 
The investigation of the existence and uniqueness of \meas is postponed until 
Section~\ref{sec:exun}.

\subsection{Hilbert modules and $C^*$--ternary rings}
In the next subsection we will introduce the Morita equivalence between \pas, 
and to do this it will be convenient
to use $C^*$--ternary rings (\cts for short). It would be possible to use 
just Hilbert bimodules, but we prefer to view a Hilbert module not as 
a space where a \cs is acting on, but rather as an object that has almost the 
status of a \cs. So, we will quickly see now some basic facts about \cts 
that will be needed later. 
\par Let us suppose that $(E,\pr{\cdot}{\cdot})$ is a full right 
Hilbert $B$-module, and let 
$A=\k{E}$ be the corresponding \cs of compact operators, that is, the 
ideal of $\adj{E}$ generated by $\{\theta_{x,y}:\, x,y\in E\}$, where 
$\theta_{x,y}(z)=x\pr{y}{z}_r$. Defining $\pr{x}{y}_l=\theta_{x,y}$, we
have that $E$ is a full left Hilbert $A$-module, and that $E$ is an 
$(A-B)$-bimodule that satisfies $\pr{x}{y}_lz=x\pr{y}{z}_r$, $\forall 
x,y,z\in E$; we will say that $E$ is a full Hilbert $(A-B$)-bimodule. 
So, defining $(x,y,z)=x\pr{y}{z}_r$, we have a ternary product 
$(\cdot,\cdot,\cdot)$ on $E$ that relates the actions of $A$ and $B$ on $E$, 
and also the left and right inner products on $E$. The object 
$\big(E,(\cdot,\cdot,\cdot)\big)$ is a \ct, and determines the 
pairs $(A,\pr{\cdot}{\cdot}_l)$ and $(B,\pr{\cdot}{\cdot}_r)$ up to 
isomorphisms. This fact was proved in \cite{z}, where the notion of \ct 
was introduced. Let us recall the exact definitions.
\begin{df}\label{df:1} 
    A *-ternary ring is a complex linear space $E$ with a 
    transformation $\mu :E\times E\times E\f E$, called ternary 
    product on $E$, such that $\mu$ is linear in the odd variables and 
    conjugate linear in the second one, and such that:  
     $\mu\big(\mu(x,y,z),u,v\big)
        =\mu\big(x,\mu(u,z,y),v\big)
        =\mu\big(x,y,\mu(z,u,v)\big),
    \ \forall x,y,z,u,v\in E$. 
    A \hm $\phi:(E,\mu )\f (F,\nu )$ of *-ternary rings is a linear 
    transformation from $E$ to $F$ such that
    $\nu \big(\phi (x),\phi (y),\phi (z)\big)
        =\phi\big(\mu(x,y,z)\big),\, \forall
                                    x,y,z\in E.$
    We will write $(x,y,z)$ instead of $\mu (x,y,z)$. 
    \par A $C^*$-norm on a *-ternary ring $E$ is a norm such that 
    $\norm{(x,y,z)}\leq\norm{x}\,\norm{y}\,\norm{z}$ and 
    $\norm{(x,x,x)}=\norm{x}^3$, $\forall x,y,z\in E$. 
    We then say that $(E,\norm{\cdot})$ is a pre-\ct, and that it is a \ct 
    if it is complete.
    \par A \rep of a $*$-ternary ring $(E,\mu)$ on the Hilbert spaces 
    $H$ and $K$ is a \hm $\pi:E\f B(H,K)$, where in the last space we 
    consider the ternary product given by $(R,S,T)=RS^*T$. $\pi(E)$ is 
    called a ternary ring of operators, or just TRO.
\end{df}
\par \cts and TRO's were studied by Zettl in \cite{z}. He proved that 
every \ct $(E,\mu)$ may be uniquely decomposed as a direct sum 
$E=E^+\oplus E^-$, where $(E^+,\mu\r{E^+})$ and $(E^-,-\mu\r{E^-})$ are 
isomorphic to closed TRO's. We will say that a \ct $E$ is \textit{positive}
if $E=E^+$. By a result of Zettl's (\cite{z}, 3.2), positive \cts correspond 
exactly to full Hilbert bimodules, that is, there exist, up to isomorphisms, 
unique pairs $(E^l,\pr{\cdot}{\cdot}_l)$ and $(E^r,\pr{\cdot}{\cdot}_r)$ 
such that $E$ is a full Hilbert $(E^l-E^r)$-bimodule, and 
$\pr{x}{y}_lz=\mu(x,y,z)=x\pr{y}{z}_r$, $\forall x,y,z\in E$.
\begin{prop}\label{prop:z(pi)}
Let $\pi :E\f F$ be a \hm of *-ternary rings 
(\ref{df:1}), where $E$ and $F$ are \cts.
Then $\pi$ is a contraction, and there exists a unique \hm 
$\pi^r:E^r\f F^r$ such that $\pi^r(\pr{x}{y}_r)
=\pr{\pi (x)}{\pi (y)}_r$, $\forall x,y\in E$. Consequently, we have 
that $\pi (xb)=\pi (x)\pi^r(b)$, $\forall x\in E$, $b\in E^r$. 
If $\pi$ is injective (surjective, an isomorphism), then  
$\pi$ is isometric, and $\pi^r$ is injective (respectively surjective, 
an isomorphism).  
\end{prop}
\begin{proof} 
If $\pi^r$ exists, it must be   
$\pi^r (\pr{x}{y}_r)=\pr{\pi (x)}{\pi (y)}_r$, $\forall x,y\in E$. 
Therefore we must see that $\sum_i\pr{x_i}{y_i}_r=0$ implies that   
$\sum_i\pr{\pi (x_i)}{\pi (y_i)}_r=0$. 
Now, if $\sum_i\pr{x_i}{y_i}_r=0$:
\[\big(\sum_i\pr{\pi (x_i)}{\pi (y_i)}_r\big)^*
\big(\sum_j\pr{\pi (x_j)}{\pi (y_j)}_r\big)
=\sum_i\pr{\pi (y_i)}{\pi (x_j\sum_j\pr{x_j}{y_j}_r)}_r
=0.\]
Thus, we have a \hm of *-algebras $\pi^r:\gen\pr{E}{E}_r\f F^r$. 
Since $\gen\pr{E}{E}_r$ is a dense *-ideal of $E^r$, then $\pi^r$ has a 
unique extension to a \hm $\pi^r:E^r\f F^r$ (\cite{fd}, VI-19.11).
\par Now, if $x\in E$: 
$\norm{\pi (x)}^2=\norm{\pr{\pi (x)}{\pi (x)}_r}
                  =\norm{\pi^r(\pr{x}{x}_r)}
                  \leq\norm{\pr{x}{x}_r}
                  =\norm{x}^2$, and hence $\pi$ is a contraction. 
In particular, $\pi$ is continuous, and therefore  
$\pi (xb)= \pi (x)\pi^r(b)$, $\forall x\in E$, $b\in E^r$, because this 
is true for each $b\in\gen\pr{E}{E}_r$, a dense subset of $E^r$.
\par If $\pi$ is surjective, it is clear that so is $\pi^r$. 
Suppose that $\pi$ is injective, and let $b\in\ker\pi^r$.
Then $\pi (xb)=0$, $\forall x\in E$, and hence $xb=0$, $\forall x\in E$, 
and therefore $\gen\pr{E}{E}_rb=0$. It follows that $E^rb=0$, and 
hence that $b=0$. Consequently $\pi^r$ is injective, and therefore isometric. 
Thus:
$\norm{\pi (x)}^2=\norm{\pr{\pi (x)}{\pi (x)}_r}
                  =\norm{\pi^r(\pr{x}{x}_r)}
                  =\norm{\pr{x}{x}_r}
                  =\norm{x}^2,$
so $\pi$ is isometric.
\end{proof}
\par Zettl's results together with Proposition \ref{prop:z(pi)} imply 
that, up to the fact that $E^r$ is determined up to isomorphisms, we have a 
functor $E\ff E^r$, 
$(E\stackrel{\pi}{\f}F)\ff(E^r\stackrel{\pi^r}{\f}F^r)$ from the category 
of \cts to the category of \css. Of course, we have a left version of this
fact.  
\begin{df}\label{df:ideals}
Let $E$ be a \ct, and $F\subseteq E$ a closed subspace. We say that
$F$ is an ideal of $E$ iff $(E,E,F)\subseteq F$ and $(F,E,E)\subseteq F)$. 
We write $F\id E$ to indicate that $F$ is an ideal of $E$ and 
$\mi(E)$ to denote the set of ideals of $E$.
\end{df}
\par It is not hard to see that $F\id E$ iff $(E,F,E)\subseteq F$. 
Let us suppose that $E$ is a positive \ct, so that $E$ is a 
full Hilbert $(E^l-E^r)$bimodule. Then it is easy to see that $F$ is 
an ideal of $E$ iff $F$ is a closed sub-$(E^l-E^r)$-bimodule of $E$.
Therefore, rephrasing a well known result in our context, we have 
(see for instance \cite{rw}, 3.22): 
\begin{prop}\label{prop:lattice}
Let $E$ be a positive \ct. Then the correspondence 
$F\ff F^r=\ov{\gen}\pr{F}{F}_r$ is a lattice isomorphism between $\mi(E)$ 
and $\mi(E^r)$, with inverse $I\ff EI$. Similarly, the correspondence 
$F\ff F^l=\ov{\gen}\pr{F}{F}_l$ is a lattice isomorphism between $\mi(E)$ 
and $\mi(E^l)$, with inverse $J\ff JE$. Consequently, there is a lattice 
isomorphism, called the Rieffel correspondence, $\rief:\mi(E^r)\to\mi(E^l)$, 
such that $\rief(I)=\ov{\gen}\pr{EI}{EI}_l$; its inverse is $\rief:\
\mi(E^l)\to\mi(E^r)$ given by $\rief(J)=\ov{\gen}\pr{JE}{JE}_r$. 
\end{prop}
\begin{cor}\label{cor:lattice}
Let $\pi:E\f F$ be a \hm of *-ternary rings between \cts $E$ and $F$. Then
$\big(\ker(\pi)\big)^r=\ker(\pi^r)$. In particular, 
$\pi$ is injective iff $\pi^r$ is injective.
\end{cor}
\ve
\subsection{Morita equivalent partial actions}
\begin{df}\label{df:patring}
Let $E$ be a positive \ct and $\al =(\{ E_t\}_{t\in G},\{ \al_t\}_{t\in G})$  
a set theoretic \pa on $E$, where $E_t\id E$, and $\al_t:E_{t^{-1}}\f E_t$ 
is an isomorphism of \cts, $\forall t \in G$. We say that 
$\al$ is a \pa of $G$ on $E$ if $\{E_t\}_{t\in G}$ is a continuous family 
(\ref{df:cf}) and, 
if $\me^{-1} =\{ (t,x):\ x\in E_{t^{-1}}\}\subseteq G\times E$ with the 
product topology, then the map $\me^{-1}\f E$ such that $(t,x)\ff\al_t(x)$ is 
continuous.   
\end{df}
\begin{prop}\label{prop:apmr}
Let $\al =(\{ E_t\}_{t\in G},\{ \al_t\}_{t\in G})$ be a \pa of the discrete
group $G$ on the \ct $E$, and consider the pairs: 
$\al^l =(\{ E_t^l\}_{t\in G},\{\al_t^l\}_{t\in G})$ and 
$\al^r =(\{ E_t^r\}_{t\in G},\{\al_t^r\}_{t\in G})$.
Then $\al^l$ is a \pa of $G$ on $E^l$, and $\al^r$ is a \pa of $G$ on $E^r$. 
\end{prop}
\begin{proof}
By \ref{prop:lattice} and \ref{prop:z(pi)}, $E_t^r\id E^r$, 
$\al_t^r:E_{t^{-1}}^r\f E_t^r$ is an isomorphism, $\forall t\in G$, and 
$\al_e^r=id_{E^r}$. It follows from \ref{prop:lattice} that, since 
$\al_{st}$ is an extension of $\al_s\al_t$, 
then $\al_{st}^r$ is an extension of $(\al_s\al_t)^r$. 
Now, the domain of $\al_s\al_t$ is 
$E_{t^{-1}}\bigcap E_{t^{-1}s^{-1}}$, and $\al_s\al_t:
E_{t^{-1}}\bigcap E_{t^{-1}s^{-1}}\f E_t\bigcap E_{st}$ is an isomorphism. 
By \ref{prop:lattice}, the domain of  
$(\al_s\al_t)^r$ is $E_{t^{-1}}^r\bigcap E_{t^{-1}s^{-1}}^r$.
Since $\al$ is a \pa, we also have by \ref{prop:lattice} that  
$\al_t^r\big(E_{t^{-1}}^r\bigcap E_r^r\big)=E_t^r\bigcap
E_{tr}^r$. It follows that the domain of $\al_s^r\al_t^r$ 
is $E_{t^{-1}}^r\bigcap E_{t^{-1}s^{-1}}^r$, and hence that  
$\al_s^r\al_t^r=(\al_s\al_t)^r$. Therefore, $\al_{st}^r$ is an extension of 
$\al_s^r\al_t^r$.  
\end{proof}

\begin{rk}\label{rk:adjpa}
It is clear that any \pa $\ga$ on a \ct $E$ is also a \pa on  
$E^*$, the adjoint \ct of $E$ (By $E^*$ we mean the \ct naturally associated 
to the adjoint bimodule of the Hilbert $E^r$-module $E$). Let $\ga$ denote 
this \pa. Then it is easy to see that $(\ga^*)^l=\ga^r$, and $(\ga^*)^r=\ga^l$.
\end{rk}

\begin{ex}\label{ex:rideal}
If $E$ is a right ideal of a \cs $A$, where $G$ acts by an
action $\beta$, then $\al :=\beta\r{E}$ is a \pa on $E$, and we have that 
$\al^r=\beta\r{\ov{\gen}E^*E}$, and $\al^l=\beta\r{\ov{\gen}EE^*}$.
\end{ex}

\begin{df}\label{df:spmr} 
Let $\al =(\{ A_t\}_{t\in G},\{ \al_t\}_{t\in G})$ and  
      $\beta =(\{ B_t\}_{t\in G},\{ \beta_t\}_{t\in G})$ be \pas  
of $G$ on the \css $A$ and $B$ respectively. 
We say that $\al$ is Morita equivalent to $\beta$ if there exists 
a \pa $\ga =(\{ E_t\}_{t\in G},\{ \ga_t\}_{t\in G})$ 
on a positive \ct $E$, such that $\ga^l=\al$, and $\ga^r=\beta$. 
We will denote this relation by $\al\sM\beta$. 
\end{df}
\begin{rk}\label{rk:combes}
In \cite{combes}, Combes defined Morita equivalence of actions. 
When in Definition \ref{df:spmr} $\al$ and $\beta$ are global actions, $\ga$ 
must also be a global action, and therefore our definition of Morita 
equivalence agrees with his in this case. Note that, in fact, $\al$, 
$\beta$ and $\ga$ are global actions if and only if one of them is a 
global action.
\end{rk}
\begin{lem}\label{lem:suerte}
Let $E$ be a \ct and $A=E^r$. If $\{ D_t\}_{t\in G}$ is a continuous 
family of ideals in $A$, and $E_t:=ED_t$, $\forall t\in G$, then 
$\{ E_t\}_{t\in G}$ is a continuous family of ideals in $E$.
\end{lem}
\begin{proof}
Let $U\subseteq E$ be an open set, $G_U=\{ s\in G:\ U\cap E_s\neq\emptyset\}$, 
and $t\in G_U$. Consider $x\in U\cap E_t$. By Cohen--Hewitt, $x=ya$, for some 
$y\in E$, and $a\in D_t$. Since the action of $A$ on $E$ is continuous, 
there exist open sets $V\subseteq E$ and $W\subseteq A$, such that $y\in V$, 
$a\in W$, and $VW\subseteq U$. Now, $a\in W\cap D_t$, and since 
$\{ D_s\}_{s\in G}$ is continuous, the set 
$G_W=\{ s\in G:\ W\cap D_s\neq\emptyset\}$ is open and contains $t$. 
For each $s\in G_W$ take $a_s\in W\cap D_s$. 
Then $xa_s\in E_s\cap VW\subseteq E$, so  
$t\in G_W\subseteq G_U$, and hence $G_U$ is open. 
\end{proof}
\par We will see next that Morita equivalence of \pas is an equivalence 
relation. Recall that if $E$ is a $(A-B)$--Hilbert bimodule and $F$ is a 
$(B-C)$--Hilbert bimodule, their inner tensor product is the 
$(A-C)$--Hilbert bimodule 
$\bts{E}{B}{F}$ constructed as follows: let  
$\bats{E}{F}$ their algebraic tensor product, and consider on $\bats{E}{F}$ 
the unique $C$--sesquilinear map $\pr{\cdot}{\cdot}_r'$ such that 
$\pr{x_1\odot y_1}{x_2\odot y_2}=\pr{y_1}{\pr{x_1}{x_2}_By_2}_C$, where 
$\pr{\cdot}{\cdot}_C$ is the $C$--inner product on $F$, and 
$\pr{\cdot}{\cdot}_B$ is the $B$--inner product on $E$. This sesquilinear map  
is a semi--inner product, that defines an inner product on the quotient 
$(\bats{E}{F})/N$, where $N=\{ z\in\bats{E}{F}:\, \pr{z}{z}_r'=0\}
=\gen\{ xb\odot y-x\odot by:\, x\in E, y\in F, b\in B\}$. Then, 
$\bts{E}{B}{F}$ is the completion of $(\bats{E}{F})/N$ with respect to this 
inner product (see for instance \cite{l} for details). We will denote 
by $x\otimes y$ the projection of $x\odot y$ on $\bts{E}{B}{F}$.

\begin{lem}\label{lem:innermaps}
Let $\mu_i:E_i\to F_i$ be \cts \hms for $i=1,2$.
Suppose that $A$, $B$ and $C$ are \css such that $E_1^l,F_1^l\id A$, 
$E_2^r,F_2^r\id C$, and $E_1^r=E_2^l\id B$. Suppose, moreover, that 
$\mu_1^r=\mu_2^l$. Then there exists a unique \hm of \cts 
$\mu_1\otimes_B\mu_2:\bts{E_1}{B}{E_2}\f\bts{F_1}{B}{F_2}$ such that  
$(\mu_1\otimes_B\mu_2)(x_1\otimes x_2)=\mu_1(x_1)\otimes\mu_2(x_2)$, 
$\forall x_1\in E_1$, $x_2\in E_2$. If $\mu_1$ and $\mu_2$ are isomorphisms, 
then $\mu_1\otimes_B\mu_2$ also is an isomorphism. Moreover, we have that 
$(\mu_1\otimes_B\mu_2)^l=\mu_1^l$, and $(\mu_1\otimes_B\mu_2)^r=\mu_2^r$.
\end{lem}
\begin{proof}
Let $\mu_1\odot\mu_2:\bats{E_1}{E_2}\f\bats{F_1}{F_2}$ be the unique
linear transformation such that $x_1\odot x_2\ff\mu_1(x_1)\odot\mu_2(x_2)$, 
$\forall x_1\in E_1$, $x_2\in E_2$; it is a \hm of *--ternary rings.
Let $\pr{\cdot}{\cdot}$ and 
$[\cdot ,\cdot ]$ be the corresponding $C$--pre--inner products on  
$\bats{E_1}{E_2}$ and $\bats{F_1}{F_2}$ respectively. 
Pick $z, z'\in\bats{E_1}{E_2}$, $z=\sum_ix_i\odot y_i$, 
$z'=\sum_jx_j'\odot y_j'$. Since 
$[(\mu_1\odot\mu_2)(z),(\mu_1\odot\mu_2)(z')]
                =\sum_{i,j}\pr{\mu_2(y_i)}
                 {\mu_1^r(\pr{x_i}{x_j'}_B^{E_1})\mu_2(y_j')}_C^{F_2}$, 
$\mu_2^r(\pr{z}{z'})
                =\sum_{i,j}\pr{\mu_2(y_i)}
                 {\mu_2^l(\pr{x_i}{x_j'}_B^{E_1})\mu_2(y_j')}_C^{F_2}$, 
and $\mu_1^r=\mu_2^l$, we conclude that  
$[(\mu_1\odot\mu_2)(z),(\mu_1\odot\mu_2)(z')]=\mu_2^r(\pr{z}{z'})$.
By taking $z=z'$ and computing norms, we have:  
$\norm{(\mu_1\odot)\mu_2(z)}^2=\norm{\mu_2^r(\pr{z}{z})}
 \leq\norm{z}^2.$ 
Thus, $\mu_1\odot\mu_2$ factors through the quotient, where it is a 
contraction, and hence extends by continuity to a \hm of \cts 
$\mu_1\otimes_B\mu_2:\bts{E_1}{B}{E_2}\f\bts{F_1}{B}{F_2}$. We have, 
$\forall z,z'\in\bts{E_1}{B}{E_2}$:
$[(\mu_1\otimes_B\mu_2)(z),(\mu_1\otimes_B\mu_2)(z')]
    =\mu_2^r(\pr{z}{z'})$, and therefore $(\mu_1\otimes_B\mu_2)^r=\mu_2^r$.
Similarly, $(\mu_1\otimes_B\mu_2)^l=\mu_1^l$. Finally, if $\mu_1$, $\mu_2$ 
are isomorphisms, we apply the first part of the proof to the maps 
$\mu_1^{-1}$ e $\mu_2^{-1}$, and we note that  
$id_{E_1}\otimes id_{E_2}=id_{\bts{E_1}{B}{E_2}}$,  
$id_{F_1}\otimes id_{F_2}=id_{\bts{F_1}{B}{F_2}}$. 
\end{proof}

\begin{prop}\label{prop:pareleq}
Morita equivalence of partial actions is an equivalence relation.  
\end{prop}
\begin{proof}
The reflexive and symmetric properties are immediately verified 
(see Remark \ref{rk:adjpa}).
\par Suppose now that $\al =(\{ A_t\},\{\al_t\})$ is a \pa  
of $G$ on $A$, $\beta =(\{ B_t\},\{\beta_t\})$ is a \pa of 
$G$ on $B$, and $\ga =(\{ C_t\},\{\ga_t\})$ is a \pa  
of $G$ on $C$, such that $\al\sM\beta$ through the \pa  
$\mu =(\{ E_t\} ,\{\mu_t\})$ of $G$ on the \ct $E$, and 
$\beta\sM\ga$ through the \pa $\nu =(\{ F_t\} ,\{\nu_t\})$ of  
$G$ on the \ct $F$. Consider the family 
$\mu\otimes_B\nu :=(\{\bts{E_t}{B}{F_t}\} ,\{\mu_t\otimes_B\nu_t\})$. 
Since $\mu_{rs}$ extends $\mu_r\mu_s$ and $\nu_{rs}$ extends $\nu_r\nu_s$, 
it follows that $\mu_{rs}\otimes_B\nu_{rs}$ extends 
$(\mu_r\otimes_B\nu_r)(\mu_s\otimes_B\nu_s)$. It is clear that 
$E_e\bigotimes_BF_e=E\bigotimes_BF$, $(\mu\otimes_B\nu)_e=
id_{\bts{E}{B}{F}}$. 
On the other hand: $(\bts{E_t}{B}{F_t})^r=C_t$:\\
$\ov{\gen}\pr{\bts{E_t}{B}{F_t}}{\bts{E_t}{B}{F_t}}_C
  =\ov{\gen}\pr{F_t}{\pr{E_t}{E_t}_BF_t}_C
  =\ov{\gen}\pr{F_t}{B_tF_t}_C
  =\ov{\gen}\pr{F_t}{F_t}_C
  =C_t.$\\
Similarly, $(\bts{E_t}{B}{F_t})^l=A_t$. 
Finally, by \ref{lem:innermaps}, every $\mu_t\otimes_B\nu_t:
\bts{E_t^{-1}}{B}{F_t^{-1}}\f\bts{E_t}{B}{F_t}$ is an isomorphism, and 
$(\mu_t\otimes_B\nu_t)^l=\mu_t^l=\al_t$, 
$(\mu_t\otimes_B\nu_t)^r=\nu_t^r=\ga_t$, so 
$\mu\otimes_B\nu$ is a set theoretic \pa on $\bts{E}{B}{F}$, and 
$(\mu\otimes_B\nu)^l=\al$, $(\mu\otimes_B\nu)^r=\ga$. 
\par It remains to show that $\bts{\mu}{B}{\nu}$ is continuous. 
First, note that the family $\{\bts{E_t}{B}{F_t}\}_{t\in G}$ is 
continuous by \ref{lem:suerte}, because $\ga$ is a \pa. Let  
$\bts{\me^{-1}}{B}{\mf^{-1}}=\{(t,z):\ z\in\bts{E_{t^{-1}}}{B}{F_{t^{-1}}}\}$. 
Note that if $f\in C_c(\me^{-1})$, 
$g\in C_c(\mf^{-1})$, then $f\oslash_Bg:G\f\bts{\me^{-1}}{B}{\mf^{-1}}$ such 
that $(f\oslash_Bg)(t)=\big(t,f(t)\otimes g(t)\big)$ is a continuous section 
of the Banach bundle $\bts{\me^{-1}}{B}{\mf^{-1}}$, and that for each 
$t\in G$, $\gen\{ (f\oslash_Bg)(t):\ f\in C_c(\me^{-1}), g\in 
C_c(\mf^{-1})\}$ is dense in $(\bts{\me^{-1}}{B}{\mf^{-1}})_t$. On the other 
hand, the map $G\f G\times (\bts{E}{B}{F})$  
such that $t\ff \big(t,\mu_t\big(f(t)\big)\otimes\nu_t\big(g(t)\big)\big)$ is  
continuous, $\forall f\in C_c(\me^{-1})$, $g\in C_c(\mf^{-1})$. So we 
conclude, by \cite{fd}, II-14.6 and II-13.16, that the application 
$\bts{\me^{-1}}{B}{\mf^{-1}}
\f G\times (\bts{E}{B}{F})$ such that $(t,x)\ff (t,(\mu_t\otimes\nu_t)(x))$ 
is a continuous \hm of Banach bundles, and therefore 
$\mu\otimes_B\nu$ is a continuous \pa. Thus $\al\sM\ga$.
\end{proof}

\begin{prop}\label{prop:equivprods}
Let $\al =(\{ A_t\}_{t\in G},\{ \al_t\}_{t\in G})$ and  
$\beta =(\{ B_t\}_{t\in G},\{ \beta_t\}_{t\in G})$ be \pas of  
$G$ on the \css $A$ and $B$ respectively. If $\al$ and $\beta$ are 
Morita equivalent, then $A\sd{\al ,r}G\sM B\sd{\beta ,r}G$. 
\end{prop}
\begin{proof}
Suppose that $\al$ and $\beta$ are Morita equivalent through a \pa 
$\ga =(\{ E_t\},\{\ga_t\})$; say that $\ga^l =\al$, 
and $\ga^r=\beta$. Consider the full left Hilbert $A$--module $A\bigoplus E$,
where $\pr{(a_1,e_1)}{(a_2,e_2)}=a_1a_2^*+\pr{e_1}{e_2}_l$. Then 
$A\bigoplus E$ establishes a Morita equivalence between $A$ and $\link{E}$, 
the linking algebra \label{vinc} of $E$ (see \cite{rw} for instance). By 
\ref{prop:lattice}, every ideal $A_t$ corresponds to an ideal 
$\rief(A_t)\id\link{E}$, and it is easy to see that 
$\rief(A_t)=\link{E_t}$, $\forall t\in G$. Since 
$\link{E_t}=\begin{pmatrix} A_t & E_t\\ E^*_t & B_t\end{pmatrix}$, and  
$(A_t)$, $(E_t)$, $(E_t^*)$ and $(B_t)$ are continuous families, so it is 
the family $(\link{E_t})_{t\in G}$. Let us define now 
$\link{\ga_t}:\link{E_{t^{-1}}}\f\link{E_t}$ by  
$\link{\ga_t}:\begin{pmatrix} a_{t^{-1}} & x_{t^{-1}}\\  
y_{t^{-1}}& b_{t^{-1}}\end{pmatrix}
\ff\begin{pmatrix} \al_t(a_{t^{-1}}) & \ga_t(x_{t^{-1}})\\  
\ga_t(y_{t^{-1}})& \beta_t(b_{t^{-1}})\end{pmatrix}.$
Since $\ga_t(a_{t^-1}x_{t^{-1}})=
\al_t(a_{t^{-1}})\ga_t(x_{t^{-1}})$ and $\ga_t(x_{t^{-1}}b_{t^{-1}})=
\ga_t(x_{t^{-1}})\beta_t(b_{t^{-1}})$, $\forall x_{t^{-1}}\in E_{t^{-1}}$, 
$a_{t^{-1}}\in A_{t^{-1}}$, $b_{t^{-1}}\in B_{t^{-1}}$, we have that  
$\link{\ga}=(\{ \link{E_t}\}_{t\in G},\{\link{\ga_t}\}_{t\in G})$ is a   
\pa of $G$ on $\link{E}$. We call $\link{\ga}$ the \textit{linking \pa} of 
$\ga$. Observe that if $\ga_1$ is the restriction of 
$\link{\ga}$ to $\begin{pmatrix}A&E\\0&0\end{pmatrix}$, then $\ga_1^l=\al$ and 
$\ga_1^r=\link{\ga}$. Similarly, it is enough to restrict $\link{\ga}$ to  
$\begin{pmatrix}0&0\\E^*&B\end{pmatrix}$ to see that it is also Morita 
equivalent to $\beta$. 
\par The considerations above show that we may assume that $A=pBp$, for some 
projection $p\in M(B)$, and that $\al$ and $\ga$ are the  
restrictions of $\beta$ to $pBp=A$ and to $E=pB$ respectively. Let  
now $\ma$ and $\mb$ be the Fell bundles corresponding to $\al$ and $\beta$ 
respectively, and let $\me =(\{ t\}\times E_t)_{t\in G}$. 
We have that $\ma\subseteq\me\subseteq\mb$, $\ma$ is a  
sub--Fell bundle of $\mb$, and that $\me$ is a sub--Banach bundle of $\mb$. 
On the other hand, if $(r,a_r)\in\ma$, $(s,x_s),\, (t,y_t)\in\me$, and  
$(u,b_u)\in\mb$, we have:
\bi
\item $(r,a_r)(s,x_s)=(rs,\beta_r\big(\beta_r^{-1}(a_r)x_s\big))\in
      (rs,\beta_r(A_{r^{-1}}E_s))\in\me$.
\item $(s,x_s)(u,b_u)=(su,\beta_s\big(\beta_s^{-1}(x_s)b_u\big))\in
       (su,\beta_s(E_{s^{-1}}B_u))\in\me$.
\item $(s,x_s)(t,y_t)^*=(s,x_s)(t^{-1},\beta_{t^{-1}}(y_t^*))
      =(st^{-1},\beta_s\big(\beta_{s^{-1}}(x_s)\beta_{t^{-1}}(y_t^*)\big))
      $, that belongs to \\ $(st^{-1},\beta_s(A_{s^{-1}}\cap A_{t^{-1}}))
      \subseteq (st^{-1},A_s\cap A_{st^{-1}})\in\ma$.
\item $(s,x_s)^*(t,y_t)=(s^{-1},\beta_{s^{-1}}(x_s^*))(t,y_t)
      =(s^{-1}t,\beta_{s^{-1}}(x_s^*y_t))$. Now, for all $t\in G$,
      we have that  
      $\ov{\gen}\{\beta_{s^{-1}}(x_s^*y_t):\, s\in G, x_s\in E_s, y_t\in E_t\}
      =B_t$: the left member of this equality contains 
      $\ov{\gen}E_e^*E_t=\ov{\gen}E^*E_t\supseteq\ov{\gen}E_t^*E_t=B_t$.      
\ei 
Therefore, we may apply Corollary \ref{cor:equimorl1}, 
from where we conclude that $A\sd{\al,r}G\sM B\sd{\beta,r}G$. 
\end{proof}
\par This is a good point to introduce the notion of \mea. 
\begin{df}\label{df:envupmor}
Let $\al$ be a \pa of $G$ on the \cs $A$. We say that  
a continuous action $\beta$ of $G$ on a \cs $B$ is a \mea  
of $\al$, if there exists an ideal $I\id B$ such that $[\beta(I)]$ is dense in 
$B$ and $\al\sM\beta\r{I}$. In other words: $\beta$ is the enveloping action 
of a \pa that is Morita equivalent to $\al$. 
\end{df}
\par We close the section with a result that is similar to 
Theorem \ref{thm:equivmor}. 
\begin{prop}\label{prop:equivmor2}
Let $(\al,A)$ be a \pa of $G$, and assume that $(\beta,B)$ 
is a \mea of $\al$. Then: 
$ A\sd{\al,r}G\sM B\sd{\beta,r}G.$
\end{prop}
\begin{proof}
Since Morita equivalence of \css is transitive, the proof follows 
immediately by combining Proposition \ref{prop:equivprods} above with Theorem 
\ref{thm:equivmor}. 
\end{proof}
\section{$\mathrm{C^*}$-algebras of kernels associated with a Fell bundle}
\label{sec:nucs}
\par In the present section we study two \css, $\nuc{\mb}$ and $\nucr{\mb}$, 
that are naturally associated with a given Fell bundle $\mb$ over the group 
$G$. Both of them are Hausdorff completions of a certain *-algebra of integral 
operators. The first one is defined by a universal property; 
the second one is a concrete \cs of adjointable operators on a Hilbert module. Indeed, $\nucr{\mb}$ is the image of $\nuc{\mb}$ under a certain natural representation on $L^2(\mb)$, which is faithful when $\mb$ is saturated. There is a natural action of the group on $\nuc{\mb}$. If $\mb$ is the \fb   
of a \pa, we will see in the following section that this natural 
action on $\nuc{\mb}$ is a \mea of the given \pa.   
\par Let $\mb =(B_t)_{t\in G}$ be a Fell bundle.
Consider a continuous function $k:G\times G\f\mb$ of compact support, 
such that $k(r,s)\in B_{rs^{-1}}$, $\forall r,s\in G$. Such a function  
will be called a {\em kernel of compact support} associated with $\mb$.
The linear space of kernels of compact support associated with $\mb$ 
will be denoted by $\nucc{\mb}$ 
We will see later that any $k\in\nucc{\mb}$ may be seen as an integral 
operator, which justifies this terminology. 
\begin{prop}\label{prop:nucc}
$\nucc{\mb}$ is a normed *-algebra with the involution  
$k^*(r,s)\! =\! k(s,r)^*,\forall k\in \nucc{\mb}$, the product  
$k_1*k_2(r,s)\!=\!\int_G\!k_1(r,t)k_2(t,s)dt,\, 
\forall k_1,k_2\in\nucc{\mb},$ and the norm $\norm{k}_2
=\left(\int_{G^2}\!\norm{k(r,s)}^2 dr ds\right)^{1/2}$\!.
\end{prop} 
\begin{proof}
Let $\nu:G\times G\f G$ be such that $\nu(r,s)=rs^{-1}$, and let 
$\mb_{\nu}$ be the retraction of $\mb$ with respect to 
$\nu$ (\cite{fd}, II-13.3). Then  
$\mb_{\nu}$ is a Banach bundle over $G\times G$, and the fiber of $\mb_{\nu}$ 
over $(r,s)$ is $(r,s,B_{rs^{-1}})$, which we may naturally identify with 
$B_{rs^{-1}}$. Therefore, $\nucc{\mb}=C_c(\mb_{\nu})$ as a linear space.   
\par Consider now the map $\mu:G\times G\times G\f\mb_{\nu}$ 
given by $\mu(t,r,s)=\big(r,s,k_1(r,t)k_2(t,s)\big)$. We have that $\mu$ is 
continuous and has compact support, and that $\mu(t,r,s)\in 
(\mb_{\nu})_{(r,s)}$, $\forall t,r,s\in G$. Thus, we may apply 
\cite{fd}, II-15.19, from where we conclude that the map 
$(r,s)\ff\int_G\mu(t,r,s)dt$ is a continuous section of compact support
of $\mb_{\nu}$. In other words, $k_1*k_2\in\nucc{\mb}$. 
\par As for $k^*$, we have that $\text{supp}(k^*)$ is compact, and $k^*(r,s)=
k(s,r)^*\in B_{sr^{-1}}^*=B_{rs^{-1}}$. Consequently, $k^*\in\nucc{\mb}$.
Routine computations show that $(k_1*k_2)*k_3=k_1*(k_2*k_3)$, 
$(k_1*k_2)*=k_2^**k_1^*$, so $\nucc{\mb}$ is a *-algebra. It is also 
immediate that $\norm{k^*}_2=\norm{k}_2$.
Finally:
\[\norm{k_1*k_2}_2^2
\leq\int_G\int_G\left[\int_G\norm{k_1(r,t)k_2(t,s)}dt\right]^2drds
\leq\int_G\int_G\left[\int_G\norm{k_1(r,u)}^2du\int_G\norm{k_2(v,s)}^2dv\right]
    drds\]%
\end{proof}
\begin{prop}\label{prop:bim}
Let $\mb=(B_t)_{t\in G}$ be a Fell bundle, and consider the action 
$\nucc{\mb}\times C_c(\mb)\f C_c(\mb)$ given by 
$k\cdot\xi\r{r}=\int_Gk(r,s)\xi(s)ds$, $\forall k\in\nucc{\mb}$, 
$\xi\in C_c(\mb)$, and $r\in G$. With this action, $C_c(\mb)$ is 
a $(\nucc{\mb}-B_e)$--bimodule. Moreover, if 
$\pr{\cdot}{\cdot}_l: C_c(\mb)\times C_c(\mb)\f \nucc{\mb}$ is such that 
$\pr{\xi}{\eta}_l\r{(r,s)}=\xi(r)\eta(s)^*$, $\forall \xi,\eta\in 
C_c(\mb)$, $r,s\in G$, we have:
\be
 \item $\pr{\xi_1}{\xi_2}_l\xi_3
       =\xi_1\pr{\xi_2}{\xi_3}_r$, $\forall \xi_1,\xi_2,\xi_3\in C_c(\mb)$.
 \item $\pr{k\cdot \xi}{\eta}_l=k*\pr{\xi}{\eta}_l$, $\forall k\in\nucc{\mb}$, 
       $\xi$, $\eta\in C_c(\mb)$. 
 \item $\pr{\xi}{\eta}_l^*=\pr{\eta}{\xi}_l$, $\forall \xi,\eta\in C_c(\mb)$.
 \item $\pr{\xi}{\xi}_l=0\iff \xi=0$.
 \item $\pr{k\xi}{\eta}_r=\pr{\xi}{k^*\eta}_r$, 
       $\forall k\in\nucc{\mb}$, $\xi$, $\eta\in C_c(\mb)$, where 
       $\pr{\cdot}{\cdot}_r:C_c(\mb)\times C_c(\mb)\f B_e$ is the right inner 
       product, that is $\pr{\xi}{\eta}_r=\int_G\xi(s)^*\eta(s)ds$.
\ee
\end{prop}
\begin{proof}
All these properties are easy to verify. As an example we prove 5., 
and leave 1.--4. to the reader. 
If $k\in\nucc{\mb}$, $\xi$, $\eta\in C_c(\mb)$, 
\[ \pr{k\xi}{\eta}_r\!
   =\!\int_G\bigg[\int_Gk(r,s)\xi(s)ds\bigg]^*\eta(r)dr\!
   =\!\int_G\int_G\xi(s)^*k(r,s)^*\eta(r)drds\!
   =\!\int_G\xi(s)^*(k^*\eta)(s)ds\!
   =\!\pr{\xi}{k^*\eta}_r. \] 
\end{proof}
\par We define $I_c(\mb):=\gen\pr{C_c(\mb)}{C_c(\mb)}_l$, where 
$\pr{\cdot}{\cdot}_l$ is the map defined in \ref{prop:bim}. Clearly, 
$I_c(\mb)$ is a two--sided *-ideal of $\nucc{\mb}$.
\par Let $E=L^2(\mb)$, and let $[\cdot,\cdot]:E\times E\f \k{E}$ be 
the corresponding left inner product (If we think of $E$ as a positive \ct, 
then $\k{E}$ is nothing but $E^l$). Note that there is a natural 
injective *-\hm $I_c(\mb)\inc\k{E}$, the only one such that $\pr{\xi}{\eta}_l
=[\xi,\eta ]$, $\forall \xi,\eta \in C_c(\mb)$. Indeed, if  
$k=\sum_i\pr{\xi_i}{\eta_i}_l=0$, then $k\zeta =0$, 
$\forall\zeta\in C_c(\mb)$, and since $k\zeta =\sum_i[\xi_i,\eta_i]\zeta$,  
and $C_c(\mb)$ is dense in $E$, we see that $\sum_i[\xi_i,\eta_i]=0\in
\k{E}$. On the other hand, since $\gen [C_c(\mb),C_c(\mb)]$ is dense in 
$\k{E}$, we have that $\k{E}$ is a $C^*$-completion of the 
*-algebra $I_c(\mb)$. We will see later that, when the Fell bundle $\mb$ is saturated, this inclusion extends to an inclusion $\Om:\nucc{\mb}\f \adj{E}$ (Theorem~\ref{thm:picpi2}).     
\par Note that, as a Banach space, the completion $HS(\mb)$ of 
$(\nucc{\mb},\norm{\cdot}_2)$ agrees with $\mathfrak{L}^2(\mb_{\nu})$ 
(see \cite{fd}, II-15.7--15.9). When $\mb=(B_t)_{t\in G}$ is the trivial Fell 
bundle over $G$ with constant fiber $\C$ (that is: the \fb associated with 
the trivial action of $G$ on $\C$), 
then $\mathfrak{L}^2(\mb_{\nu})=L^2(G\times G)$ is naturally identified 
with the Hilbert--Schmidt operators on $L^2(G)$. Hence we may think of 
$HS(\mb)$ as the algebra of ``Hilbert--Schmidt operators'' on $L^2(\mb)$.
\begin{df}\label{df:nuc}
The universal \cs of $HS(\mb)$ will be called the \cs of kernels of 
the Fell bundle $\mb$, and will be denoted by $\nuc{\mb}$. The closure of 
$I_c(\mb)$ in $\nuc{\mb}$ will be denoted by $I(\mb)$.
\end{df}
\ve
\subsection{Natural action on the kernels}\label{subsection:acnucs}
There is a natural action of $G$ on $\nucc{\mb}$:
$\beta:G\times\nucc{\mb}\f\nucc{\mb}$ such that  
$\beta_t(k)(r,s)=\Del(t)k(rt,st)$, where $\Del$ is the modular function on $G$.
We have:
\[\beta_t(k_1*k_2)\r{(r,s)}
=\Del(t)^2\int_Gk_1(rt,vt)k_2(vt,st)dv
=\int_G\beta_t(k_1)(r,v)\beta_t(k_2)(v,s)dv
=\beta_t(k_1)*\beta_t(k_2)\r{(r,s)}\]
\[\big(\beta_t(k^*)\big)\r{(r,s)}
=\Del(t)k^*(rt,st)
=\Del(t)k(st,rt)^*
=\beta_t(k)(s,r)^*
=\beta_t(k)^*\r{(r,s)}.\]
This action may be extended to $HS(\mb)$ and hence to $\nuc{\mb}$: 
doing $u=rt$, $v=st$ in the integral below:
\[\norm{\beta_t(k)}_2^2
=\int_G\int_G\Del(t)^2\norm{k(rt,st)}^2drds
=\int_G\int_G\Del(t)^{-1}\Del(t)^{-1}\Del(t)^2\norm{k(u,v)}^2dudv
=\norm{k}_2^2.\]
Note that $\beta$ is a continuous action on $\nucc{\mb}$ with the \ilt  
(recall from \ref{prop:nucc} that $\nucc{\mb}=C_c(\mb_\nu)$), and therefore 
is continuous on $HS(\mb)$. Since $\nuc{\mb}$ is the universal \cs of 
$HS(\mb)$, $\beta$ also extends to a continuous action on $\nuc{\mb}$. 
All these actions will be denoted by $\beta$.
\begin{lem}\label{lem:bolu}
Let $B$ be a Banach bundle over the \lc space $X$. Suppose that  
$\Theta\subseteq C_c(X)$ is dense in $C_c(X)$ in the \ilt, and that 
$F\subseteq C_c(B)$ is a linear subspace such that $\Theta F\subseteq F$. 
Then the closure of $F$ in the \ilt is the space: 
$\{ g\in C_c(B):\ g(x)\in\ov{F(x)}, \, \forall x\in X\}$. 
In particular, if $\ov{F(x)}=B_x$, $\forall x\in X$, then $F$ is dense in 
$C_c(B)$.
\end{lem}
\begin{proof}
For $x\in X$, the map $e_x:C_c(B)\f B_x$ such that $f\ff f(x)$ 
is a continuous linear map in the inductive limit topology . 
Therefore $e_x(\ov{F})\subseteq\ov{F(x)}$.
Conversely, suppose that $g\in C_c(B)$ is such that $g(x)\in\ov{F(x)}$, 
$\forall x\in X$. Since $X$ is \lc, there exists a compact subset $K$
of $X$ whose interior contains $\text{supp}(g)$. Now, given 
$\epsilon >0$ and $x\in K$, there exists $f_x\in F$ such that 
$\norm{g(x)-f_x(x)}<\epsilon$. Since $g,f_x$ and the norm on $B$ are  
continuous maps, there exists a precompact open neighborhood $V_x$ of $x$, 
such that $\norm{g(y)-f_x(y)}<\epsilon$, $\forall y\in V_x$. The family 
$(V_x)_{x\in K}$ is an open covering of $K$, so it has a finite 
subcovering $V_{x_1},\ldots ,V_{x_m}$. Let $(\psi_j)_1^m$
be a partition of unity of $K$ subordinated to $(V_{x_j})$, 
and let $f=\sum_{j=1}^m\psi_jf_{x_j}$. We have that $f\in C_c(B)$ and 
$\norm{g-f}_\infty<\epsilon$. Therefore, it is enough to show that 
$f\in\ov{F}$. For this, it is sufficient to show that $\psi f'\in\ov{F}$, 
$\forall \psi\in C_c(G)$ and  
$f'\in F$. But, since $\Theta$ is dense in $C_c(G)$, there exists 
$\psi'\in\Theta$ such that $\norm{\psi-\psi'}_{\infty}\leq\varepsilon$, 
for a given $\varepsilon$, and hence 
$\norm{\psi f'-\psi'f'}_{\infty}\leq\varepsilon\norm{f'}_{\infty}$. 
Since $\Theta F\subseteq F$, we have that $\psi'f'\in F$, and therefore 
$\psi f'\in \ov{F}$. 
\end{proof}

\begin{prop}\label{prop:envnucc}
The linear $\beta$--orbit of $I_c(\mb)$ is dense in $\nucc{\mb}$ in the \ilt.
\end{prop}
\begin{proof}
Recall that $\nucc{\mb}=C_c(\mb_{\nu})$, where $\mb_{\nu}$ is the  
retraction of $\mb$ with respect to the map $\nu:G\times G\f G$ 
such that $\nu(r,s)=rs^{-1}$. By Lemma \ref{lem:bolu}, it suffices  
to show that, if $Z$ is the linear $\beta$--orbit of $I_c(\mb)$, then: 
1. $Z(r,s)$ is dense in $B_{rs^{-1}}$, $\forall r,s\in G$ and  
2. $\forall \varphi_i,\psi_i\in C_c(G)$, and $\forall \zeta\in Z$, the  
function $(r,s)\ff\sum_i\varphi_i(r)\psi_i(s)\zeta(r,s)$ belongs to $Z$. 
\par 1. If $b\in B_{rs^{-1}}$, by Cohen--Hewitt there exist $b_1\in B_e$, 
       $b_2\in B_{rs^{-1}}$ such that $b_1b_2=b$. There also exist sections 
       $\xi,\eta\in C_c(\mb)$ such that $\xi(e)=\Del(r)b_1$, 
       $\eta(sr^{-1})=b_2^*$. Thus we have: 
       $\beta_{r^{-1}}(\pr{\xi}{\eta}_l)\r{(r,s)}
       =\Del(r)^{-1}\xi(rr^{-1})\eta(sr^{-1})^*
       =\Del(r)^{-1}\Del(r)b_1(b_2^*)^*=b$. 
\par 2. Let $\varphi$, $\psi\in C_c(G)$, $\xi$, $\eta\in C_c(\mb)$, 
       $r,s,t\in G$. If $\phi:G\f\C$, let $\phi_t:G\f\C$ be given by  
       $\phi_t(s)=\phi(st)$. Then: 
       $(\varphi\otimes\psi)\big(\beta_t(\pr{\xi}{\eta}_l)\big)\r{(r,s)}
       =\varphi(r)\psi(s)\Del(t)\xi(rt)\eta(st)^*
       =\varphi_{t^{-1}}(rt)\psi_{t^{-1}}(st)\Del(t)\xi(rt)\eta(st)^*
       =\beta_t\big(\pr{\varphi_{t^{-1}}\xi}
                        {\bar{\psi}_{t^{-1}}\eta}_l\big)\r{(r,s)},$
       and therefore 
       $(\varphi\otimes\psi)\big(\beta_t(\pr{\xi}{\eta}_l)\big)\in Z$, 
       because $\varphi_{t^{-1}}\xi$, $\bar{\psi}_{t^{-1}}\eta\in C_c(\mb)$.
\end{proof}

\par Recall that a Fell bundle $\mb$ is called \textit{saturated} if for all 
$r,s\in G$ we have that $\gen B_rB_s$ is dense in $B_{rs}$.  
The ideal $I_c(\mb)$ measures the level of saturation of the bundle:

\begin{prop}\label{prop:sat}
$\mb$ is saturated if and only if $I_c(\mb)$ is dense in $\nucc{\mb}$ in
the \ilt.
\end{prop}
\begin{proof}
First Suppose that $\mb$ is saturated. If  
$\varphi$, $\psi\in C_c(G)$, $\xi$, $\eta\in C_c(\mb)$, then  
$(\varphi\otimes\psi)\pr{\xi}{\eta}_l=\pr{\varphi\xi}{\bar{\psi}\eta}_l
\in I_c(\mb)$. On the other hand, since $\mb$ is saturated, given  
$a\in B_{rs^{-1}}$ and $\epsilon >0$, there exist $b_1,\ldots ,b_n\in B_r$, 
$c_1,\ldots ,c_n\in B_s$, such that $\norm{b-\sum_{i=1}^nb_ic_i^*}<\epsilon$. 
Moreover, there exist sections $\xi_1,\ldots,\xi_n,\eta_1,\ldots ,\eta_n
\in C_c(\mb)$ such that $\xi_j(r)=b_j$ and $\eta_j(s)=c_j$, 
$\forall j=1,\ldots,n$. Therefore: 
$\norm{b-\sum_{i=1}^n\pr{\xi_i}{\eta_i}_l\r{(r,s)}}
=\norm{b-\sum_{i=1}^n\xi_i(r)\eta_i(s)^*}
=\norm{b-\sum_{i=1}^nb_ic_i^*}<\epsilon$. It follows from \ref{lem:bolu} that 
$I_c(\mb)$ is dense in $\nucc{\mb}$ in the \ilt.
\par Conversely, assume that $\mb$ is not saturated: there exist $r,s\in G$ 
such that $\ov{\gen}B_rB_{s^{-1}}\neq B_{rs^{-1}}$. Let $b\in B_{rs^{-1}}$ be 
such that $b\notin\ov{\gen}B_rB_{s^{-1}}$, and let $d$ be the distance  
from $b$ to $\ov{\gen}B_rB_{s^{-1}}$. There exists a continuous section of 
compact support of $\mb_{\nu}$ that takes the value $b$ in $(r,s)$.  
In other words, there exists $k\in\nucc{\mb}$ such that 
$k(r,s)=b$. Now, if $\xi_i$, $\eta_i\in C_c(\mb)$:
$\norm{k-\sum_i\pr{\xi_i}{\eta_i}_l}_{\infty}
\geq\norm{k(r,s)-\sum_i\pr{\xi_i}{\eta_i}_l(r,s)}
=\norm{b-\sum_i\xi_i(r)\eta(s)^*}
\geq d,$
and therefore $k\notin \ov{I_c(\mb)}$, because $d>0$.    
\end{proof}
\ve 
\subsection{The reduced \cs of kernels of a Fell bundle}
\par Our next goal will be to show that the action of $\nucc{\mb}$ 
on $C_c(\mb)$ extends to all of $E=L^2(\mb)$ and, from this, that
there is a *-\hm $\Om:\nuc{\mb}\f\adj{E}$. This homomorphism is
injective if $\mb$ is saturated. Besides, we will show that every
representation $\pi$ of $\mb$ on a Hilbert space $H$ induces a
representation of $\nuc{\mb}$ on $L^2(G,H)$ which is faithful whenever
$\pi|_{B_e}$ is injective. To do this, we will represent represent $E$
as a ternary ring of operators, and use the uniqueness of the
enveloping action.    
\begin{prop}\label{prop:kenv}
   $(\Bbbk(\mb),\beta)$ is the enveloping action of $(I(\mb),\beta|_{I(\mb)})$.  
\end{prop}
\begin{proof}
This is a direct consequence of Proposition~\ref{prop:envnucc}.
\end{proof}

\begin{lem}\label{lem:niguais}
Let $\mathfrak{c}:=\ov{\{\sum_{i=1}^nk_i*k_i^*:\, n\geq 1,\, k_i\in 
I_c(\mb)\}}\subseteq\nucc{\mb}$, where the closure is taken in the \ilt.
Then, if $\xi\in C_c(\mb)$, we have that 
$\pr{\xi}{\xi}_l\in\mathfrak{c}\cap I_c(\mb)$.  
\end{lem}
\begin{proof}
Let $\xi\in C_c(\mb)$. Since $C_0(\mb)$ is a non--degenerate right Banach 
module over $B_e$, the Cohen--Hewitt theorem implies that there exist 
$\eta\in C_c(\mb)$, and $b\in B_e$, such that  
$\xi=\eta b$. Let $\zeta\in C_c(\mb)$ such that $\zeta(e)=b^*$. The function  
$G\f B_e$ such that $t\ff \zeta(t)^*\zeta(t)-bb^*$ is continuous and vanishes  
at $t=e$. Thus, given $\epsilon >0$, there exists a neighborhood 
$V_{\epsilon}$ of $e$, such that if $t\in V_{\epsilon}$, then 
$\norm{\zeta(t)^*\zeta(t)-bb^*}<\epsilon$. 
Let $(f_V)_{V\in\mv}$ be an approximate unit of $L^1(G)$ as in Lemma 
\ref{lem:lema}, and for each $V\in\mv$ let $\zeta_V=f_V^{1/2}\zeta$. Then  
$\zeta_V\in C_c(\mb)$, and $\pr{\zeta_V}{\zeta_V}_r\f bb^*$, because 
if $V\subseteq V_{\epsilon}$:
\[\norm{\pr{\zeta_V}{\zeta_V}_r-bb^*}
=\norm{\int_Gf_V(t)\big(\zeta (t)^*\zeta (t)-bb^*\big)dt}
\leq\int_Gf_V(t)\norm{\zeta (t)^*\zeta (t)-bb^*}dt
\leq\epsilon.\]
\par Now, for each $V\in\mv$ consider the kernel $k_V\in I_c(\mb)$ 
given by $k_V=\pr{\eta}{\zeta_V}_l$. If $V\subseteq V_{\epsilon}$: 
\begin{align*}
\norm{\big(\pr{\xi}{\xi}_l-k_V*k_V^*\big)\r{(r,s)}}
&=\norm{\xi(r)\xi(s)^*-\int_Gk_V(r,t)k_V^*(t,s)dt}\\
&=\norm{\eta(r)b\big(\eta(s)b\big)^*-\int_Gk_V(r,t)k_V(s,t)^*dt}\\
&=\norm{\eta(r)bb^*\eta(s)^*
        -\int_G\eta(r)\zeta_V(t)^*\big(\eta(s)\zeta_V(t)^*\big)^*dt}\\
&=\norm{\eta(r)bb^*\eta(s)^*
        -\int_G\eta(r)\zeta_V(t)^*\zeta_V(t)\eta(s)^*dt}\\
&=\norm{\eta(r)\big(bb^*-\pr{\zeta_V}{\zeta_V}_r\big)\eta(s)^*}\\
&\leq\epsilon\norm{\eta}_{\infty}^2
\end{align*}
Thus $k_V*k_V^*\f\pr{\xi}{\xi}_l$ in the \ilt, and hence
$\pr{\xi}{\xi}_l\in\mathfrak{c}$.
\end{proof}

\par Given an ideal $I$ of a $C^*$-algebra $A$, we set
$I^\perp:=\{a\in A:ax=0\,\forall x\in I\}$. Note that $I^\perp=0$ if
and only if $I$ is an 
essential ideal of $A$. On the other hand, if $\pi:I\to\mathcal{L}(F)$
is a non-degenerate representation of $I$ on the Hilbert
module $F$, then $\pi$ has a unique extension to a representation
$\tilde{\pi}:A\to\mathcal{L}(F)$, such that
$\tilde{\pi}(a)\pi(x)\xi=\pi(ax)\xi$, $\forall a\in A,x\in I$ and
$\xi\in F$. Moreover $\ker\tilde{\pi}=I^\perp$ if $\pi$ is faithful.      
\begin{thm}\label{thm:picpi2}
There exists a unique representation
$\Omega:\Bbbk(\mb)\to\ml(L^2(\mb))$ such that $\Omega(k)\xi=k\xi$
(recall Proposition~\ref{prop:bim}) $\forall k\in \Bbbk(\mb)$ and 
$\xi\in C_c(\mb)$, that is: 
\[\Omega(k)\xi|_r=\int_Gk(r,s)\xi(s)ds, \forall k\in \Bbbk_c(\mb),
\xi\in C_c(\mb) \textrm{ and }r\in G.\] Moreover we have:
\be
\item $\ker\Omega=I(\mb)^\perp$, so $\Omega$ is faithful if and only
  if $I(\mb)$ is an essential ideal.  
\item The restriction $\Omega|_{I(\mb)}:I(\mb)\to\ml(L^2(\mb))$ is a
  faithful representation.   
\item If $\mb$ is saturated,  
  then $\Omega$ is faithful. 
\ee
\end{thm}
\begin{proof}
We observe first that if $k\in\Bbbk_c(\mb)$ has null norm in
$\Bbbk(\mb)$, then $k(r,s)=0$, $\forall r,s\in G$. In fact, suppose
$\pi$ is a representation of $\mb$ whose restriction to $B_e$ is
faithful. Then, according to Proposition~\ref{prop:pic}, the corresponding
representation $\pi_c$ of $\Bbbk_c(\mb)$ is faithful and can be
extended to a representation of $HS(\mb)$. Thus the universal
$C^*$-seminorm defined by $HS(\mb)$ on $\Bbbk_c(\mb)$ is in fact a
$C^*$-norm, so $\norm{k}_{\Bbbk(B)}=0$ implies $k=0$.  \par Now, by
Proposition~\ref{prop:bim} and Lemma~\ref{lem:niguais} the map
$\langle\,,\,\rangle_l:C_c(\mathcal{B})\times C_c(\mathcal{B})\to
\Bbbk_c(\mathcal{B})$ is a pre-inner product on $C_c(\mathcal{B})$, so
we can complete $C_c(\mb)$ with respect to
$|\xi|:=\sqrt{\|\langle\xi,\xi\rangle_l\|}$, thus obtaining a left
Hilbert $\Bbbk(\mathcal{B})$-module $F$. Let us see that the right
action of $B_e$ on $C_c(\mb)$ is continuous with respect to $|\
|$. Suppose first that $b\in B_e$ is such that
$b=\langle\eta,\zeta\rangle_r$ for some $\eta,\zeta\in
C_c(\mathcal{B})$. Then if $\xi\in C_c(\mathcal{B})$ we have  
\[|\xi b|=|\xi\langle\eta,\zeta\rangle_r|
         =|\langle\xi,\eta\rangle_l\zeta|
         \leq|\xi|\, |\eta|\, |\zeta|.  
\]  
Moreover it is immediate that if $\xi,\eta\in C_c(\mathcal{B})$ and
$b\in B_e$, we have $\langle\xi b,\eta\rangle_l=\langle\xi,
\eta b^*\rangle_l$, so every $b\in
J:=\textrm{span}\{\pr{\eta}{\zeta}_r:\eta,\zeta\in C_c(\mb)\}$ acts as
an adjointable operator on~$F$. This gives a representation $J\to
\ml(F)^{\textsf{op}}$, which extends (\cite[VI-19.11]{fd}) to a representation of
$B_e$ because $J$ is an ideal of $B_e$. In particular each map
$C_c(\mb)\ni\xi\mapsto \xi b$ defines an adjointable operator on $F$,
for every $b\in B_e$. Note that this representation is isometric, for
$\xi b=0$ for every $\xi\in C_c(\mb)$ entails $\norm{\pr{\xi b}{\xi
    b}_l}=0$, that is $\xi(r)bb^*\xi(s)^*=0$ $\forall r,s\in G$,
$\xi\in C_c(\mb)$, and
since this implies  $b'b=0$ $\forall b'\in\mb$, we conclude that
$b=0$. Then the norm induced on $B_e$ by this representation agrees
with the original norm of $B_e$. Then for $\xi\in C_c(\mb)$
we have:  
\[ \norm{\pr{\xi}{\xi}_l}_{\Bbbk(\mb)}=|\xi|^2=\norm{\pr{\xi}{\xi}_r}_{B_e}=\norm{\xi}_{L^2(\mb)}^2.\]
Therefore the norms $|\ |$ and $\norm{\ }_{L^2(\mb)}$ agree on
$C_c(\mb)$, which implies that $F={L^2(\mb)}$, and that ${L^2(\mb)}$
is a $I(\mb)-B_e$-equivalence bimodule. Then the left action of $I(\mb)$
on ${L^2(\mb)}$ gives a faithful and non-degenerated representation 
$I(\mb)\to\ml({L^2(\mb)})$, whose unique extension
$\Omega:\Bbbk(\mb)\to\ml({L^2(\mb)})$ to a representation of
$\Bbbk(\mb)$ has kernel~$I(\mb)^\perp$.  
\par Next we show that $\Omega$ extends the action of $\Bbbk_c(\mb)$
on $C_c(\mb)$ provided by Proposition~\ref{prop:bim}. Fix $k\in
\Bbbk_c(\mb)$. Let $\xi,\eta,\zeta\in C_c(\mb)$. Then, using Proposition~\ref{prop:bim}:    
\begin{gather*}\Omega(k)(\xi\pr{\eta}{\zeta}_r)
=\Omega(k*\pr{\xi}{\eta}_l)(\zeta)
=\Omega(\pr{k  \xi}{\eta}_l)(\zeta)
=\pr{k  \xi}{\eta}_l\zeta
=k \xi\pr{\eta}{\zeta}_r
\end{gather*}
Thus $\Omega(k)(\xi b)=k  \xi b$ for any $\xi\in C_c(\mb)$ and $b\in
J$. Let $(e_i)$ be an approximate unit of $B_e$ contained in its dense
ideal $J$. Then given $\xi\in C_c(\mb)$ we have $\xi e_i\to \xi$, 
$k \xi e_i\to k \xi$ and $\Omega(k)(\xi e_i)\to\Omega(k)(\xi)$ in
${L^2(\mb)}$. 
Thus $\Omega(k)\xi=
k \xi.$
Finally, (3) follows from Proposition~\ref{prop:sat} and (1) above, since  the
saturation of $\mb$ implies $I(\mb)^\perp=0$. 
 \end{proof}
\begin{df}\label{df:nucr}
The $C^*$-subalgebra $\Bbbk_r(\mb):=\Omega(\Bbbk(\mb))$ of
$\ml(L^2(\mb))$ will be called the reduced $C^*$-algebra of kernels of
$\mb$.  
\end{df}
\par Thus by Theorem~\ref{thm:picpi2} we can identify $\Bbbk(\mb)$ with $\Bbbk_r(\mb)$ via $\Omega$ whenever $\mb$ is saturated.
\medskip
\par Recall that the right \rr of $G$ on $L^2(G,H)$ is 
$\rho :G\times L^2(G,H)\f L^2(G,H)$, such that 
$\rho_t(x)\r{r}=\Del(t)^{1/2}x(rt)$, $\forall r,t\in G$, $x\in L^2(G,H)$, 
where $\Delta$ is the modular function on $G$.

\begin{prop}\label{prop:pic}
Let $\pi:\mb\f B(H)$ be a \rep of the Fell bundle $\mb$ on a Hilbert space 
$H$. For $k\in\nucc{\mb}$, let $\pi_c(k):C_c(G,H)\f C_c(G,H)$ be such that  
$\pi_c(k)x\r{r}=\int_G\pi\big(k(r,s)\big)x(s)ds$
$\forall x\in C_c(G,H)$. Then we have that  
$\norm{\pi_c(k)x}_2\leq\norm{k}_2\,\norm{x}_2$, and hence $\pi_c(k)$ extends  
to an operator $\pi_c(k):L^2(G,H)\f L^2(G,H)$, with $\norm{\pi_c(k)}\leq
\norm{k}_2$. Moreover, $\pi_c:\nucc{\mb}\f B\big(L^2(G,H)\big)$ is a bounded  
\rep of the normed *--algebra $\nucc{\mb}$ such that 
$\pi_c\big(\beta_t(k)\big)=\rho_t\pi_c(k)\rho_{t^{-1}}$, where $\rho$ is the 
right \rr of $G$.  
If $\pi\r{B_e}$ is faithful, then so is $\pi_c$; if $\pi$ is non 
degenerate, then $\pi_c$ is non--degenerate.  
\end{prop}
\begin{proof}
Since $x\in C_c(G,H)$ and the map $\mb\times H\f H$ such that 
$(b,h)\ff\pi(b)h$ is continuous (\cite{fd}, VIII-8.8), we have that the map  
$G\times G\f H$ such that $(r,s)\ff \pi\big(k(r,s)\big)x(s)$ is continuous, 
and since it has compact support, we see that 
$r\ff\int_G\pi\big(k(r,s)\big)x(s)ds$ is a continuous map of compact support  
from $G$ to $H$, so in particular it belongs to $L^2(G,H)$. On the other hand:
\[\norm{\pi_c\big(k\big)x}_2^2
=\int_G\norm{\left[\int_G\pi\big(k(r,s)\big)x(s)ds\right]}^2dr
\leq\int_G\left[\left(\int_G\norm{\pi\big(k(r,s)\big)}^2ds\right)^{1/2}
 \,\norm{x}_2\right]^2dr
\leq\norm{k}_2^2\,\norm{x}_2^2,\]
so $\norm{\pi_c(k)}\leq\norm{k}_2$. It follows that $\pi_c(k)$ extends  
to a bounded operator $\pi_c(k):L^2(G,H)\f L^2(G,H)$, and  
$\pi_c$ extends by continuity to a bounded linear map 
$HS(\mb)\f B\big(L^2(G,H)\big)$. In addition, 
$\pi_c(k_1*k_2)=\pi_c(k_1)\pi_c(k_2)$: if $x\in C_c(G,H)$, 
$k_1,k_2\in\nucc{\mb}$ we have that
\[\pi_c(k_1*k_2)x\r{r}=
\int_G\pi\left[\int_Gk_1(r,t)k_2(t,s)dt\right]x(s)ds=
\int_G\pi\big(k_1(r,t)\big)\left[\int_G\pi\big(k_2(t,s)\big)x(s)ds
                                                 \right]dt \]
\[\pi_c(k_1)\pi_c(k_2)x\r{r}
=\int_G\pi\big(k_1(r,s)\big)\big[\pi_c(k_2)x\big](t)dt
=\int_G\pi\big(k_1(r,t)\big)\left[\int_G\pi\big(k_2(t,s)\big)x(s)ds
                                                 \right]dt\]
As for the involution, we have that 
$\pi_c(k^*)=\pi_c(k)^*$. In fact, if $k\in\nucc{\mb}$, $x$, $y\in C_c(G,H)$:
\[\pr{\pi_c(k)x}{y}
=\int_G\pr{\pi_c(k)x(r)}{y(r)}dr
=\int_{G^2}\pr{\pi\big(k(r,s)\big)x(s)}{y(r)}dsdr
=\int_{G^2}\pr{x(s)}{\pi\big(k(r,s)^*\big)y(r)}dsdr\]
\[\pr{x}{\pi_c(k^*)y}\!
=\!\int_G\!\pr{x(s)}{\pi_c(k^*)y(s)}ds\!
=\!\int_G\!\!\pr{x(s)}{\!\!\int_G\!\!\pi\big(k^*(s,r)^*\big)y(r)dr}ds\!
=\!\int_{G^2}\!\!\!\pr{x(s)}{\pi\big(k(r,s)^*\big)y(r)}dsdr.\]
It follows that $\pi_c:\nucc{\mb}\f 
B\big(L^2(G,H)\big)$ is a bounded \rep.
\par Now, if $k\in\nucc{\mb}$, $x\in C_c(G,H)$, $t,r\in G$:
\[\pi_c\big(\beta_t(k)\big)x\r{r}
=\int_G\pi\big(\beta_t(k)(r,s)\big)x(s)ds
=\int_G\Del(t)\pi\big(k(rt,st)\big)x(s)ds
=\int_G\pi\big(k(rt,u)\big)x(ut^{-1})du\]
\[\big(\rho_t\pi_c(k)\rho_{t^{-1}}\big)x\r{r}
=\int_G\Del(t)^{1/2}\pi\big(k(rt,u)\big)\rho_{t^{-1}}(x)(s)ds
=\int_G\pi\big(k(rt,u)\big)x(ut^{-1})du,\]
and hence $\pi_c\big(\beta_t(k)\big)=\rho_t\pi_c(k)\rho_{t^{-1}}$.
\par Suppose that $\pi\r{B_e}$ is faithful. If $\pi_c(k)=0$, for  
$k\in\nucc{\mb}$, we would have that $\int_G\pi\big(k(r,s)\big)x(s)ds=0$ 
$\forall x\in C_c(G,H)$.
In particular, if $(f_V)_{V\in\mv}$ is an approximate unit  
of $L^1(G)$ like in \ref{lem:lema}, and if $h\in H$, 
choosing $x(s)=f_V(t^{-1}s)h$ and taking limit with respect to $V$, by 
\ref{lem:lema} we have that $\pi\big(k(r,t)\big)h=0$, and therefore that  
$\pi\big(k(r,s)\big)=0$, $\forall r,s\in G$. Hence $\pi_c$ is also faithful. 
\par Finally, let us suppose that $\pi$ is non--degenerate, and let  
$y\in L^2(G,H)$ such that $\pr{\pi_c(k)x}{y}=0$, 
$\forall x\in L^2(G,H)$, $k\in\nucc{\mb}$. Then  
$\int_G\pr{\pi_c(k)x(r)}{y(r)}dr=0$, $\forall k\in\nucc{\mb}$, 
$\forall x\in C_c(G,H)$. Thus  
$\int_G\int_G\pr{\pi\big(k(r,s)\big)x(s)}{y(r)}dsdr=0$, $\forall k\in 
\nucc{\mb}$, $x\in C_c(G,H)$, so  
$\pr{\pi\big(k(r,s)\big)x(s)}{y(r)}=0$ almost everywhere 
in $G\times G$, and therefore $y(r)=0$ almost everywhere 
in $G$, because $\{\pi\big(k(r,s)\big)x(s):\ k\in\nucc{\mb}, 
(r,s)\in G\times G, x\in C_c(G,H)\} =H$. In fact, since $\pi$ is  
non--degenerate, $\pi\r{B_e}$ is also non--degenerate, and hence, 
given $h\in H$, there exist $b\in B_e$ and $h'\in H$, such that 
$\pi(b)h'=h$. Now there exist $b_1$, $b_2\in B_e$ 
such that $b=b_1b_2^*$, and sections $\eta_1$, $\eta_2\in C_c(\mb)$ such that 
$\eta_i(e)=b_i$, $i=1,2.$ Define $k:G\times G\f \mb$ such that 
$k(r,s)=\eta_1(r)\eta_2(s)^*$. It is clear that $k\in\nucc{\mb}$. 
Finally, taking $x\in C_c(G,H)$ such that $x(e)=h'$ we have:
$\pi\big(k(e,e)\big)x(e)
=\pi\big(\eta_1(e)\eta_2(e)^*\big)h'
=\pi(b_1b_2^*)h'
=\pi(b)h'=h .$
Consequently $\eta =0$ in $L^2(G,H)$, and therefore $\pi_c$ is 
non--degenerate. 
\end{proof}
\begin{lem}\label{lem:pi2}
Let $\pi:\mb\f B(H)$ be a \rep of a Fell bundle $\mb$. 
Then there exists a unique \rep $\pi_2:L^2(\mb)\f B\big(H,L^2(G,H)\big)$, 
such that $\forall \xi\in C_c(\mb)$ and $h\in H$ we have: 
$\pi_2(\xi)h\r{r}=\pi\big(\xi(r)\big)h$. 
If $\pi\r{B_e}$ is faithful, then $\pi_2$ also is faithful. 
\end{lem}
\begin{proof}
Let $\xi\in C_c(\mb)$ and $h\in H$; then:
\[ 
\norm{(\pi_2\xi)h}_2^2
=\int_G\pr{\pi\big(\xi(r)\big)h}{\pi\big(\xi(r)\big)h}dr
=\pr{\pi\left[\int_G\xi(r^*)\xi(r)dr\right]h}{h}
=\pr{\pi(\pr{\xi}{\xi}_r)h}{h}.\]
It follows that $\pi_2\xi\in L^2(G,H)$. On the other hand, since 
$\pi(\pr{\xi}{\xi}_r)$ is a positive operator, the equality above 
implies that $\norm{\pi_2\xi}=\norm{\pi(\pr{\xi}{\xi}_r)}^{1/2}$. Thus, 
$\norm{\pi_2\xi}\leq\norm{\xi}$, and $\norm{\pi_2\xi}=\norm{\xi}$ if  
$\pi\r{B_e}$ is faithful. It follows that we may extend 
$\pi_2$ to $L^2(\mb)$, and this extension is an isometry if 
$\pi\r{B_e}$ is faithful. 
\par An easy computation shows that $(\pi_2\xi)^*$ is determined by the 
formula $(\pi_2\xi)^*x=\int_G\pi\big(\xi(r)^*\big)x(r)dr\in H$. From this we 
see that $\pi_2(\xi)\pi_2(\eta)^*\pi_2(\zeta)
=\pi_2\big(\xi\pr{\eta}{\zeta}_r\big)$, $\forall \xi,\eta,\zeta\in C_c(\mb)$.
Indeed, if $h\in H$, $r\in G$:
\[ \pi_2(\xi)\pi_2(\eta)^*\pi_2(\zeta)h\r{r}
=\pi\big(\xi(r)\big)\int_G\pi\big(\eta(s)^*\big)\pi_2(\zeta)h\r{s}ds
=\int_G\pi\big(\xi(r)\eta(s)^*\big)\pi\big(\zeta(s)\big)hds\]
\[\pi_2\big(\xi\pr{\eta}{\zeta}_r\big)h\r{r}
=\pi\big(\xi(r)\pr{\eta}{\zeta}_r\big)h
=\int_G\pi\big(\xi(r)\eta(s)^*\zeta(s)\big)hds
=\int_G\pi\big(\xi(r)\eta(s)^*\big)\pi\big(\zeta(s)\big)hds\]
Therefore, $\pi_2$ is a \hm of positive $C^*$-trings. 
\end{proof}
\begin{prop}\label{prop:picpi2}
Let $\pi:\mb\to B(H)$ be a representation of the Fell bundle $\mb$,
and let $\pi_c:\Bbbk_c(\mb)\to B(L^2(G,H))$ and $\pi_2:L^2(\mb)\to
B(H,L^2(G,H))$ be the representations provided by
Propostion~\ref{prop:pic} and Lemma~\ref{lem:pi2}, respectively. Then we have, for all
$k\in\Bbbk_c(\mb)$, $b\in B_e$, and $\xi,\eta\in E$: 
\be
   \item[(i)] $\pi_2(k\xi)=\pi_c(k)\pi_2(\xi)$
   \item[(ii)] $\pi_c(\pr{\xi}{\eta}_l)=\pi_2(\xi)\pi_2(\eta)^*$
   \item[(iii)] $\pi_2(\xi b)=\pi_2(\xi)\pi(b)$
   \item[(iv)] $\pi(\pr{\xi}{\eta}_r)=\pi_2(\xi)^*\pi_2(\eta)$
\ee
\end{prop}
\begin{proof}
\par Let $h\in H$, $r\in G$, $x\in L^2(G,H)$:\\
(a) $\pi_2(k\xi)=\pi_c(k)\pi_2(\xi)$: \[\big(\pi_2(k\xi)\big)h\r{t}
  =\int_G\pi\big(k(t,s)\big)\pi\big(\xi(s)\big)hds
  =\int_G\pi\big(k(t,s)\big)(\pi_2\xi)h\r{s}ds
  =\pi_c(k)(\pi_2\xi)h\r{t}.\]
(b)$\pi_c(\pr{\xi}{\eta}_l)=\pi_2(\xi)\pi_2(\eta)^*$:
  \[\pi_c(\pr{\xi}{\eta}_l)x\r{t}
  =\int_G\pi\big(\xi(t)\big)\pi\big(\eta(s)^*\big)x(s)ds
  =\pi\big(\xi(t)\big)\int_G\pi\big(\eta(s)^*\big)x(s)ds
  =(\pi_2\xi)(\pi_2\eta)^*x\r{t}.\]
(c)$\pi_2(\xi b)=\pi_2(\xi)\pi(b)$: 
 $\pi_2(\xi b)h\r{r}
  =\pi\big(\xi b(r)\big)h
  =\pi\big(\xi(r)b\big)h
  =\pi\big(\xi(r)\big)\pi(b)h\r{r}
  =\pi_2(\xi)\pi(b)h\r{r}.$\\
(d)$\pi(\pr{\xi}{\eta}_r)=\pi_2(\xi)^*\pi_2(\eta)$:
  \[\pi_2(\xi)^*\pi_2(\eta)h
  =\int_G\pi\big(\xi(r)^*\big)\pi_2(\eta)h\r{r}ds
  =\int_G\pi\big(\xi(r)^*\big)\pi\big(\eta(r)\big)hdr
  =\int_G\pi\big(\xi(r)^*\eta(r)\big)hdr
  =\pi(\pr{\xi}{\eta}_r)h\] 
\end{proof}

\begin{prop}
Let $\mb$ be a Fell bundle and $I=I(\mb)$. Suppose that $\pi:\mb\to
B(H)$ is a representation of the Fell bundle $\mb$. Then:  
 \be
   \item $\pi_c$ extends to a representation $\pi_H:\Bbbk(\mb)\to
     B(L^2(G,H))$, such
     that \[\pi_H(\beta_t(k))=Ad_{\rho_t}(\pi_H(k)),\ \forall k\in
     \Bbbk(\mb),\] where $\rho$ is the right regular representation of
     $G$ on $L^2(G,H)$.    
   \item Let $\Bbbk_\pi(\mb)$ be the closure of $\pi_H(\Bbbk(\mb))$
     and $I_\pi:=\pi_{H}(I)$, where $\pi_H$ is as in part (1). Then
     $(\Bbbk_\pi(\mb),Ad_{\rho})$ is the enveloping action of
     $(I_\pi,Ad_{\rho}|_{I_\pi})$.    
   \item If $\pi|_{B_e}$ is injective, then
     $\pi_H:(\Bbbk(\mb),\beta)\to (\Bbbk_\pi(\mb),Ad_\rho)$ is an
     isomorphism.  
 \ee
\end{prop}
\begin{proof}
The first statement follows at once from Proposition~\ref{prop:pic}, and the
second statement is clearly a consequence of (1) together with
Proposition~\ref{prop:kenv}. To see that (3) also holds
note that, since $\pi_c$ is injective, the map
$\pr{\,}{\,}_\pi:=\pi_c\circ\pr{\,}{\,}_l:C_c(\mb)\times C_c(\mb)\to
\pi_c(I_c(\mb))$ is a pre-inner product on $C_c(\mb)$. Let $F_\pi$ be
the corresponding completion of $C_c(\mb)$. Then, proceeding as in the proof
of Theorem~~\ref{thm:picpi2}, we see that $F_\pi=L^2(\mb)$, and this proves that
$\pi_H|_I:I\to I_\pi$ is an isomorphism. Then, as we can identify $I$
and $I_\pi$, it follows that $\pi_H:\Bbbk(\mb)\to\Bbbk_\pi(\mb)$ is an
isomorphism by Proposition~\ref{prop:kenv} and the uniqueness of the enveloping action.    
\end{proof}
\begin{rk}\label{rk:nucreds}
Let $k\in\nucc{\mb}$ and $f\in C_c(\mb)$, and consider $h:G\times G\f\mb$  
such that $h(r,t)=\int_Gf(s)k(s^{-1}r,t)ds$. Note that 
$h(r,t)\in B_{rt^{-1}}$, $\forall r,t\in G$, because for all $r,s,t\in G$, 
we have that $f(s)k(s^{-1}r,t)\in B_sB_{s^{-1}rt^{-1}}\subseteq B_{rt^{-1}}$. 
Hence $h\in\nucc{\mb}$. Now, considering 
$\xi\in C_c(\mb)\subseteq L^2(\mb)$, we have:
$
\La_f\Om_k\xi\r{r}
=\int_Gf(s)\Om_k(\xi)(s^{-1}r)ds
=\int_G\int_Gf(s)k(s^{-1}r,t)\xi(t)dtds 
=\int_Gh(r,t)\xi(t)dt
=\Om_h\xi\r{t}.$
It follows that $C^*_r(\mb)\subseteq M\big(\nucr{\mb}\big)$. 
In particular, any \rep of $\nucr{\mb}$ defines a \rep of $C^*_r(\mb)$. 
\end{rk}
\par In \cite{dualpaut}, Quigg introduced the notion of   
\textit{ideal property} of a \cs: a property $\mathcal{P}$ of a \cs is 
ideal if every \cs has a largest ideal with property $\mathcal{P}$, and if 
this property is invariant by Morita equivalence and inherited by ideals. 
\medskip
\par From now on we will maintain the notation $I=I(\mb)$ unless there exists some kind of ambiguity.  
\begin{cor}\label{cor:pideal}
Let $\mathcal{P}$ be an ideal property. Then $B_e$ has property  
$\mathcal{P}$ iff $\nuc{\mb}$ has property $\mathcal{P}$. Therefore 
$B_e$ is liminal, antiliminal, postliminal or nuclear, iff $\nuc{\mb}$ is 
respectively liminal, antiliminal, postliminal or nuclear. 
\end{cor}
\begin{proof}
Since $B_e$ and $I$ are Morita equivalent, $B_e$ has property $\mathcal{P}$ if 
and only if $I$ has property $\mathcal{P}$. On the other hand, by 
Proposition~\ref{prop:class}, $I$ has property $\mathcal{P}$ if and only if 
$\nuc{\mb}$ has property $\mathcal{P}$. 
\end{proof}

\par We would like to emphasize the functorial 
nature of $\nuc{\mb}$. Suppose that $\ma=(A_t)$ and $\mb=(B_t)$ 
are Fell bundles, and that $\phi:\ma\f\mb$ is a Fell bundle \hm. For instance, 
$\phi$ could be a \hm induced by a morphism between \pas. 
If $k\in\nucc{\ma}$, it is clear that  
$\nucc{\phi}(k)\in\nucc{\mb}$, where 
$\nucc{\phi}(k)\r{(r,s)}=\phi(k(r,s))$. Moreover, 
it is easy to check that the map  
$\nucc{\phi}:\nucc{\ma}\f\nucc{\mb}$ is a \hm of *-algebras, which is 
injective if $\phi$ is injective. 
In addition   
$\norm{\nucc{\phi}(k)}_2\leq\norm{k}_2$, so $\nucc{\phi}$ extends to a \hm  
$HS(\phi):HS(\ma)\f HS(\mb)$ of Banach *--algebras. Hence, there is an induced 
\hm of \css $\nuc{\phi}:\nuc{\ma}\f\nuc{\mb}$. It is straightforward to 
verify that $\ma\ff\nuc{\ma}$ and $(\ma\stackrel{\phi}{\to}\mb)\ff
\big(\nuc{\ma}\stackrel{\nuc{\phi}}{\to}\nuc{\mb}\big)$ is a functor from the 
category of Fell bundles to the category of \css.
\par If $\phi:\ma\f\mb$ is such that $\phi|_{A_e}$ is injective, then  
$\nuc{\phi}:\nuc{\ma}\f\nuc{\mb}$ also is injective. To see this, let  
$\pi:\mb\f B(H)$ be a \rep such that $\pi\r{B_e}$ 
is faithful. Then $\pi\phi:\ma\f B(H)$ is a \rep such that $\pi\phi\r{A_e}$ 
is faithful. Then $(\pi\phi)_H:\nuc{\ma}\to\Bbbk_{\pi\phi}(\ma)$ and
$\pi_H:\nuc{\mb}\to\Bbbk_\pi(\mb)$ are isomorphisms. Now a
straightforward computation shows that $\pi_H\Bbbk(\phi)=(\pi\phi)_H$,
and since the latter is injective, so must be $\Bbbk(\phi)$.  
\par Next we will refine Theorem~\ref{thm:picpi2} in the case the group $G$ is discrete. First we introduce a definition: 
\begin{df}\label{df:esssat}
We say that a Fell bundle $\mb$ is \textit{essentially saturated}
whenever every $B_tB_t^*$ is an essential ideal of $B_e$. 
\end{df}
It is no hard to see that $\mb$ is essentially saturated if and only
if forall $s,t\in G$ we have $(B_sB_t)^\perp=0$ in the right
$B_e$-module $B_{st}$.    
As we shall see next, if $\mb$ is an essentially saturated Fell
bundle over a discrete group $G$, the elements of $\Bbbk(\mb)$ can be
represented by functions defined on $G\times G$. This will allow us to  
see that the faithfulness of $\Omega$ is equivalent
to the essential saturation of~$\mb$. 

\begin{lem}\label{lem:esssat}
Let $\mb$ be a Fell bundle over the discrete group $G$, such that
$(B_tB_t^*)^\perp=0$, $\forall t\in G$. Then for every
$k\in\Bbbk(\mb)$ there exists a unique function $G\times G\to \mb$,
which we also denote $k$, such that $\Omega(k)\xi(r)=\sum_{s\in
  G}k(r,s)\xi(s)$, $\forall \xi\in L^2(\mb)$, $r\in G$. 
\end{lem}  
\begin{proof}
First note that for every pair of elements $u, v\in G$  we have a map $B_u\to\ml(B_v,B_{uv})$ given by $b_u\mapsto M_{b_u}$ such that
$M_{b_u}(b_v)=b_ub_v$. It is easy to check that it is in fact
a homomorphism of $C^*$-trings, so $\norm{M_{b_u}}\leq\norm{b_u}$. The assumption 
$(B_{v}B_{v}^*)^\perp=0$ ensures that it is injective, for 
$M_{b_u}=0$ implies $b_u^*b_u\in (B_vB_v^*)^\perp$, that is,
$b_u=0$. Then $b_u\mapsto M_{b_u}$ is an isometry (by Proposition~\ref{prop:z(pi)}).          
\par Now let $k\in\Bbbk(\mb)$, $r,s\in G$. Consider the map $k_s^r:B_s\to B_r$
such that $k_s^r(b_s):=\Omega(k)(b_s\delta_s)(r)$, where
$b_s\delta_s(t)=b_s$ if $t=s$, and $b_s\delta_s(t)=0$ otherwise. Then
\[\norm{k_s^r(b_s)}=\norm{\Omega(k)(b_s\delta_s)(r)}\leq\norm{\Omega(k)(b_s\delta_s)}\leq\norm{\Omega(k)}\,\norm{b_s}\leq\norm{k}\,\norm{b_s}.\] 
Then $k_s^r$ is bounded, and is easily checked that it is adjointable,
with $(k_s^r)^*=(k^*)_r^s$. In case
$k\in\Bbbk_c(\mb)$, $k_s^r$ is just multiplication by $k(r,s)$ on
the left, so $\norm{k_s^r}=\norm{k(r,s)}$ as shown at the beginning of
the proof. Suppose now that
$k\in\Bbbk(\mb)$ and $(k_n)\subseteq\Bbbk_c(\mb)$ is a sequence
converging to $k$ in $\Bbbk(\mb)$. Then $(k_n)_s^r\to k_s^r$ in
$\mathcal{L}(B_s,B_r)$. Since this implies that $(k_n(r,s))$ is a
Cauchy sequence, there exists an element in $B_{rs^{-1}}$, which we
call $k(r,s)$, such that $k_n(r,s)\to k(r,s)$. But this implies that
the operator $k_s^r$ is multiplication by $k(r,s)$ on the
left. Therefore, if $k\in \Bbbk(\mb)$ and $\xi=\sum_s\xi(s)\delta_s\in
L^2(\mb)$: $
\Omega(k)\xi(r)
=\sum_s\Omega(k)(\xi(s)\delta_s)(r)
=\sum_sk_s^r(\xi(s))
=\sum_s k(r,s)\xi(s). 
$
\end{proof}

\begin{prop}\label{prop:esssat}
Let $\mb$ be a Fell bundle over the discrete group $G$. Then $\Omega$
is faithful if and only if $\mb$ is essentially saturated.  
\end{prop}
\begin{proof}
Suppose $(B_tB_t^*)^\perp=0$, $\forall t\in G$, and let
$k\in\ker\Omega$, which we identify with a function $G\times G\to\mb$
as in Lemma~\ref{lem:esssat}. For every $B_t$ and $r\in G$ we have
$0=\Omega(k)(B_t\delta_t)=\sum_sk(r,s)B_t=k(r,t)B_t$, so 
$k(r,t)^*k(r,t)\in (B_tB_t^*)^\perp=0$. Then $k=0$, which implies that
$\Omega$ is injective. Suppose conversely that there exist $t\in G$ and $a\in (B_tB_t^*)^\perp$ such that $a\neq 0$. Consider $k_a:G\times G\to\mb$ such that $k_a(t,t)=a$, and $k_a(r,s)=0$ otherwise. Then
$k_a\neq 0$, but $\Omega(k_a)=0$, as it is easily checked.
\end{proof}

\par Proposition~\ref{prop:esssat} applies in particular to bundles
associated to certain partial actions. For instance, suppose $h$ is a
homeomorphism between two open and dense subsets of a locally compact
and Hausdorff space $X$. Then the Fell bundle of the partial action
of $\Z$ associated to $h$ is essentially saturated.  

\section{Existence and uniqueness of the Morita enveloping action} 
\label{sec:exun}
\par In the sequel we prove the main result of the paper: any \pa has a \mea, 
which is unique up to Morita equivalence. 
To do this, we apply the previous results to the case when the Fell 
bundle $\mb$ is the one associated with a \pa $\al$. We show that  
the natural action on $\nuc{\mb}$ is a \mea of $\al$. 
More precisely, assume that $\mb$ is the Fell bundle of a \pa 
$\al=(\{ D_t\}_{t\in G},\{\al_t\}_{t\in G})$. Then define  
on $E$ the \pa\ $\ga=(\{ E_t\}_{t\in G},\{\ga_t\}_{t\in G})$, where 
$E_t=ED_t$, and $\ga_t:E_{t^{-1}}\f E_t$ is such that $\ga_t(\xi)\r{r}=
\Del(t)^{1/2}\xi(rt)$, for all $\xi\in E_{t^{-1}}\cap C_c(\mb)$. In Theorem 
\ref{thm:super}, we prove that  
$\ga^r=\al$ and $\ga^l=\nal$, where $\nal$ is the natural action of $G$ on 
$\nuc{\mb}$. It follows that $(\nuc{\mb},\nal)$ is a \mea
of $\al$, and therefore any \pa has a \mea. We prove later in Proposition 
\ref{prop:super} that the \mea is unique up to Morita equivalence.

\begin{prop}\label{prop:restc}
Let $\mb$ be a \fb, $\beta$ the natural action of $G$ on  
$\nucc{\mb}$, and $I_c(\mb)$ the two--sided 
*-ideal of $\nucc{\mb}$ defined after Proposition \ref{prop:bim}, that is,
$I_c(\mb)=\gen\pr{C_c(\mb)}{C_c(\mb)}_l$. 
For $t\in G$, consider the two--sided *-ideal 
$I_t^c=I_c(\mb)\bigcap \beta_t(I_c(\mb))$ of $\nucc{\mb}$.
Then the closure of $I_t^c$ in the \ilt is the set  
\[\{ k\in \nucc{\mb}:\, k(r,s)\in 
\ov{(\gen B_rB_s^*)\cap(\gen B_{rt}B_{st}^*)},\, \forall r,s\in G\}.\]  
\end{prop}
\begin{proof}
Let $k\in I_c(\mb)$. There exist $\xi_i$, $\eta_i\in C_c(\mb)$, such that 
$k=\sum_i\pr{\xi_i}{\eta_i}_l$. In particular, 
$k(r,s)=\sum_i\xi_i(r)\eta_i(s)^*\in\gen B_rB_s^*$, $\forall r,s\in G$. 
Now, if $k\in I_t^c$, then $\beta_{t^{-1}}(k)\in I_c(\mb)$, and therefore   
$\beta_{t^{-1}}(k)(r,s)\!\in \gen B_rB_s^*$, $\forall r,s\in G$, that is, 
$\Del(t)^{-1} k(rt^{-1},st^{-1})\in\gen B_rB_s^*$, $\forall r,s\in G$. 
Equivalently, $k(r,s)\in\gen B_{rt}B_{st}^*$, $\forall r,s\in G$. 
Let $\Theta=\bats{C_c(G)}{C_c(G)}$, which is a  
dense subalgebra of $C_c(G\times G)$ in the \ilt. If $\phi$, $\psi\in C_c(G)$, 
$\xi$, $\eta\in C_c(\mb)$, then $(\phi\otimes\psi)\pr{\xi}{\eta}_l=
\pr{\phi\xi}{\bar{\psi}\eta}_l$, and hence  
$\Theta I_t^c\subseteq I_t^c$. The result follows now from Lemma 
\ref{lem:bolu}, 
\end{proof}

\par Suppose now that $\al=(\{ D_t\}_{t\in G},\{\al_t\}_{t\in G})$ 
is a \pa of $G$ on the \cs $A$, and let  
$\mb$ be its associated \fb. For each $t\in G$, consider 
the subset $E_t:=\ov{\gen}ED_t$ of $E$. By the Cohen--Hewitt theorem, 
$E_t=\{ \xi a:\, \xi\in E, a\in D_t\}$. By Proposition 
\ref{prop:lattice}, $E_t$ is an ideal of $E$. 
If $\xi:G\f\mb$ is a section, there exists a unique function  
$\xi':G\f A$ such that $\xi'(r)\in D_r$, $\forall r\in G$, and such that  
$\xi(r)=(r,\xi'(r))$. $\xi$ is continuous with compact support if and  
only if $\xi'$ is continuous and of compact support.
Suppose that $\xi\in C_c(\mb)$, and $a\in D_t$. Then  
$\xi a\r{r}=\xi(r)(e,a)=(r,\xi'(r))(e,a)
=(r,\al_r\big(\al_{r^{-1}}(\xi'(r))a\big))$. Since $\al_{r^{-1}}(\xi'(r))a\in
D_{r^{-1}}\cap D_t$, we have that $\xi a\r{r}\in D_r\cap D_{rt}$, $\forall 
r\in G$. Thus $E_t$ is the closure of  
$E_t^c:=\{\xi\in C_c(\mb):\, \xi'(r)\in D_r\cap D_{rt},\, \forall r\in G\}$. 

\begin{prop}\label{prop:restr}
Let $\al=(\{ D_t\}_{t\in G},\{\al_t\}_{t\in G})$ be a \pa of $G$ on $A$, 
with associated Fell bundle $\mb$. Then: 
\be 
 \item $B_rB_s^*=\{(rs^{-1},D_r\cap D_{rs^{-1}})\}$, $\forall r,s\in G$.
 \item $B_rB_s^*\cap B_{rt}B_{st}^* 
       =\{(rs^{-1},D_r\cap D_{rt}\cap D_{rs^{-1}})\}$, $\forall r,s,t\in G$.
 \item The closures of $\pr{E_t^c}{E_t^c}_l$ and $I_t^c$ (defined in 
       \ref{prop:restc}) in the \ilt agree, 
       and they are equal to the set:
       $\{ k\in \nucc{\mb}:\ k(r,s)\in B_rB_s^*\cap B_{rt}B_{st}^*,\, 
       \forall r,s\in G\}$
 \item $E_t^l=I_t$, where $I_t$ is the closure $I_t^c$ in $I$. 
 \item If $\nal$ is the natural action of $G$ on $\nuc{\mb}$, then  
       $\nal\r{I}=\big(\{ I_t\}_{t\in G},\{\nal_t\r{I_t}\}_{t\in G}\big)$.  
\ee
\end{prop}
\begin{proof}
To prove 1., recall that $B_r=(r,D_r)$, $\forall r\in G$. Then 
       \begin{align*}
       (r,D_r)(s,D_s)^*
       &=(r,D_r)(s^{-1},\al_{s^{-1}}(D_s))\\
       &=\big(rs^{-1},\al_r\big(\al_{r^{-1}}(D_r)\al_{s^{-1}}(D_s)\big)\big)\\
       &=\big(rs^{-1},\al_r(D_{r^{-1}}\cap D_{s^{-1}})\big)\\
       &=(rs^{-1}, D_r\cap D_{rs^{-1}})
       \end{align*}
2. $B_rB_s^*\cap B_{rt}B_{st}^* 
       =\{(rs^{-1},D_r\cap D_{rs^{-1}})\}\bigcap 
       \{(rs^{-1},D_{rt}\cap D_{rs^{-1}})\}
       =\{(rs^{-1},D_r\cap D_{rt}\cap D_{rs^{-1}})\}$.\\
3. Let $\xi$, $\eta\in E$, $\xi(r)=\big(r,\xi'(r)\big)$ and 
   $\eta(r)=\big(r,\eta'(r)\big)$, $\forall r\in G$. We have:
       \begin{align*}
       \pr{\xi}{\eta}_l(r,s)
       &=\xi(r)\eta(s)^*\\
       &=\big(r,\xi'(r)\big)\big(s,\eta'(s)\big)^*\\
       &=\big(r,\xi'(r)\big)\big(s^{-1},\al_{s^{-1}}(\eta'(s)^*)\big)\\
       &=\big(rs^{-1},\al_r\big(\al_{r^{-1}}(\xi'(r))
         \al_{s^{-1}}(\eta'(s))\big)\big) 
       \end{align*}
       Hence, if $\xi$, $\eta\in E_t^c$, then 
       \begin{align*}
       \pr{\xi}{\eta}_l\in 
       \big(rs^{-1},\al_r\big(\al_{r^{-1}}(D_r\cap D_{rt})
         \al_{s^{-1}}(D_s\cap D_{st})\big)\big)\\ 
       &=\big(rs^{-1},\al_r(D_{r^{-1}}\cap D_t\cap D_{s^{-1}})\big)\\
       &=\big(rs^{-1},D_r\cap D_{rt}\cap D_{rs^{-1}}\big)\\
       &=B_rB_s^*\cap B_{rt}B_{st}^*,\ \ \text{ by 2.}
       \end{align*}
       Since $B_rB_s^*\cap B_{rt}B_{st}^*$ is a closed linear space,  
       Lemma \ref{lem:bolu} shows that the closure of
       $\pr{E_t^c}{E_t^c}_l$ in the \ilt is $\{ k\in\nucc{\mb}:\ k(r,s)\in 
       B_rB_s^*\cap B_{rt}B_{st}^*,\, \forall r,s\in G\}$, and by  
       \ref{prop:restc}, this set agrees with the closure of $I_t^c$ in the 
       \ilt. So 3. is proved.\\
4. This is an immediate consequence of 3.\\
5. We must show that $I\cap\nal_t(I)=I_t$. Since  
       $I_t^c=I_c(\mb)\cap\nal_t(I_c(\mb))\subseteq I\cap\nal_t(I)$, we have
       that $I_t\subseteq I\cap\nal_t(I)$.
       \par To prove the converse inclusion, by \ref{prop:lattice} 
       it is enough to show that $\big(I\cap\nal_t(I)\big)E\subseteq E_t$. 
       Let $x\in I\cap\nal_t(I)$. Then $x=\nal_t(y)$, for some $y\in I$. 
       For a given $\epsilon>0$, let $k\in\nucc{\mb}$ be such that  
       $\norm{y-k}<\epsilon$. Then $\norm{x-\nal_t(y)}<\epsilon$, because  
       $\nal_t$ is an isometry. Consider now $\xi\in C_c(\mb)$. 
       We have:
       \[\nal_t(k)\xi\r{r}=\int_G\Del(t)k(rt,st)\xi(s)ds\in B_{rt}B_{st}B_s.\]
       By 1., $B_{rt}B_{st}B_s=(rs^{-1},D_{rt}\cap D_{rs^{-1}})(s,D_s)$, 
       and therefore 
       \begin{align*}
       \nal_t(k)\xi\r{r}\in
       &\big(r,\al_{rs^{-1}}(\al_{sr^{-1}}(D_{rt}\cap D_{rs^{-1}})D_s)\big)\\
       &=\big(r,\al_{rs^{-1}}(D_{sr^{-1}}\cap D_{st}\cap D_s)\big)\\
       &=(r,D_{rs^{-1}}\cap D_{rt}\cap D_r),
       \end{align*}
       from where it follows that $\nal_t(k)\xi\in E_t^c$. By the continuity of 
       the product, it follows now that $xE\subseteq E_t$, 
       $\forall x\in I\cap\nal_t(I)$, and this ends the proof. 
\end{proof}
\par Consider now, for each $t\in G$, the map 
$\ga_t':C_c(G,A)\f C_c(G,A)$ 
such that $\ga_t'(\xi')\r{r}=\Del(t)^{1/2}\xi'(rt)$. Suppose that  
$\xi\in E_{t^{-1}}^c$,  where $\xi(r)=\big(r,\xi'(r)\big)$, $\forall r\in G$. 
Then the map $r\ff (r,\ga_t'(\xi')(r))$ is an element of $E_t^c$:
$\ga_t'(\xi')(r)=\Del(t)^{1/2}\xi'(rt)\in D_{rt}\cap D_{rtt^{-1}}
=D_r\cap D_{rt}$. Therefore, we may define $\ga_t:E_{t^{-1}}^c\f E_t^c$, 
such that $\ga_t(\xi)(r)=(r,\ga_t'(\xi')(r))=(r,\Del(t)^{1/2}\xi'(rt))$.   

\begin{thm}\label{thm:super}
Any \pa on a \cs has a \mea. More precisely, let 
$\al=(\{ D_t\}_{t\in G},\{\al_t\}_{t\in G})$ be a \pa of $G$ on a \cs $A$, 
and let $\mb$ the \fb associated with $\al$. 
Let $E_t^c$ and $\ga_t$ be the ones defined previously. Then:
\be
 \item For each $t$, $\ga_t:E_{t^{-1}}^c\f E_t^c$ is an isometry  
       that extends by continuity to an isomorphism  
       $\ga_t:E_{t^{-1}}\f E_t$ of $C^*$-trings.
 \item $\ga=(\{ E_t\}_{t\in G},\{\ga_t\}_{t\in G})$ is a \pa of $G$ on $E$.
 \item $\ga^r=\al$, and $\ga^l=\nal\r{I}$, where $\nal$ is the natural 
       action of $G$ on $\nuc{\mb}$.
\ee
\end{thm}
\begin{proof}
Let $\xi$, $\eta\in C_c(\mb)$. Identifying $A$ with the unit fiber of  
$\mb$ over the unit element of $G$, and  
the kernel $k\in\nucc{\mb}$ with $k':G\times G\f A$ such that  
$k(r,s)=\big(rs^{-1},k'(r,s)\big)$, we have the identities: 
\begin{equation}\label{e:izquierda}
\pr{\xi}{\eta}_l(r,s)
=\al_r\big(\al_{r^{-1}}(\xi'(r))\al_{s^{-1}}(\eta'(s))\big)
\end{equation}
\begin{equation}\label{e:derecha}
\pr{\xi}{\eta}_r=\int_G\al_{r^{-1}}\big(\xi'(r)^*\eta'(r)\big)dr
\end{equation}
If $\xi$, $\eta\in E_{t^{-1}}^c$, by (\ref{e:derecha}) we have that: 
\begin{align*}
\pr{\ga_t(\xi)}{\ga_t(\eta)}_r
&=\int_G\al_{r^{-1}}\big(\ga_t'(\xi')(r)^*\ga_t'(\eta')(r)\big)dr\\
&=\int_G\al_{r^{-1}}
\big(\Del(t)^{1/2}\xi'(rt)^*\Del(t)^{1/2}\eta'(rt)\big)dr\\
&=\int_G\Del(t^{-1})\Del(t)\al_{ts^{-1}}\big(\xi'(s)^*\eta'(s)\big)ds\\
&=\int_G\al_t\al_{s^{-1}}\big(\xi'(s)^*\eta'(s)\big)ds\\
&=\al_t(\pr{\xi}{\eta}_r),
\end{align*}
Then $\ga_t$ is an isometry with dense range, and hence it extends to a 
bijective isometry $\ga_t:E_{t^{-1}}\f E_t$, such that 
$\pr{\ga_t(\xi)}{\ga_t(\eta)}_r=\al_t(\pr{\xi}{\eta}_r)$, $\forall \xi$, 
$\eta\in E$. On the other hand, if $\xi\in E_{t^{-1}}^c$ and $a\in D_{t^{-1}}$, 
we have that  
$\xi a(r)=(r,\xi'(r))(e,a)=\big(r,\al_r(\al_{r^{-1}}\xi'(r)a\big)$. 
Thus:
\begin{align*}
\ga_t(\xi a)(r)
&=\big(r,\Del(t)^{1/2}\al_{rt}\big(\al_{{rt}^{-1}}(\xi'(rt)\big)a\big)\\
&=\big(r,\Del(t)^{1/2}\al_{rt}\big(\al_{{t}^{-1}r^{-1}}(\xi'(rt)\big)
\al_{t^{-1}}(\al_t(a))\big)\\
&=\big(r,\Del(t)^{1/2}\al_r\big(\al_{{r}^{-1}}(\xi'(rt)\big)\al_t(a)\big)\\
&=\big(r,\Del(t)^{1/2}\xi'(rt)\big)\big(e,\al_t(a)\big)\\
&=\ga_t(\xi)\al_t(a)(r).
\end{align*}
It follows that $\ga_t$ is an isomorphism of \cts, and $\ga_t^l=\al_t$,  
which proves 1.
\par Let us see that $\ga$ is a set theoretic \pa: it is clear that 
$E_e=E$ and $\ga_e=id_E$. Suppose now that $\xi\in E_{t^{-1}}^c$ is 
such that $\ga_t(\xi)\in E_{s^{-1}}^c$. Then:
\[\ga_s\ga_t(\xi)\r{r}
=\big(r,\Del(s)^{1/2}\big(\ga_t'(\xi')\big)(rs)\big)
=\big(r,\Del(s)^{1/2}\Del(t)^{1/2}\xi'(rst)\big)
=\big(r,\Del(st)^{1/2}\xi'(rst)\big)
=\ga_{st}(\xi)\r{r}.\]
Then $\ga_{st}$ is an extension of $\ga_s\ga_t$, and therefore $\ga$ is  
a set theoretic \pa on $E$. So, 2. is proved, up to the continuity of $\ga$, 
that will be proved later. 
\par Since $E_t^l=D_t$ by \ref{prop:lattice}, 
the computation that showed that $\ga_t$ was an isometry also implies that 
$\ga^r=\al$. To see that $\ga^l=\nal\r{I}$, observe that if $\xi$, 
$\eta\in E_{t^{-1}}^c$, by (\ref{e:izquierda}) we have:
\begin{align*}
\pr{\ga_t(\xi)}{\ga_t(\eta)}_l(r,s)
&=\al_r\big(\al_{r^{-1}}(\ga_t'(\xi')(r))\al_{s^{-1}}(\ga_t'(\eta')(s))\big)\\
&=\al_r\big(\al_{r^{-1}}(\Del(t)^{1/2}\xi'(rt))
            \al_{s^{-1}}(\Del(t)^{1/2}\eta'(st))\big)\\
&=\Del(t)\al_{rt}\big(\al_{{rt}^{-1}}(\xi'(rt))
                       \al_{{st}^{-1}}(\eta'(st))\big)\\
&=\Del(t)\pr{\xi}{\eta}_l(rt,st)\\
&=\nal_t(\pr{\xi}{\eta}_l)(r,s)
\end{align*}
Hence $\nal_t$ is an extension of $\ga_t^l$, $\forall t\in G$. 
By Proposition \ref{prop:restr}, we have that $E_t^l=I_t$, 
$\forall t\in G$, from where we conlude that $\ga^l=\nal\r{I}$. 
Then 3. is proved.
\par It remains to show the continuity of $\ga$. We will prove that if 
$\xi_i$ is a net in $E$ that converges to $\xi$, with 
$\xi_i\in E_{t_i^{-1}}$, $\xi\in E_{t^{-1}}$, then  
$\ga_{t_i}(\xi_i)\f\ga_t(\xi)$. Let $\epsilon >0$. By the Cohen--Hewitt theorem, 
there exist $\xi'\in E_{t^{-1}}$, 
$a\in D_{t^{-1}}$, such that $\xi=\xi'a$, and we may suppose that $a\neq 0$. 
Since $E_{t^{-1}}^c$ is dense in  
$E_{t^{-1}}$, there exists $\eta\in E_{t^{-1}}^c$ such that $\norm{\xi'-\eta}
<\epsilon/4\norm{a}$. Therefore, 
\begin{equation}\label{e:xieta}
\norm{\xi-\eta a}\leq\epsilon/4
\end{equation}
Since the family $\{ D_r\}_{r\in G}$ is continuous, there exists 
$d\in C_c(G,A)$ such that $d(t^{-1})=a$. Let us define $\eta_s=\eta 
d(s^{-1})$. In particular, $\eta_t=\eta a$. Note that $\eta_s\in 
E_{t^{-1}}^c\bigcap E_{s^{-1}}^c$, because $\eta\in E_{t^{-1}}^c$, and  
$d(s^{-1})\in D_{s^{-1}}$, $\forall s\in G$. Since the action of $A$ on $E$ is  
continuous, and since $d$ is continuous, we have that $\eta_s\f\eta_t$ in $E$. 
Therefore, there exists $i_1$ such that $\forall i\geq i_1$ we have: 
\begin{equation}\label{e:xietas}
\norm{\xi -\eta_{t_i}}\leq\epsilon/4
\end{equation}
(In fact, we even have that $\eta_s\f\eta_t$ in the \ilt, because   
$\norm{\eta_s(r)-\eta_t(r)}
=\norm{\eta(r)d(s^{-1})-\eta(r)d(t^{-1})}
\leq\norm{\eta}_{\infty}\norm{d(s^{-1})-d(t^{-1})}$).
\par Let $\eta(r)=\big(r,\eta'(r)\big)$. Since $\eta_s=\eta d(s^{-1})$ we  
must have: \[\big(r,\eta_s'(r)\big)
=\big(r,\eta'(r)\big)\big(e,d(s^{-1})\big)
=\big(r,\al_r\big(\al_{r^{-1}}(\eta'(r)d(s^{-1}))\big)\big).\]
Then $\ga_s\eta_s(r)=\big(r,\ga_s'\eta_s'(r)\big)
=\big(r,\Del(s)^{1/2}\al_{rs}\big(
\al_{s^{-1}r^{-1}}(\eta'(rs)d(s^{-1}))\big)\big)$. 
Hence the map $G\times G\f A$ such that  
$(s,r)\ff\ga_s'\eta_s'(r)$ is continuous and has compact support. Therefore  
the map $G\f C_0(A)$ given by  
$s\ff \ga_s'\eta_s'$ is continuous. It follows that  
$\norm{\ga_s\eta_s-\ga_t\eta_t}_{\infty}
\f 0$ when $s\f t$. So we also have that 
$\norm{\ga_s\eta_s-\ga_t\eta_t}
\f 0$ if $s\f t$. Hence there exists $i_2$ such that $\forall i\geq i_2$:
\begin{equation}\label{e:etas}
\norm{\ga_{t_i}\eta_{t_i}-\ga_t\eta_t}\leq\epsilon/4
\end{equation}
On the other hand, since $\xi_i\f\xi$, there exists $i_3$ such that 
$\forall i\geq i_3$ we have:
\begin{equation}\label{e:xis}
\norm{\xi_i-\xi}<\epsilon/4
\end{equation}
Let $i_0\geq i_1,i_2,i_3$. Since each $\ga_s$ is an isometry, 
if $i\geq i_0$, by (\ref{e:xieta}), 
(\ref{e:xietas}), (\ref{e:etas}) and (\ref{e:xis}) we have that:
\begin{align*}
\norm{\ga_{t_i}(\xi_i)-\ga_t(\xi)}
&\leq\norm{\ga_{t_i}(\xi_i)-\ga_{t_i}(\eta_{t_i})}
 +\norm{\ga_{t_i}(\eta_{t_i})-\ga_t(\eta_t)}
 +\norm{\ga_t(\eta_t)-\ga_t(\xi)}\\
&\leq\norm{\xi_i-\eta_{t_i}}+\epsilon/4 +\norm{\eta a-\xi}\\
&\leq\norm{\xi_i-\xi}+\norm{\xi-\eta_{t_i}}+\epsilon/2\\
&<\epsilon
\end{align*}
This proves that $\ga$ is continuous, which ends the proof.
\end{proof}

\par Our next task is to show that, up to Morita equivalence, 
the \mea is unique. This will be proved in 
Proposition \ref{prop:super}. 

\begin{thm}\label{thm:equimork}
Let $\mb=(B_t)_{t\in G}$ be a \fb over $G$, $\me=(E_t)_{t\in G}$ a right ideal 
of $\mb$, and  
$\ma=(A_t)_{t\in G}$ a sub--\fb of $\mb$ contained in $\me$. Define  
$\nucc{\me}=\{ k\in\nucc{\mb}:\ k(r,s)\in E_{rs^{-1}},\, 
\forall r,s\in G\}$, and let $\nuc{\me}$ be the closure of $\nucc{\me}$ in  
$\nuc{\mb}$. Suppose that $\ma\me\subseteq\me$ and $\me\me^*\subseteq 
\ma$. Then:
\be
 \item $\nucc{\ma}\nucc{\me}\subseteq \nucc{\me}$, and 
        $\nuc{\ma}\nuc{\me}\subseteq \nuc{\me}$
 \item $\nucc{\me}\nucc{\mb}\subseteq \nucc{\me}$, and
       $\nuc{\me}\nuc{\mb}\subseteq \nuc{\me}$.
 \item $\nucc{\me}\nucc{\me}^*\subseteq\nucc{\ma}$, and
       $\nuc{\me}\nuc{\me}^*=\nuc{\ma}$. In particular, $\nuc{\ma}$ is  
       a hereditary sub-\cs of $\nuc{\mb}$.
 \item If $\gen(\me^*\me\bigcap B_t)$ is dense in $B_t$, for every $t\in G$, 
       then the closure of $\gen\nucc{\me}^*\nucc{\me}$ 
       in $\nucc{\mb}$ in the \ilt is $\nucc{\mb}$, and 
       $\ov{\gen}\,\nuc{\me}^*\nuc{\me}=\nuc{\mb}$. In particular,  
       $\nuc{\ma}$ and $\nuc{\mb}$ are Morita equivalent.
\ee 
Moreover, $\nuc{\me}$ is invariant under the natural action of $G$ on 
$\nuc{\mb}$. In the hipotheses of {\rm 4.} above, the natural actions of 
$G$ on $\nuc{\ma}$ and on $\nuc{\mb}$ are Morita equivalent. 
\end{thm}
\begin{proof}
As we have already remarked in Section \ref{sec:nucs}, if 
$\ma\subseteq\mb$, then $\nuc{\ma}\subseteq\nuc{\mb}$.  
\par Let $a\in\nucc{\ma}$, $\xi$, $\eta\in\nucc{\me}$, $b\in\nucc{\mb}$, 
$r$, $s\in G$. 
Then:\\ 
1. $a*\xi (r,s)=\int_Ga(r,t)\xi(t,s)dt\in A_{rt^{-1}}E_{ts^{-1}}
       \subseteq\me\bigcap B_{rs^{-1}}=E_{rs^{-1}}$. Therefore, 
        $\nucc{\ma}\nucc{\me}\subseteq \nucc{\me}$, and hence  
        $\nuc{\ma}\nuc{\me}\subseteq \nuc{\me}$.\\ 
2. Note that the second assertion is a consequence of the first one. On 
   the other hand: $\xi*b(r,s)=\int_G\xi(r,t)b(t,s)dt\in E_{rt^{-1}}B_{ts^{-1}}
       \subseteq\me\bigcap B_{rs^{-1}}=E_{rs^{-1}}$, which proves the first  
       assertion of 2. \\ 
3. $\xi*\eta^*(r,s)=\int_G\xi(r,t)\eta(s,t)^*dt\in E_{rt^{-1}}E_{ts^{-1}}
       \subseteq\ma\bigcap B_{rs^{-1}}$. Thus, 
       $\nucc{\me}\nucc{\me}^*\subseteq\nucc{\ma}$ e  
       $\nuc{\me}\nuc{\me}^*\subseteq \nuc{\ma}$. The equality follows from  
       the Cohen--Hewitt theorem and from the inclusion 
       $\nuc{\ma}\nuc{\ma}^*\subseteq\nuc{\me}\nuc{\me}^*.$\\ 
4. Let $F=\gen \nucc{\me}^*\nucc{\me}\subseteq\nucc{\mb}$, and consider the   
       dense subalgebra $\Theta=\bats{C_c(G)}{C_c(G)}$ of $C_c(G\times G)$.
       If $\varphi$, $\psi\in C_c(G)$, $\xi$, $\eta\in \nucc{\me}$, then: 
       \[(\varphi\otimes\psi)\xi^**\eta(r,s)
       =\int_G\varphi(r)\psi(s)\xi(t,r)^*\eta(t,s)dt
       =\int_G\big((\bar{\varphi}\xi)(t,r)\big)^*(\psi\eta)(t,s)dt%
       =(\bar{\varphi}\xi)^**(\psi\eta)(r,s),\]
       where $\varphi\xi (r,s)=\varphi (r)\xi(r,s)$. Since  
       $\varphi\xi\in\nucc{\me}$, $\forall \varphi\in C_c(G)$, 
       $\xi\in\nucc{\me}$, we see that $\Theta F\subseteq F$.
       \par Now, let $a\in B_{rs^{-1}}$ be such that $a=b^*c$, with  
       $b\in E_{tr^{-1}}$ and $c\in E_{ts^{-1}}$, for some $t\in G$. Let 
       $\xi$, $\eta\in\nucc{\me}=C_c(\me_{\nu})$ be such that $\xi(t,r)=b$, 
       $\eta(t,s)=c$, so $\xi(t,r)^*\eta(t,s)=b^*c=a$. 
       Let $\mv$ be a base of neighborhoods of the identity of $G$, and  
       $(f_V)_{V\in\mv}$ an approximate identity of $L^1(G)$ as in  
       \ref{lem:lema}. For each $V\in\mv$ consider 
       $\eta_V^t\in\nucc{\me}$ such that $\eta_V^t(r,s)=f_V(t^{-1}r)\eta(r,s)$.
       The map $G\f B_{rs^{-1}}$ given by  
       $u\ff\xi(u,r)^*\eta(u,s)$ is continuous with compact support. 
       Thus, if $\epsilon >0$, there  
       exists a neighborhood $V_{\epsilon}$ of the identity of $G$ such that,   
       if $t^{-1}u\in V_{\epsilon}$, then  
       $\norm{\xi(u,r)^*\eta(u,s)-\xi(t,r)^*\eta(t,s)}<\epsilon$. Hence, if  
       $V\subseteq V_{\epsilon}$: 
       \begin{align*}
       \norm{\xi^**\eta_V^t(r,s)-a}
       &=\norm{\int_Gf_V(t^{-1}u)[\xi(u,r)^*\eta(u,s)-\xi(t,r)^*\eta(t,s)]du}\\
       &\leq\int_{tV}f_V(t^{-1}u)
          \norm{\xi(u,r)^*\eta(u,s)-\xi(t,r)^*\eta(t,s)}du\\
       &<\int_Gf_V(t^{-1}u)\epsilon du\\
       &=\epsilon
       \end{align*}
       Since $B_{rs^{-1}}
       =\ov{\gen}\{E_{ur^{-1}}^*E_{us^{-1}}: u\in G\}$, 
       we conclude that  
       $\gen\nucc{\me}^*\nucc{\me}\bigcap B_{rs^{-1}}$ is dense in  
       $B_{rs^{-1}}$, $\forall r,s\in G$, and hence that 
       $\gen\nucc{\me}^*\nucc{\me}$ is dense in $\nucc{\mb}$ in the \ilt by 
       \ref{lem:bolu}. 
       This, together with 3., proves that $\nuc{\ma}$ and $\nuc{\mb}$ are  
       Morita equivalent.   
\par Let $\beta$ be the natural action of $G$ on $\nuc{\mb}$, and  
$k\in\nucc{\me}$, $t\in G$. Then $\beta_t(k)\r{(r,s)}=
\Del(t)k(rt,st)\in E_{rt(st)^{-1}}=E_{sr^{-1}}$, and therefore $k\in\nucc{\me}$. 
It follows that $\beta_t(\nucc{\me})=\nucc{\me}$, thus   
$\beta_t(\nuc{\me})=\nuc{\me}$. This completes the proof. 
\end{proof}

\begin{cor}\label{cor:equimork}
Let $\al$ and $\beta$ be Morita equivalent \pas with associated \fbs  
$\mb_{\al}$ and $\mb_{\beta}$ respectively. Let 
$\nal$ and $\nbeta$ be the natural actions on $\nuc{\mb_{\al}}$ and  
$\nuc{\mb_{\beta}}$. Then $\nal$ and $\nbeta$ are Morita equivalent. 
In particular, $\nuc{\mb_{\al}}$ and $\nuc{\mb_{\beta}}$ are Morita 
equivalent.  
\end{cor}
\begin{proof}
Suppose that $\ga$ is a \pa that implements the equivalence between $\al$ 
and $\beta$. Since the Morita equivalence of \pas is transitive, we may replace  
$\beta$ by the linking \pa of $\ga$ (see the proof of \ref{prop:equivprods} 
in page \pageref{vinc}). In the proof of Proposition  
\ref{prop:equivprods}, we constructed a right ideal $\me$ of $\mb_{\beta}$ such 
that $\mb_{\al}$, $\me$ and $\mb_{\beta}$ satisfy the conditions 1.--4. of 
Theorem \ref{thm:equimork} (with $\mb_{\al}$ in place of $\ma$ and $\mb_{\beta}$ 
instead of $\mb$). Now the result follows from \ref{thm:equimork}.     
\end{proof}

\begin{prop}\label{prop:super}
The \mea is unique, up to Morita equivalence.
\end{prop}
\begin{proof}
Since the Morita equivalence of \pas is transitive, by Corollary   
\ref{cor:equimork} above it is enough to show that if $(\beta,B)$ is an  
enveloping action of the \pa $(\al,A)$ and $\nal$ is the natural action on  
$\nuc{\ma}$, where $\ma$ is the associated Fell bundle of $\al$, then $\beta$  
and $\nal$ are Morita equivalent.
\par Now, let $\mb$ be the Fell bundle associated with the \ea $\beta$, 
and $\nbeta$ the natural action of $G$ on $\nuc{\mb}$. Consider the right ideal 
$\me$ of $\mb$ constructed in the proof of Theorem 
\ref{thm:equivmor}, that is, $\me =\{ (t,x)\in\mb:\ x\in A, \forall t\in G\}$. 
Applying Theorem \ref{thm:equimork}, it follows that  
$\nal$ and $\nbeta$ are Morita equivalent. Since $\mb$ is a saturated \fb. 
by \ref{prop:sat} we have that $\nuc{\mb}=I$, and therefore  
$\nbeta\sM\beta$ by \ref{thm:super} 2. and 3. In conclusion: $\nal\sM\beta$.   
\end{proof}

\begin{cor}\label{cor:pideal2}
Let $(\al,A)$ be a \pa with \mea $(\beta,B)$, and assume that $\mathcal{P}$ 
is an ideal property. Then $A$ has property  
$\mathcal{P}$ if and only if $B$ has property $\mathcal{P}$. 
In particular, $B$ is nuclear, liminal, antiliminal or postliminal, 
if and only if $A$ is respectively nuclear, liminal, antiliminal or 
postliminal.
\end{cor}
\begin{proof}
This is a direct consequence of \ref{cor:pideal}, \ref{thm:super} and 
\ref{prop:super}.
\end{proof}

\section{Partial actions induced on $\hat{A}$ and 
$\p(A)$}\label{sec:apps} 
Let $\al$ be a \pa on the \cs $A$ that has a Morita enveloping 
action $\beta$ acting on a \cs $B$. One can see that $\al$ induces \pas 
on $\hat{A}$ and $\p (A)$, the spectrum of $A$ and the primitive ideal 
space of $A$ respectively. In this part we will 
show that the enveloping actions of these \pas are precisely the actions 
induced by $\beta$ on $\hat{B}$ and $\p (B)$ respectively.
\par Let us fix some notation. If $I\id A$ let  
$\mo_I=\{ P\id A:\, \text{ $P$ is primitive and }P\not\supseteq I\}$, that 
is, $\{\mo_I:\ I\id A\}$ is the Jacobson topology of $\p (A)$. Recall that 
there is a natural map $\ka:\hat{A}\f\p(A)$, given by 
$\ka([\pi])=\ker\pi$, that is surjective but in general not injective.
The topology of $\hat{A}$ is the initial topology defined by $\ka$, so the 
open sets are of the form $\mv_I=\ka^{-1}(\mo_I)$. 
With this topology $\hat{A}$ is \lc, and it is compact if $A$ is unital 
(\cite{fd}, VII-6.11).
The restriction maps $r:\mo_I\f\p(I)$, such that $r(P)=P\cap I$, and 
$r:\mv_I\f\hat{I}$, such that $r([\pi])=[\pi\r{I}]$, are 
homeomorphisms. 

\begin{prop}\label{prop:induzidas}
Let $\beta:G\times B\f B$ be a continuous action of $G$ on the \cs $B$. 
Then:
\be
\item $\hat{\beta}:G\times \hat{B}\f\hat{B}$ such that $\hat{\beta}_t([\pi])=
      [\pi\circ\beta_{t^{-1}}]$ is a continuous action.
\item $\tilde{\beta}:G\times \p(B)\f\p(B)$ such that $\tilde{\beta}_t(P)=
      \beta_t(P)$ is a continuous action.  
\ee
\end{prop}
\begin{proof}
A proof of 1. may be found in \cite{rw}, 7.1. To prove 2., note that the 
following diagram commutes:
\[ \xymatrix
{G\times\hat{B}\ar@{->}[r]^-{id\times\ka}\ar[d]_-{\hat{\beta}}  
&G\times\p(B)\ar[d]^-{\tilde{\beta}}\\
\hat{A}\ar@{->}[r]_-{\ka} &\p(B)} \]
Since $\hat{B}$ has the initial topology induced by $\kappa$, 
then $\p(B)$ is precisely the topological quotient space 
of $\hat{B}$ with respect to $\ka$, and therefore $G\times\p(B)$ is 
the topological quotient space of $G\times\hat{B}$ with respect to 
$id\times\ka$. Thus, 
$\tilde{\beta}$ is continuous if and only if $\tilde{\beta}(id\times\ka)$ 
is continuous; but $\tilde{\beta}(id\times\ka)=\ka\hat{\beta}$, and  
$\ka\hat{\beta}$ is continuous.   
\end{proof}

It is possible to give a direct proof of the following result, but we will 
obtain it indirectly from \ref{prop:esprimen} and \ref{prop:ache}. 

\begin{prop}\label{prop:primap}
Let $\al=(\{ D_t\}_{t\in G}, {\al_t}_{t\in G})$ be a \pa $G$ on the \cs $A$. 
\be
 \item For each $t\in G$ let $\mo_t:=\mo_{D_t}
       =\{ P\in\p (A):\, P\not\supseteq D_t\}$. Then 
       $\tal =(\{\mo_t\}_{t\in G},\{\tal_t\}_{t\in G})$ is  
       a set theoretic \pa of $G$ on $\p (A)$, where $\tal_t(P)$ is the  
       unique primitive ideal of $A$ such that $\tal_t(P)\cap D_t=
       \al_t(P\cap D_{t^{-1}})$.
 \item For each $t\in G$, let  
       $\mv_t:=\{ [\pi]\in\hat{A}:\, \pi\r{D_t}\neq 0\}$, that is, 
       $\mv_t=\ka^{-1}(\mo_t)$, where $\mo_t$ is like in 1. Then 
       $\hat{\al}=(\{ \mv_t\}_{t\in G}, \{\hat{\al}_t\}_{t\in G})$ is a 
       set theoretic \pa of $G$ on $\hat{A}$, where, if $[\pi]\in\mv_t$, then  
       $\hat{\al}_t([\pi])$ is the class of the unique extension of 
       $\pi\circ\al_{t^{-1}}:D_t\f B(H_{\pi})$ to all of $A$.
\ee
\end{prop}
\par Suppose that the \pa $(\al,A)$ has an enveloping action 
$(\beta,B)$. Since $A\id B$, $\hat{A}$ is naturally homeomorphic to an 
open subset of $\hat{B}$, and therefore the action $\hat{\beta}$ restricted 
to this open set defines, via this homeomorphism, a \pa of $G$ on 
$\hat{\al}$. We show next that this \pa agrees with $\hat{\al}$. 
Similarly, the restriction of $\tilde{\beta}$ to $\p(A)$ agrees with $\tal$.

\begin{prop}\label{prop:esprimen}
Let $\al$ be a \pa on the \cs $A$, with enveloping action $\beta$ acting  
on the \cs $B$. Then $(\tal ,\p(A))$ and $(\hat{\al},\hat{A})$ are 
\pas, and: 
\be
 \item $(\tilde{\beta},\p(B))$ is the enveloping action of 
       $(\tal ,\p(A))$, and 
 \item $(\hat{\beta},\hat{B})$ is the enveloping action of 
       $(\hat{\al},\hat{A})$. 
\ee
\end{prop}
\begin{proof}  
We identify $\p(A)$ with $\mo_A$ through the homeomorphism $r^{-1}$. 
Thus $\tal=(\{\mo_t\},\{\tal_t\})$ becomes: 
$\mo_t=\{ P\in\p(B):\, P\not\supseteq D_t\}$, and if $P\in\mo_{t^{-1}}$, 
then $\tal_t(P)\in\mo_t$ is the unique primitive ideal of $B$ such that 
$\tal_t(P)\cap D_t=\al_t(P\cap D_{t^{-1}})$. 
Let us see first that $\dom(\tilde{\beta}_t\r{\mo_A})=\dom\tal_t$: 
\[ \dom(\tilde{\beta}_t)\r{\mo_A}
=\{ P\in\p(B):\, P\in\mo_A\text{ e }\tilde{\beta}_t(P)\in\mo_A\}
=\mo_A\cap\mo_{\tilde{\beta}_{t^{-1}}(A)} 
=\mo_{A\cap\tilde{\beta}_{t^{-1}}(A)} 
=\dom({\tal_t}).\]
Now, if $P\in\mo_{t^{-1}}$, we have that $\tilde{\beta}_t(P)=\beta_t(P)
\not\supseteq D_t$, and hence $\tilde{\beta}_t(P)\in\mo_t$. But  
$\beta_t(P)\cap D_t=\beta_t(P\cap\beta_{t^{-1}}(D_t))=
\al_t(P\cap D_{t^{-1}})$. Then $\tilde{\beta}_t(P)$ agrees with 
$\tal_t(P)$; in particular $\tal$ is a \pa. 
\par It remains to verify that the $\tilde{\beta}$-orbit of $\mo_A$ is 
$\p(B)$. Suppose that there exists $P\in\p(B)$ such that 
$P\neq\tilde{\beta}_t(Q)$, $\forall Q\in\mo_A$. Then 
$\tilde{\beta}_t(P)\notin\mo_A$, $\forall t\in G$, that is, 
$\beta_t(P)\supseteq A$, $\forall t\in G$.
But then $P\supseteq\beta_t(A)$, $\forall t\in G$, and therefore 
$P\supseteq \ov{[\beta (A)]}=B$, because $\beta$ is the enveloping action  
of $\al$. The contradiction implies that every primitive ideal belongs to  
the $\tilde{\beta}$-orbit of some element of $\mo_A$. 
\par As for $\hat{\beta}$ and $\hat{\al}$, identify $\hat{A}$ with  
$\mv_A$ via $r^{-1}$. Then 
$\hat{\al}=(\{\mv_t\} ,\{\hat{\al}_t\} )$ becomes: 
$\mv_t=\{ [\pi]\in\hat{B}:\, \pi\r{D_t}\neq 0\}$, and for  
$[\pi]\in\mv_{t^{-1}}$, $\hat{\al}_t([\pi])$ is the class of the unique  
extension to $B$ of the irreducible \rep $\pi\circ\al_{t^{-1}}$ of $D_t$. 
From the computations above it follows that   
$\dom(\hat{\beta})\r{\mv_A}=\mv_{t^{-1}}=\dom(\hat{\al}_t)$. On the other 
hand, if $[\pi]\in\mv_{t^{-1}}$, then $(\pi\circ\beta_{t^{-1}})$ is  
an extension to $B$ of the \rep $\pi\circ\al_{t^{-1}}$, and therefore 
$\hat{\beta}_t([\pi])$ agrees with $\hat{\al}_t([\pi])$. It remains to show 
that the $\hat{\beta}$-orbit of $\mv_A$ is all of $\hat{B}$, and this is 
similar to which has been done previously for $\tal$ and $\tilde{\beta}$: if 
$[\pi]\notin\hat{\beta}_t(\mv_A)$, $\forall t\in G$, then  
$\pi\r{\beta_t(A)}=0$, $\forall t\in G$, and therefore $\pi=0$, what  
is a contradiction. 
\end{proof}

\begin{lem}\label{lem:ache}
Let $\ga =(\{ E_t\}_{t\in G},\{ \al_t\}_{t\in G})$ be a \pa of 
$G$ on a positive \ct $E$, and let $(\beta,B)=(\ga^l,E^l)$, 
$(\al,A)=(\ga^r,E^r)$. Consider the set theoretic \pas $(\tilde{\beta},\p(B))$ 
and $(\tal,\p(A))$ induced by $\beta$ and $\al$ on $\p(B)$ and $\p(A)$ 
respectively (see Proposition \ref{prop:primap}). Then the Rieffel 
homeomorphism $\rief:\p(B)\f\p(A)$ is an isomorphism of set theoretic \pas 
$\rief:\tilde{\beta}\f\tal$. A similar statement stands  
for the set theoretic \pas $\hat{\beta}$ and $\hat{\al}$ induced by $\beta$  
and $\al$ on the corresponding spectra $\hat{B}$ and $\hat{A}$ of $B$ and $A$. 
\end{lem}
\begin{proof}
It is well known that if $P\in\p(B)$, then $\rief(P)\in\p(A)$, and also that 
$\rief$ is a homeomorphism (see for instance Corollary 3.33 of \cite{rw}).  
Since $\rief(E_t^l)=E_t^r$, it follows that $\rief(\mo_t^{B})=\mo_t^{A}$. 
Let $I\id B$. Then $\beta_t\big(E_{t^{-1}}^l\cap I\big)$ is an ideal  
of $E_t^l$, and therefore there exists a unique $F\id E_t$ such that  
$F^l=\beta_t\big(E_{t^{-1}}^l\cap I\big)$. We claim that  
$F=\ga_t(E_{t^{-1}}\cap IE)$. Indeed, $\pr{F}{F}_l=
\pr{\ga_t(E_{t^{-1}}\cap IE)}{\ga_t(E_{t^{-1}}\cap IE)}_l=
\beta_t(\pr{E_{t^{-1}}\cap IE}{E_{t^{-1}}\cap IE}_l)$, so $F^l=
\beta_t(E_{t^{-1}}^l\cap I)$. Similarly, if $J\id A$, 
we have that $\al_t\big(E_{t^{-1}}^r\cap J\big)
=\big(\ga_t(E_{t^{-1}}\cap EJ)\big)^r$.
Now, if $P\in\mo_{t^{-1}}^B$, then $\tilde{\beta}_t(P)=Q$, where  
$Q\in\p(B)$ is the unique primitive ideal such that  
$Q\cap E_t^l=\beta_t\big(E_{t^{-1}}^l\cap P\big)$.  
Since $\rief$ is a lattice isomorphism,
it follows that $\rief(Q)\cap \rief(E_t^l)
=\rief\big(\beta_t(E_{t^{-1}}^l\cap P\big).$ That is, 
\[ \rief(Q)\cap E_t^r
=\big(E_{t^{-1}}\cap PE\big)^r
=\al_t\big(E_{t^{-1}}^r\cap (PE)^r\big)
=\al_t\big(E_{t^{-1}}^r\cap \rief(P)\big).\]
Thus, since $\rief(Q)$ is a primitive ideal, it must agree with  
$\tal_t\big(\rief(P)\big)$. It follows that $\rief:\tilde{\beta}\f\tal$ 
is a morphism of set theoretic \pas. Similarly, its inverse map 
$\rief:\tal\f\tilde{\beta}$ is a morphism, so $\rief$ is an 
isomorphism between the set theoretic \pas $\tilde{\beta}$ and $\tal$. 
\par The proof of the corresponding statement for $\hat{\beta}$ and 
$\hat{\al}$ is similar and it is left to the reader.
\end{proof}
\begin{prop}\label{prop:ache}
Let $(\al,A)$ be a \pa of $G$ on the \cs $A$, and let $(\beta,B)$ be its 
Morita enveloping action. Then $(\tal ,\p(A))$ and $(\hat{\al},\hat{A})$ are 
\pas, and:
\be
 \item $(\tilde{\beta},\p(B))$ is the enveloping action of 
       $(\tal ,\p(A))$, and 
 \item $(\hat{\beta},\hat{B})$ is the enveloping action of 
       $(\hat{\al},\hat{A})$. 
\ee
\end{prop}
\begin{proof}
Since $\beta$ is the \mea of $\al$, this one is Morita equivalent to 
$\beta\r{I}$, for some ideal $I$ of $B$. By \ref{prop:esprimen}, 
$\beta\r{I}$ is a \pa on 
$\p(I)$. On the other hand, $\rief$ is a homeomorphism, and therefore 
$\tal$ is continuous, that is, $\tilde{\al}$ is a \pa on $\p(A)$. Since 
$(\tilde{\beta},\p(B))$ is the enveloping action of 
$(\tilde{\beta\r{I}},\p(I))$ and $\rief$ is an isomorphism of \pas 
between $(\tilde{\beta},\p(B))$ and $(\tal ,\p(A))$, it follows that  
$(\tilde{\beta},\p(B))$ is the enveloping action of $(\tal ,\p(A))$. 
The proof of 2. is similar.
\end{proof}

\begin{cor}\label{cor:iguais1}  
If $(\al,A)$ is a \pa with \mea $(\beta,B)$, 
then: \[\hat{A}=\p(A)\iff\hat{B}=\p(B).\]
\end{cor}
\begin{proof}
By \ref{prop:ache} we have that  
$(\hat{\beta},\hat{B})$ is the enveloping action of $\hat{\al}$ and  
that $(\tilde{\beta},\p(B))$ is the enveloping action of $\tal$. 
Suppose that $\hat{\al}=\tal$. Since the enveloping action is  
unique, we have that $(\hat{\beta},\hat{B})=(\tilde{\beta},\p(B))$. 
The converse is clear.  
\end{proof}

The following result may be thought of as a non--commutative version of the  
well known fact that the integral curves of a vector field on a 
compact manifold $X$ are defined on all of $\R$.

\begin{cor}\label{cor:conti2}
Suppose that $\al$ is a \pa of $G$ on the \cs $A$, with \mea  
$(\beta,B)$. If $\p(A)$ is compact (this is true if $A$ is unital), then 
there exists an open subgroup $H$ of $G$ such that the restriction of 
$\al$ to $H$ is a global action. In particular, if $G$ is a connected group, 
then $\al$ is a global action. 
\end{cor}
\begin{proof}
If $A$ is unital, then $\p(A)$ is compact by \cite{fd}, VII-6.11.  
Now, if $\p(A)$ is compact, by \ref{prop:conti} there exists an open subgroup 
$H$ of $G$ such that $\tal$ restricted to $H$ is a global action on 
$\p(A)$, and therefore every primitive ideal of $A$ is in the domain of 
$\tal_s$, $\forall s\in H$. By the definition of $\tal$, this implies that 
there is no primitive ideal of $A$ containing the ideal $D_{s^{-1}}$, 
and hence $D_{s^{-1}}=A$, $\forall s\in H$. Therefore $\al\r{H\times A}$ 
is a global action. If $G$ is connected, then $G=H$.
\end{proof}

\section{Takai duality for partial actions}\label{sec:duality}
\par In this last section of the paper we relate our previous results
on \eas with Takai duality for \pas. If we tried to translate naively 
Takai duality from the context of global actions to our case of partial ones, 
it should be expressed as follows: if $\al$ is a \pa of $G$ on $A$ and $\del$ 
is the dual coaction of $G$ on $A\sd{\al,r}G$, then $A$ and   
$A\sd{\al,r}G\sd{\del}\hat{G}$ are Morita equivalent. However, as pointed out 
by Quigg in \cite{dualpaut}, this version of Takai duality fails for \pas, and 
its failure is proportional to the ``partialness'' of $\al$. 
\par In what follows we show that we still have a Takai duality for \pas, 
which may be briefly expressed in this manner: 
if $\hat{\del}$  is the dual action of $G$ on 
$A\sd{\al,r}G\sd{\del}\hat{G}$, then 
$\big(\hat{\del},A\sd{\al,r}G\sd{\del}\hat{G}\big)$ is the \mea of $(\al,A)$.
In other words, $A$ is no longer Morita equivalent to 
$A\sd{\al,r}G\sd{\del}\hat{G}$, but to an ideal $I$ of this double crossed 
product, in such a way that this ideal together with the dual action allow us 
to recover the whole double crossed product. 
\par To be more precise, we will prove that if $\mb_{\al}$ is the \fb 
associated 
to the \pa $\al$, and $\nal$ is the natural action on $\nuc{\mb_{\al}}$, then 
the dynamical systems $\big(\nal,\nuc{\mb_{\al}}\big)$ and 
$\big(\hat{\del},A\sd{\al,r}G\sd{\del}\hat{G}\big)$ are isomorphic.    

\par Before proceeding, let us recall some basic facts about crossed products 
by coactions; for more details, the reader is referred to \cite{kattak} and 
\cite{lprs}. The tensor product we consider is the minimal tensor product. 
If $A$ and $B$ are \css, set   
$\tilde{M}(\bts{A}{}{B}):=\{ m\in M(\bts{A}{}{B}):\, 
m(1\otimes B)+(1\otimes B)m 
\subseteq \bts{A}{}{B}\}$. Let $w_G:L^2(G\times G)\f L^2(G\times G)$ be the 
unitary operator such that $w_G(\xi)\r{(r,s)}=\xi(r,r^{-1}s)$. 
Its adjoint is $w_G^*(\xi)\r{(r,s)}=\xi(r,rs)$. The \textit{comultiplication} 
on $C^*_r(G)$ is $\del_G:C^*_r(G)\f 
\tilde{M}\big(\bts{C^*_r(G)}{}{C^*_r(G)}\big)$ given by $\del_G(x)
=w_G(x\otimes 1)w_G^*$, $\forall x\in C^*_r(G)$.
A \textit{coaction} of the \lc group $G$ on a \cs $A$ is an injective \hm 
$\del:A\f\tilde{M}\big(\bts{A}{}{C^*_r(G)}\big)$ such that for any  
approximate unit $(e_i)_{i\in I}$ of $A$ we have that $\del(e_i)$ converges 
strictly to $1$ in $\tilde{M}\big(\bts{A}{}{C^*_r(G)}\big)$, and such that  
$(\del\otimes id)\del =(id\otimes\del_G)\del.$  
In particular, $\del_G$ is a coaction of $G$ on $C^*_r(G)$. 
\par Suppose that $\del :A\f\tilde{M}\big(\bts{A}{}{C^*_r(G)}\big)$ 
is a coaction of $G$ on $A$. The \textit{reduced crossed product} 
$A\sd{\del,r}\hat{G}$ of $A$ by the coaction $\del$ is:
\[ A\sd{\del,r}\hat{G}:=\big\{\del(a)(1\otimes \varphi):\, a\in A,\, 
\varphi\in C_0(G)\big\}\subseteq M\big(\bts{A}{}{\mathcal{K}}\big),\]
where $\varphi\in C_0(G)$ acts by multiplication on $L^2(G)$. 

\par If $\del$ is a coaction of $G$ on $A$, then there exists a canonical action 
$\hat{\del}$ of $G$ on $A\sd{\del,r}\hat{G}$. This action is called  
the \textit{dual action} of $G$ on $A\sd{\del,r}\hat{G}$, and  
it is characterized by the fact that it verifies
$\hat{\del}_t\big(\del(a)(1\otimes\varphi)\big)=\del(a)(1\otimes \varphi_t)$
where, if $\varphi\in C_0(G)$, then $\varphi_t(s)=\varphi(st)$, 
$\forall s\in G$. 
\par Now consider a \fb $\mb=(B_t)_{t\in G}$, and let $\mb_G$ be 
the trivial \fb over $G$, that is, the \fb associated with the trivial action 
of $G$ on $\C$. We have a \fb $\bts{\mb}{}{\mb_G}$ over $G\times G$, 
whose reduced \cs is naturally isomorphic to $\bts{C^*_r(\mb)}{}{C^*_r(G)}$. 
We also have that $L^2(\bts{\mb}{}{\mb_G})\cong\bts{L^2(\mb)}{}{L^2(G)}$ (see 
\cite{fav0} or \cite{fav3}). Let  
$w_{\mb}:L^2(\bts{\mb}{}{\mb_G})\f L^2(\bts{\mb}{}{\mb_G})$ be such that 
$w_{\mb}(\xi)\r{(r,s)}=\xi(r,r^{-1}s)$. Then $w_{\mb}$ is an adjointable 
operator, in fact unitary, with adjoint $w_{\mb}^*$ such that 
$w_{\mb}^*(\xi)\r{(r,s)}=\xi(r,rs)$. Now consider the map  
$\mathcal{L}\big(L^2(\mb)\big)\f\mathcal{L}\big(L^2(\bts{\mb}{}{\mb_G})\big)$ 
such that $x\ff w_{\mb}(x\otimes 1)w_{\mb}^*$. Let 
$\del_{\mb}:C^*_r(\mb)\f\tilde{M}\big(\bts{C^*_r(\mb)}{}{C^*_r(G)}\big)$ be 
the restriction of this map to $C^*_r(\mb)$. Then $\del_{\mb}$ is a coaction 
of $G$ on $C^*_r(\mb)$ (\cite{exeng}), called the \textit{dual coaction} of 
$G$ on $C^*_r(\mb)$.

\par Let $f\in C_c(\mb)$, $\xi\in C_c(\mb)\subseteq L^2(\mb)$, and 
$x\in C_c(G)\subseteq L^2(G)$. Then we have:
\begin{align*}
\del_{\mb}(\La_f)(\xi\otimes x)\r{(r,s)}
&=w_{\mb}(\La_f\otimes 1)w_{\mb}^*(\xi\otimes x)\r{(r,s)}\\
&=(\La_f\otimes 1)w_{\mb}^*(\xi\otimes x)\r{(r,r^{-1}s)}\\
&=\int_G(\La_{f(t)}\otimes 1)w_{\mb}^*(\xi\otimes x)
\r{(t^{-1}r,r^{-1}s)}dt\\
&=\int_G(\La_{f(t)}\otimes 1)(\xi\otimes x)\r{(t^{-1}r,t^{-1}s)}dt\\
&=\int_Gf(t)\xi(t^{-1}r)\otimes x(t^{-1}s)dt\\
&=\int_G(\La_{f(t)}\otimes\la_t)(\xi\otimes x)\r{(r,s)}dt.
\end{align*}
Consider the \rep of $\mb$ given by $b_t\f\La_{b_t}\otimes\la_t\in 
\mathcal{L}\big(\bts{L^2(\mb)}{}{L^2(G)}\big)$. It follows from the 
computations above that the integrated form of this \rep factors through 
$\del$. Thus if $\pi:C^*_r(\mb)\f B(H)$ is a faithful \rep, then 
$\del_{\mb}(\La_f)=\pi_{\la}$, where $\pi_{\la}$ is the integrated \rep of 
$\mb\f B\big(L^2(G,H)\big)$ such that $b_t\ff\la_t\otimes\pi(b_t)$.
 
\par Let $\al$ be a \pa, and $\mb_{\al}$ its 
associated \fb. The next result shows the form that Takai duality has for 
\pas: $(\beta ,\nuc{\mb})$ and  
$(\hat{\del}_{\mb}$, $C^*_r(\mb)\sd{\del_{\mb},r}\hat{G})$ 
are isomorphic dynamical systems. 

\begin{prop}\label{prop:nadual}
Let $\mb=(B_t)_{t\in G}$ be a Fell bundle over $G$, $\del$ the  
dual coaction of $G$ on $C^*_r(\mb)$, $\hat{\del}$ the dual action of $G$ 
on $C^*_r(\mb)\sd{\del,r}\hat{G}$, and $\beta$ the natural action of $G$ on 
$\nuc{\mb}$. Then there exists an isomorphism $\sigma:\nuc{\mb}\f
C^*_r(\mb)\sd{\del,r}\hat{G}$ such that  
$\sigma\beta_t=\hat{\del}_t\sigma$, $\forall t\in G$. 
\end{prop} 
\begin{proof}
We may assume, without loss of generality, that $C^*_r(\mb)\subseteq B(H)$ 
non--degenerately, for some Hilbert space $H$. Therefore also  
$\mb\subseteq B(H)$, and $C^*_r(\mb)\sd{\del,r}\hat{G}\subseteq 
B\big(L^2(G,H)\big)$. 
\par On the other hand, by \ref{prop:pic}, the inclusion $\mb\inc B(H)$ defines 
a faithful \rep $\sigma:\nuc{\mb}\f B\big(L^2(G,H)\big)$ such that, if 
$k\in\nucc{\mb}$, $x\in L^2(G,H)$, then  
$\sigma(k)\r{r}=\int_Gk(r,s)x(s)ds$. 
\par Now, if $f\in C_c(\mb)$, $\varphi\in C_c(G)\subseteq C_0(G)$, and
$\xi\in C_c(G)\subseteq L^2(G)$, $h\in H$, we have:
\begin{align*}
\del(f)(1\otimes\varphi)(\xi\otimes h)\r{r}
&=\int_G\la_s(\varphi\xi)(r)f(s)hds\\
&=\int_G\varphi(s^{-1}r)\xi(s^{-1}r)f(s)hds\\
&=\int_G\varphi(t^{-1})\xi(t^{-1})f(rt)hdt
  \ \ \\ 
&=\int_G\Del(s)^{-1}\varphi(s)\xi(s)f(rs^{-1})hds
  \ \ \\ 
&=\sigma(k_{\varphi,f})(\xi\otimes h)\r{r},
\end{align*}
where $k_{\varphi,f}(r,s)=\Del(s)^{-1}\varphi(s)f(rs^{-1})$. Thus, 
$\sigma(k_{\varphi,f})=\del(f)(1\otimes \varphi)$, and hence  
$\sigma\big(\nuc{\mb}\big)\supseteq C^*_r(\mb)\sd{\del,r}\hat{G}$. 
To prove the converse inclusion, it is enough to show that
$F:=\gen\{ k_{\varphi,f}:\,\varphi\in C_c(G),\, f\in C_c(\mb)\}$ is  
dense in $\nucc{\mb}$ in the \ilt. This will follow from Lemma \ref{lem:bolu}. 
For given $\phi$, $\psi\in C_c(G)$, let $\phi\diamond\psi:G\times G\f \C$ be 
such that $\phi\diamond\psi(r,s)=\phi(s)\psi(rs^{-1})$. It is clear that 
$\phi\diamond\psi$ is continuous and has compact support. Let $\Theta:=
\gen\{\psi_1\diamond\psi_2:\, 
\psi_1,\psi_2\in C_c(G)\}$. Then $\Theta$ is dense in $C_0(G\times G)$ by  
the Stone--Weierstrass theorem. In particular, $\Theta$ is a dense subspace 
$C_c(G\times G)$ in the \ilt. Let us see that $\Theta F\subseteq F$:
if $\psi_1\diamond\psi_2\in\Theta$, $k_{\varphi,f}\in F$, we have:
\[ (\psi_1\diamond\psi_2)k_{\varphi,f}(r,s)
=\Del(s)^{-1}\psi_1(s)\psi_2(rs^{-1})\varphi(s)f(rs^{-1})
=k_{\psi_1\varphi,\psi_2f}(r,s).\]  
Now, it is clear that $F(r,s)=B_{rs^{-1}}$, $\forall r,s\in G$. 
This shows that $F$ is dense in $\nuc{\mb}$, and therefore  
$\sigma :\nuc{\mb}\f C^*_r(\mb)\sd{\del,r}\hat{G}$ is an isomorphism.
Finally, if $t\in G$ we have:
\[\sigma^{-1}\big(\hat{\del}_t(\del(f)(1\otimes\varphi)\big)
=k_{\varphi_t,f}.\]
On the other hand:
\[
\beta_t(k_{\varphi,f})\r{(r,s)}
=\Del(t)k_{\varphi,f}(rt,st)
=\Del(t)\Del(st)^{-1}\varphi(st)f(rs^{-1})
=\Del(s)^{-1}\varphi_t(s)f(rs^{-1})
=k_{\varphi_t,f}\r{(r,s)},\]
so $\sigma\beta_t(k_{\varphi,f})=\hat{\del}_t\sigma(k_{\varphi,f})$, 
and since $F$ is dense in $\nuc{\mb}$, it follows that
$\sigma\beta_t=\hat{\del}_t\sigma$.
\end{proof}
\begin{rk}\label{rk:pideal}
Since $\nuc{\mb}$ is isomorphic to $C^*_r(\mb)\sd{\del,r}\hat{G}$, we may 
apply Corollary \ref{cor:pideal2} (or \ref{cor:pideal}) 
to $C^*_r(\mb)\sd{\del,r}\hat{G}$. 
In particular, we have that  
$C^*_r(\mb)\sd{\del,r}\hat{G}$ is nuclear, liminal, antiliminal or 
postliminal, 
if and only if $B_e$ is respectively nuclear, liminal, antiliminal or 
postliminal. This was already known for discrete groups 
(\cite{dcoq},\cite{dcong}). 
\par We also conclude that if $\ma$ and $\mb$ satisfy conditions 1.--3. 
in Theorem \ref{thm:equimork}, then $C^*_r(\ma)\sd{\del,r}\hat{G}$ is an  
hereditary sub-\cs of $C^*_r(\mb)\sd{\del,r}\hat{G}$. If they also satisfy   
condition 4. in \ref{thm:equimork}, then 
$\big(C^*_r(\ma)\sd{\del,r}\hat{G},\hat{\del}\big)$ and  
$\big(C^*_r(\mb)\sd{\del,r}\hat{G},\hat{\del}\big)$ are Morita equivalent 
dynamical systems. In particular, if $(\al,A)$ and $(\beta,B)$ are Morita   
equivalent \pas, then  
$\big(A\sd{\al,r}G\sd{\del,r}\hat{G},\hat{\del}\big)$ and  
$\big(B\sd{\beta,r}G\sd{\del,r}\hat{G},\hat{\del}\big)$ 
are Morita equivalent.
\end{rk}

\vspace*{1cm}
\end{document}